\documentclass[a4paper,12pt]{report}
\usepackage[british]{babel}
\usepackage{latexsym,amsfonts,amsmath,amssymb,amsthm}
\usepackage[latin1]{inputenc}
\usepackage[T1]{fontenc}
\usepackage{dsfont}
\usepackage{makeidx}
\usepackage[pdfpagelabels,bookmarksnumbered]{hyperref}
\hypersetup{
pdfauthor = {Hans Adler},
pdftitle = {Explanation of Independence},
pdfsubject = {Stability theory (MSC 03C45)},
pdfkeywords = {stability theory, model theory, independence, thorn-forking}}

\newtheorem{theorem}{Theorem}[chapter]
\newtheorem{corollary}[theorem]{Corollary}
\newtheorem{lemma}[theorem]{Lemma}
\newtheorem{remark}[theorem]{Remark}
\newtheorem{fact}[theorem]{Fact}
\newtheorem{proposition}[theorem]{Proposition}
\newtheorem{Definition}[theorem]{Definition}
\newenvironment{definition}{\begin{Definition}}{\end{Definition}}

\newtheorem{Example}[theorem]{Example}
	\newenvironment{example}{\begin{Example}\normalfont}{\end{Example}}
\newtheorem{Exercise}[theorem]{Exercise}
	\newenvironment{exercise}{\begin{Exercise}\normalfont}{\end{Exercise}}
\newcommand{\sir}{strict independence relation}
\newenvironment{exercises}{\subsection*{Exercises}\small}{\par}
\newenvironment{notes}{\subsection*{Notes}\small}{\par}
\newtheorem{question}{Question}[chapter]

\renewcommand{\phi}{\varphi}
\renewcommand{\theta}{\vartheta}

\newcommand{\card}[1]{\left|#1\right|}
\newcommand{\restrict}{{\upharpoonright}}
\newcommand{\concat}{^\smallfrown}

\newcommand{\Q}{\mathds Q}
\newcommand{\R}{\mathds R}
\newcommand{\dist}{\operatorname{d}}

\def\Ind#1#2{#1\setbox0=\hbox{$#1x$}\kern\wd0\hbox to 0pt{\hss$#1\mid$\hss}\lower.9\ht0\hbox to 0pt{\hss$#1\smile$\hss}\kern\wd0}
\def\Notind#1#2{#1\setbox0=\hbox{$#1x$}\kern\wd0\hbox to 0pt{\mathchardef
	\nn=12854\hss$#1\nn$\kern1.4\wd0\hss}\hbox to
	0pt{\hss$#1\mid$\hss}\lower.9\ht0 \hbox to
	0pt{\hss$#1\smile$\hss}\kern\wd0}
\newcommand{\ind}[1][]{\mathop{\mathpalette\Ind{}^{\!\!\!\!\rlap{$\scriptscriptstyle\textnormal{#1}$}\,\,\,\,}}}
\newcommand{\nind}[1][]{\mathop{\mathpalette\Notind{}^{\!\!\!\rlap{$\scriptscriptstyle\textnormal{#1}$}\,\,\,}}}
\newcommand{\indb}[1][]{\mathop{\mathpalette\Ind{}^{\!\!\!\!\rlap{#1}\,\,\,\,}}}
\newcommand{\nindb}[1][]{\mathop{\mathpalette\Notind{}^{\!\!\!\rlap{#1}\,\,\,}}}
\def\thind{\ind[\th]}       \def\nthind{\nind[\th]}
\def\dind{\ind[d]}          \def\ndind{\nind[d]}
\def\find{\ind[f]}          \def\nfind{\nind[f]}
\def\aind{\ind[a]}          
\def\eqind{\ind[eq]}        
\def\mind{\ind[M]}          \def\nmind{\nind[M]}
\def\pind{\indb[$\scriptstyle\prime$]}  \def\npind{\nindb[$\scriptstyle\prime$]}
\def\sind{\indb[$\scriptstyle*$]}       \def\nsind{\nindb[$\scriptstyle*$]}
\def\omind{\ind[$\Omega$]}  \def\nomind{\nind[$\Omega$]}
\def\omsind{\indb[${\scriptscriptstyle\Omega}{\scriptstyle*}$]}\def\nomsind{\nindb[${\scriptscriptstyle\Omega}{\scriptstyle*}$]}

\def\xind{\ind[$\Xi$]}      

\def\equ{\mathop{\equiv}}

\newcommand{\acl}{\operatorname{acl}}
\newcommand{\dcl}{\operatorname{dcl}}
\newcommand{\bdd}{\operatorname{bdd}}
\newcommand{\tp}{\operatorname{tp}}
\newcommand{\stp}{\operatorname{stp}}

\renewcommand{\S}{\operatorname{S}}
\newcommand{\cl}{\operatorname{cl}}

\newcommand{\typedef}d
\newcommand{\eq}{^{\operatorname{eq}}}

\newcommand{\heq}{^{\operatorname{heq}}}
\newcommand{\cb}{\operatorname{cb}}
\newcommand{\wcb}{\operatorname{wcb}}
\newcommand{\aker}{\operatorname{aker}}
\newcommand{\monster}{\mathcal M}

\newcommand{\D}{\operatorname{D}}

\def\Xim{\Xi_{\textnormal{M}}}

\makeindex
\providecommand*{\phantomsection}{}	

\author{\href{mailto:hans.adler@math.uni-freiburg.de}{Hans Adler}}

\begin{document}
	\begin{titlepage}
	\begin{center}
	\ \\[2cm]
 	{\LARGE\bf EXPLANATION\\OF\\INDEPENDENCE\\[5mm]}
	\ \\[2cm]
	{\LARGE 
	{\sc Hans Adler}}\vfill
	Dissertation zur Erlangung des Doktorgrades\\
	der Fakult\"at f\"ur Mathematik und Physik\\
	der Albert-Ludwigs-Universit\"at Freiburg im Breisgau\\
	June 2005
	\end{center}
	\end{titlepage}
	\pagestyle{plain}
	\begin{tabular}{ll}
		Dekan:&Prof.{} Dr.{} Josef Honerkamp\\
		Gutachter:&Prof.{} Dr.{} Martin Ziegler\\
		&Prof.{} Dr.{} Frank Olaf Wagner\\
		Datum der m\"undl.{} Pr\"ufung:&16.~November 2005\thispagestyle{empty}
	\end{tabular}\vfill\newpage
	\pagenumbering{roman}\setcounter{page}{1}
	\chapter*{Preface}
	The original inspiration and motivation for this thesis came
	from~\cite{Hans Scheuermann: Unabhaengigkeitsrelationen}---an inexperienced author's
	attempt to write down in a short period of time
	everything he knew and many things he did
	not understand, yet containing a few new ideas.
	With all due respect for the author of~\cite{Hans Scheuermann: Unabhaengigkeitsrelationen},
	I always felt that I could do better than that.
	Over many years of work---some of it on mathematics, but much more
	on a paid job as a software developper---the original inspiration
	together with input from more recent research by Alf Onshuus
	and Itay Ben-Yaacov was gradually transformed into the present
	dissertation.
	\tableofcontents\newpage\pagenumbering{arabic}
\chapter*{Introduction}
\addcontentsline{toc}{chapter}{Introduction}

Independence relations\footnote{An independence relation is the same thing as
a notion of independence as defined by Kim and Pillay.}
generalise concepts such as linear independence in
vector spaces or transcendence in fields to much more general complete
theories. (In an even more general context than covered in this thesis, they also
generalise orthogonality in Hilbert spaces and stochastic independence.)

More complicated notions such as orthogonality of types are built on top
of independence. Also, important dividing lines between theories are often
defined by means of additional axioms which an independence relation may
or may not satisfy. Therefore the notion of independence relation is certainly
among the most fundamental notions in stability theory.
Independence relations satisfying certain additional properties were studied
in~\cite{Michael Makkai: A survey of basic stability theory}, \cite{Harnik + Harrington: Fundamentals of Forking}
and~\cite{Kim + Pillay: Simple theories}.
An independence relation for o-minimal theories is more or less explicit
in~\cite{Pillay + Steinhorn: Definable sets in ordered structures I} and
in~\cite{Anand Pillay: Some remarks on definable equivalence relations in O-minimal structures}.
Recently all important independence relations were unified
by work of Thomas Scanlon, Alf Onshuus and Clifton Ealy.

In this thesis I study independence relations in a systematic way.
For the purposes of this introduction let us say that a relation
$\ind$ is a pre-independence relation if $\ind$ satisfies the first five axioms
for independence relations
as well as strong finite character (cf.~Definitions~\ref{DefLocalSir} and~\ref{DefSFC}).
I generalise Saharon Shelah's idea of passing from dividing to
forking as follows:
Whenever $\ind$, a candidate for being an independence relation,
is in fact a pre-independence relation, then $\sind$, derived
from $\ind$ as in Definition~\ref{DefSind}, is a better candidate in the
sense that $\ind=\sind$ if $\ind$ is already an independence relation,
and $\sind$ is always a pre-independence relation satisfying the extension axiom.
(Hence local character is the only axiom that need not hold for $\sind$.)

I show that Shelah-dividing independence $\dind$
and M-dividing independence $\mind$ from~\cite{Hans Scheuermann: Unabhaengigkeitsrelationen}
are pre-independence relations. For certain sets $\Omega$
of pairs of formulas I also define $\Omega$-dividing $\omind$
(Shelah-dividing localised to $\Omega$, cf.~Definition~\ref{DefOmega})
and show that it is a pre-independence relation.
$\Xi$ and $\Xim$, two such sets $\Omega$, are defined in Definitions~\ref{DefXi} and~\ref{DefinitionXim}.
I show that, if $\find$ denotes Shelah-forking independence and $\thind$
denotes thorn-forking independence, we have
$\find=\ind[d*]=\ind[$\Xi*$]$\; and $\thind=\ind[M*]=\ind[$\Xim*$]$\;\;\;.
Since $\mind$\, and $\ind[$\Xim$]$\; have simpler definitions than thorn-dividing,
this helps to understand thorn-forking.

In the following I go into more detail and also describe some other related results.

\subsubsection*{M-dividing, and a simple definition of thorn-forking}

M-dividing from~\cite{Hans Scheuermann: Unabhaengigkeitsrelationen} is defined as follows:
$A\mind_CB$ holds iff for all sets $C'$ such that $C\subseteq C'\subseteq\acl(BC)$
we have $\acl(AC')\cap\acl(BC')=\acl C'$.
Note that with this definition we need to take care whether we evaluate
algebraic closure $\acl$ in~$T$ or in~$T\eq$.
(The same is true for strong dividing and thorn-dividing, but not
for Shelah-dividing.)
In~\cite{Hans Scheuermann: Unabhaengigkeitsrelationen} it was observed that M-dividing is closely
related to `modular pairs' in the lattice of algebraically closed sets,
and that $\mind$\, as evaluated on the real elements of a pregeometric
theory (e.\,g., strongly minimal or o-minimal) is the familiar notion of independence.
Here I continue the study of M-dividing by proving:

\begin{itemize}
\item $\mind$ is an independence relation iff $\mind$ is symmetric.
	(This answers a question implicit in~\cite{Hans Scheuermann: Unabhaengigkeitsrelationen}.)
\item Thorn-forking is the notion of forking corresponding to M-dividing,
	i.\,e., $\thind=\ind[M*]$\;.
\item If there is an independence relation $\ind$ for $T\eq$
	satisfying the condition $a\ind_Ca\implies a\in\acl\eq C$, then $T$ is rosy and
	$\thind$ is the coarsest independence relation for $T\eq$ satisfying this condition.
\end{itemize}
The last result is dual to the situation with Shelah-forking, which, in a
simple theory, is the finest independence relation.

\subsubsection*{Thorn-forking via localised Shelah-dividing}

Forking is traditionally defined via local dividing, i.\,e., dividing of formulas.
The notion of $k$-dividing of a formula as it appears in Byunghan Kim's
treatment of simple theories can be seen as a more thorough localisation.
In the same sense, dividing of a formula $\phi$ with respect to a $k$-inconsistency
witness $\psi$ for $\phi$ (introduced by Itay Ben-Yaacov) is even more radically local.

Back-porting some of Ben-Yaacov's ideas into the elementary context,
I examine generalised local dividing with respect to a set $\Omega$ of `inconsistency
pairs' $(\phi,\psi)$ where $\psi$ is an inconsistency witness for~$\phi$.
The only interesting cases I know are $\Omega=\Xi$ (the set of \emph{all} inconsistency pairs)
and $\Omega=\Xim$ as defined in Definition~\ref{DefinitionXim}. In the second chapter I show:
\begin{itemize}
	\item $\ind[$\Xi$]=\dind$, hence $\ind[$\Xi*$]=\find$,
		so Shelah-forking is a special case of $\Omega$-forking.
	\item $\ind[$\Xim*$]\;\;=\thind$,
		so thorn-forking is a special case of $\Omega$-forking.
\end{itemize}
I also define local $\D_\Delta$-ranks of types for finite $\Delta\subseteq\Omega$,
isolate two technical conditions which $\Omega$ may satisfy (`transitivity' and `normality')
and prove:
\begin{itemize}
	\item $\Xi$ and $\Xim$ are transitive and normal.
	\item If $\Omega$ is transitive and normal, then $\omind$ is a pre-independence relation.
	\item If $\Omega$ is transitive and normal and all $\D_\Delta$-ranks are finite,
	then $\bar a\omsind_CB$ holds iff $\D_\Delta(\bar a/BC)=\D_\Delta(\bar a/C)$.
\end{itemize}

\subsubsection*{Kernels and canonical bases}

The weak canonical base of a type over an algebraically closed set is the smallest
algebraically closed
subset over which the type is free. In the third chapter I show that this concept
is closely related to thorn-forking:
\begin{itemize}
\item If an independence relation (satisfying the anti-reflexivity condition
	$a\ind_Ca\implies a\in\acl C$) has weak canonical bases, then it is thorn-forking independence.
\end{itemize}
As important tools for studying (weak) canonical bases I define the kernel and the
algebraic kernel of a sequence of indiscernibles:
The (algebraic) kernel of an infinite sequence of indiscernibles is the greatest subset of its
definable (resp.~algebraic) closure over which it is still indiscernible.
I show:
\begin{itemize}
\item Every infinite sequence of indiscernibles has a kernel and an algebraic kernel, and they are 
	invariant under `collinearity.'
\item For $\thind$ to have weak canonical bases the following condition is sufficient:
	$A\thind_{C_1}B$, $A\thind_{C_2}B$ and $C_1,C_2\subseteq B$ together imply $A\thind_{\acl C_1\cap\acl C_2}B$.
\item If $\thind$ has weak canonical bases, then the weak canonical base of a type can be computed as the algebraic kernel
	of an arbitrary Morley sequence.
\item In a stable theory the canonical base of a stationary type can be computed
	as the kernel of its Morley sequence.
	In simple theories the situation is slightly more complicated.
\item If a sequence of indiscernibles has a canonical base in the sense
	of Buechler, then the canonical base coincides with the kernel.
\end{itemize}

\subsubsection*{Some old results}

The author of~\cite{Hans Scheuermann: Unabhaengigkeitsrelationen} never formally published his results.
Some of them are generalised in Chapter~\ref{Chapter3} or treated in exercises.
I feel that two of them should be mentioned here:
\begin{itemize}
\item Let $T$ be a simple theory and let $T'$ be a reduct of $T$ that has elimination of hyperimaginaries.
If $C=\acl\eq C$ in $T$ and $A\ind_CB$ holds in $T$, then $A\ind_CB$ holds in $T'$
(Exercise~\ref{ExcCanonicalReduct}).
\item
Let $T$ be a simple theory with elimination of hyperimaginaries.
$T$ is 1-based iff the lattice of algebraically closed sets is modular (Exercise~\ref{ExcOneBased}).
\end{itemize}

\subsubsection*{Acknowledgements}

I thank Nina Frohn, Immanuel Halupczok, Markus Junker, Herwig N\"ubling,
Jan-Christoph Puchta, Mark Weyer and of course my supervisor Martin Ziegler
for helpful discussions and comments on drafts of this thesis.
I particularly thank Anand Pillay, who was so extraordinarily fast to send his comments
on a very late draft of this thesis (answering three of my open questions)
that I was still able to take them into account.
Needless to say, all errors, inaccuracies and excentricities that have survived are my own.

\subsubsection*{Preliminaries}

Readers are assumed to be acquainted with the culture of stability theory.
It should be no surprise to them that we work in a big saturated model
of a complete consistent theory. That we compute definable closure $\dcl$
and algebraic closure $\acl$ either by means of definable or algebraic
formulas, or, equivalently, by means of automorphisms of the big model.
That we work with many-sorted theories, such as $T\eq$, whenever we feel
like it. And that we work with indiscernible sequences, which are always
implicitly assumed to be infinite.

It should not be hard to learn these things from the first pages
of a book like~\cite{Anand Pillay: Geometric Stability Theory}.
\chapter{Abstract independence}
\label{Chapter1}

\noindent In this chapter an axiomatic treatment of independence relations for a
complete consistent first-order theory is presented.
Some properties of forking and thorn-forking are proved in this context.
Thorn-forking is defined in a new way---via M-dividing,
a new notion introduced in Section~\ref{SectionThornForking}.

While the geometric theory of forking is usually based on a combinatorial
foundation, we will see that a geometric treatment is possible from the beginning.
This approach is faster and probably more comprehensible than the usual
combinatorial one, but we get slightly weaker results.
We will improve them in the next chapter by means of the usual, more combinatorial, methods.

The exposition in this chapter is self-contained apart from the general
cultural background of stability theory:
Except in the notes at the end of each section, no other knowledge
from stability theory is assumed.
Tuples of elements ($\bar a,\bar b,\bar c,\dots$)  or of variables
($\bar x,\bar y,\bar z,\dots$) are allowed to be infinite unless mentioned
otherwise. When a formula is written $\phi(\bar x)$ it means that each of
its free variables appears in the (possibly infinite) tuple~$\bar x$.
There is an endless supply of formal variables that we can use in types.
$\S^*(B)$ denotes the class of complete types over $B$ in arbitrarily long
sequences of variables.

I write $AB$ for $A\cup B$, and for any tuple $\bar a=(a_0,a_1,\ldots,)$
I will abuse notation by writing $\bar a$ for the set
$\{a_0,a_1,\ldots,\}$ as well.
$A\equiv_BA'$ means that there is an automorphism of the big models
that fixes $B$ pointwise and maps the set $A$ to the set $A'$.
$\bar a\equiv_B\bar a'$ means that there is an automorphism fixing $B$
pointwise and mapping the tuple $\bar a$ to the tuple $\bar a'$.
$(A,B)\equiv_C(A',B')$ means that there is an automorphism fixing $C$
pointwise that maps $A$ to $A'$ and $B$ to $B'$.


\section{Just a bunch of silly axioms?}

We will call a ternary relation $\ind$
between (small) subsets of the big model an \emph{independence relation} if it
satisfies the axioms of the following definition:

\begin{definition}\label{DefLocalSir}\index{independence relation}\index{006@$\ind$}
	The following are the \emph{axioms for independence relations:}
	\begin{description}
		\item[(invariance)]
			\index{axiom!invariance}
			\index{invariance}\ \\
			If $A\ind_CB$ and $(A',B',C')\equ (A,B,C)$, then $A'\ind_{C'}B'$.
		\item[(monotonicity)]
			\index{axiom!monotonicity}
			\index{monotonicity}\ \\
			If $A\ind_CB$, $A'\subseteq A$ and $B'\subseteq B$, then $A'\ind_CB'$.
		\item[(base monotonicity)]
			\index{axiom!base monotonicity}
			\index{base monotonicity}\ \\
			Suppose $D\subseteq C\subseteq B$.
			If $A\ind_DB$, then $A\ind_CB$.
		\item[(transitivity)]
			\index{axiom!transitivity}
			\index{transitivity}\ \\
			Suppose $D\subseteq C\subseteq B$.
			If $B\ind_CA$ and $C\ind_DA$, then $B\ind_DA$.
		\item[(normality)]
			\index{axiom!normality}
			\index{normality}\ \\
			$A\ind_CB$ implies $AC\ind_CB$.
		\item[(extension)]
			\index{axiom!extension}
			\index{extension}\ \\
			If $A\ind_CB$ and $\hat B\supseteq B$, then there is $A'\equ_{BC}A$ such that $A'\ind_C\hat B$.\\
			(Equivalently, by invariance, there is $\hat B'\equ_{BC}\hat B$ such that $A\ind_C\hat B'$.)
		\item[(finite character)]
			\index{axiom!finite character}
			\index{finite character}\ \\
			If $A_0\ind_CB$ for all finite $A_0\subseteq A$, then $A\ind_CB$.
		\item[(local character)]
			\index{axiom!local character}
			\index{local character}\ \\
			For every $A$ there is a cardinal $\kappa(A)$ with the following property:\\
			For any set $B$ there is a subset $C\subseteq B$ of cardinality $|C|<\kappa(A)$
			such that $A\ind_CB$.
	\end{description}
\end{definition}

\begin{definition}
	An independence relation is \emph{strict} if it also satisfies the
	following axiom:\index{independence relation!strict}\index{strict independence relation}
	\begin{description}
		\item[(anti-reflexivity)]
			\index{axiom!anti-reflexivity}
			\index{anti-reflexivity}\ \\
			$a\ind_Ba$ implies $a\in\acl B$.
	\end{description}
\end{definition}

\noindent For a first example, we need look no further than the coarsest\footnote{Generalising
common usage for equivalence relations and topologies,
we say that a relation $\ind$ is finer\index{finer relation} than $\pind$,
and $\pind$ coarser\index{coarser relation} than $\ind$, if $A\ind_CB$ implies $A\pind_CB$.}
relation of all:  the trivial relation that holds for all triples.
It satisfies all of the above axioms except anti-reflexivity.
So the trivial relation is always a (non-strict) independence relation.

In practice, when I say `by transitivity' I will often mean the following
stronger property (or a variant thereof):

\begin{remark}\label{RemarkTransitivity}
	Let $\ind$ be a relation satisfying monotonicity, transitivity and normality.
	Then $B\ind_{CD}A$ and $C\ind_DA$ together imply $BC\ind_DA$.
\end{remark}

\begin{proof}
	If $B\ind_{CD}A$ and $C\ind_DA$, then
	$BCD\ind_{CD}A$ and $CD\ind_DA$ by normality.
	Hence $BCD\ind_DA$ by transitivity, so $BC\ind_DA$ by monotonicity.
\end{proof}

\begin{example}\label{ExForest1}
	\index{example!everywhere infinite forest}
	(Everywhere infinite forest)\\
	Let $T$ be the theory, in a signature with one binary relation~$E$,
	of a non-empty undirected tree that branches infinitely in every node.
	Then $T$ is complete, and the models of $T$ are precisely the non-empty
	forests that branch infinitely in every node.
	In this theory, $\acl A$ is the set of all nodes that lie on
	a path between two elements of~$A$.
	
	Consider the following relation:
	\[
		A\ind_CB \quad\iff\quad
			\textsl{every path from $A$ to $B$ meets $\acl C$.}
	\]
	It is easy to see that $A\ind_CB$ implies $AC\ind_CB$ and $\acl A\ind_CB$.
	Using this, it is not hard to check that $\ind$ is a \sir.
	The details are left to the reader as an exercise
	(Exercise~\ref{ExcForest1}).
\end{example}

If we look for more general \sir s with Exercise~\ref{ExcModular} in mind,
we will tend to find \emph{coarse} \sir s such as thorn-forking.
In the next section we will introduce Morley sequences. These will
motivate us to try another approach that will lead us
to \emph{fine} \sir s such as Shelah-forking.

But first we observe that $A\ind_CB\iff B\ind_CA$ holds both for
$\aind$ from Exercise~\ref{ExcModular}
below and for $\ind$ from Example~\ref{ExForest1}.
In the next section we will examine whether this is an accident.

\begin{exercises}
	Solutions for all exercises can be found in the appendix in
	Section~\ref{SectionExercises}.
	\begin{exercise}(relations between the axioms, existence and symmetry)\label{ExcExistence}
	
		Consider the following additional properties which a relation $\ind$
		may satisfy:
		\begin{description}
			\item[(existence)]\index{axiom!existence}\index{existence}
				For any $A$, $B$ and $C$ there is $A'\equ_CA$ such that $A'\ind_CB$.
			\item[(symmetry)]
				$A\ind_CB\iff B\ind_CA$
		\end{description}
		
		(i) Any relation that satisfies invariance, extension and symmetry also satisfies normality.
		
		(ii) Any relation that satisfies invariance, extension and local character
		also satisfies existence.
		
		(iii) Any relation that satisfies invariance, monotonicity, transitivity, normality, existence
		and symmetry also satisfies extension. 
	\end{exercise}
	
	\begin{exercise}(local character)\label{ExcLocalCharacter}
	
		(i) Suppose the relation $\ind$ satisfies invariance and the existence
		condition from Exercise~\ref{ExcExistence}.
		Let $\kappa(A)=(\card T + \card A)^+$.
		For any sets $A$ and $B$ there is $C_1\subseteq B$ such that
		$A\ind_{C_1}B$ and also a set $C_2$ such that $\card{C_2}\leq\kappa(A)$
		and $A\ind_{C_2}B$.
		
		(ii) Suppose $\ind$ is an independence relation.
		Let $\mathcal A$ be a set of finite subsets of the big model such that
		for every finite subset $A$ of the big model there is $A'\in\mathcal A$
		such that $A\equiv A'$. Let $\kappa=\operatorname{sup}_{A\in\mathcal A}\kappa(A)$.
		Show that for any sets $A$ and $B$ there is $C\subseteq B$ such that
		$|C|<\kappa+|A|^+$ and $A\ind_CB$.
	\end{exercise}
	
	\begin{exercise}\label{ExcModular}(modularity and distributivity)\index{008@$\aind$}
	
		Consider the following relation:
		\[ A\aind_CB\iff\acl(AC)\cap\acl(BC)=\acl C. \]
		
		(i) The relation $\aind$ satisfies
		the existence condition from Exercise~\ref{ExcExistence} (hard).
		
		(ii) The relation $\aind$ satisfies all axioms for \sir s except base monotonicity.
		
		(iii) The relation $\aind$ satisfies base monotonicity
		(and is a \sir) iff the lattice of algebraically closed\index{lattice!modular or distributive}\index{modular lattice}
		subsets of the big model is \emph{modular,} i.\,e.,
		whenever $A,B,C$ are algebraically closed sets such that
		$B\supseteq C$, the equation $B\cap\acl(AC)=\acl((B\cap A)C)$ holds.
		
		(iv) An independence relation is
		\emph{perfectly trivial}\index{independence relation!perfectly trivial}\index{perfectly trivial}
		if $A\ind_CB$ implies $A\ind_{C'}B$ for all $C'\supseteq C$.
		Suppose $\aind$ is an independence relation.
		Show that $\aind$ is a perfectly trivial independence relation if and only if
		the lattice of algebraically closed sets is \emph{distributive,}\index{distributive lattice}
		i.\,e., $\acl((A\cap B)C)=\acl((A\cap C)(B\cap C))$ holds for all algebraically closed
		sets $A$, $B$ and $C$.
	\end{exercise}
	
	\begin{exercise}\ \label{ExcForest1}
	
		For Example~\ref{ExForest1}, check that $A\ind_CB$ implies $\acl A\ind_CB$
		and that $\ind$ is a \sir.
	\end{exercise}
\end{exercises}

\begin{notes}
	The symbol `$\ind$' was first used for independence by Michael Makkai in~\cite{Michael Makkai: A survey of basic stability theory},
	but the symbol $\perp$ was used in a similar context in lattice theory much earlier.\footnote{If
		we apply the definition of $A\perp B$ from~\cite{John von Neumann (ed. I. Halperin): Continuous Geometries},
		a book based on John von Neumann's work on lattice theory in the 1930s,
		to the lattice of algebraically closed sets it means $A\aind_\emptyset B$.
		Note that the the original context was a very special type of lattices which were,
		in particular, modular.
		$(A,B)\perp$ in~\cite{Lee Roy Wilcox: Modularity in Lattices}, when translated in the same way,
		means $A\mind_\emptyset B$, cf.~Definition~\ref{DefMThorn}.}
	A~possible pronunciation of $\ind$ is `anchor'.
	
	Independence is often expressed using another notation that is equivalent
	to the $\ind$ notation: For $C\subseteq B$
	and complete types $p(\bar x)\in \S(C)$ and $q(\bar x)\in \S(B)$,
	write $p\sqsubseteq q$ if $p\subseteq q$ and $\bar a\ind_CB$, where $\bar a$ realises~$q$.
	An axiomatic characterisation of (classical) forking independence
	in a stable theory
	in terms of such a relation was discovered by Victor Harnik and
	Leo~A.~Harrington~\cite{Harnik + Harrington: Fundamentals of Forking}.
	The next big break-through in this direction was the core result
	of Byunghan Kim's dissertation~\cite{Byunghan Kim: Forking in simple unstable theories} (see also~\cite{Kim + Pillay: Simple theories}):
	an axiomatic characterisation of (Shelah-)forking independence in a simple theory.
	
	The terms `independence relation' and `notion of independence' are not (yet) standardised.
	Some authors include certain axioms of varying strength (boundedness or the amalgamation property)
	that make sure that if there is an independence relation at all, then it is unique
	and coincides with classical forking independence (and the theory is stable or simple,
	respectively). Apart from that, the axiomatic systems are usually equivalent to
	the system of Definition~\ref{DefLocalSir}.
	
	Some axioms appear in a slightly unusual form here.
	The transitivity\index{axiom!transitivity}\index{transitivity} axiom from~\cite{Michael Makkai: A survey of basic stability theory}, for instance,
	was dualised ($A$ on the right-hand side) and separated into
	transitivity and normality.\index{axiom!normality}\index{normality} (The term `normality' is new.)
	Many authors use another variant sometimes called `full transitivity'.
	The terminology around extension\index{axiom!extension}\index{extension}
	and existence\index{axiom!existence}\index{existence} is also far from unified.
	This probably dates back to~\cite{Michael Makkai: A survey of basic stability theory}, where `existence'
	combines existence and extension into one axiom.
	Local character\index{axiom!local character}\index{local character} was strengthened so it is
	more useful when we do not have finite character.
	
	For Example~\ref{ExForest1} cf.{} the note to Example~\ref{ExReductNotMSymmetric}
	below. Exercise~\ref{ExcModular} is from \cite{Hans Scheuermann: Unabhaengigkeitsrelationen};
	the definition of $\aind$ was inspired
	by~\cite{Lee Fong Low: Lattice of algebraically closed sets in one-based theories}. Exercise~\ref{ExcExistence} presents three easy relations
	that hold between the axioms of independence relations and existence and symmetry.
	Proving a fourth relation is the main object of the next section, while
	the fact that these are the only relations is the subject of Section~\ref{SectionExamples}.
\end{notes}


\section{Morley sequences and symmetry}

We now prove that every independence relation is symmetric.
But we need some preparations first.

\begin{proposition}\label{PropShifting}
	Let $\ind$ be an independence relation.\\
	If $A\ind_DBC$ and $B\ind_DC$, then $AB\ind_DC$.\\
	Actually, it is sufficient that $\ind$ satisfies
	the first five axioms.
\end{proposition}

\begin{proof}
	$A\ind_DBC$ implies $A\ind_DBCD$ by extension and invariance,
	hence $A\ind_{BD}BCD$ by base monotonicity,
	hence $A\ind_{BD}C$ by monotonicity.
	Combining this with $B\ind_DC$, we get $AB\ind_DC$
	by transitivity (i.\,e., by monotonicity, transitivity,
	normality and Remark~\ref{RemarkTransitivity}).
\end{proof}

\begin{definition}
	Let $\ind$ be a ternary relation.\index{Morley sequence}\\
	A \emph{$\ind$-Morley sequence in} a type $p(\bar x)\in\S^*(B)$ over a set $C\subseteq B$ is a sequence of
	$B$-indiscernibles $(\bar a_i)_{i<\omega}$ such that $(\bar a_i)_{i<n}\ind_C\bar a_n$
	for every $n<\omega$, and every $\bar a_i$ realises $p(\bar x)$.\\
	A \emph{$\ind$-Morley sequence for} a complete type $p(\bar x)\in\S^*(B)$ is
	a $\ind$-Morley sequence in $p(\bar x)$ over~$B$.
\end{definition}

Recall our convention that tuples may be infinite. In most cases just a
convenience, this is crucial in this section and the next one.
The following consequence of the Erd\H{o}s-Rado theorem is proved
in~\cite[Lemma 1.2]{Itay Ben-Yaacov: Simplicity in compact abstract theories} in the more general context
of compact abstract theories (cats):

\begin{fact}\label{FactExtraction} \textbf{(Extracting a sequence of indiscernibles)}\\
	Let $B$ be a set of parameters and $\kappa$ a cardinal.\index{indiscernibles!extracting}\index{extracting indiscernibles}
	Then for any sequence $(\bar a_i)_{i<\beth_{(2^{\card T+\card B+\kappa})^+}}$ consisting of sequences
	of length $\card{\bar a_i}=\kappa$
	there is a  $B$-indiscernible sequence $(\bar a'_j)_{j<\omega}$ with
	the following property:\\
	For every $k<\omega$ there are $i_0<i_1<\ldots<i_k<\kappa$ such that
	$\bar a_{i_0},\bar a_{i_1},\ldots,\bar a_{i_k}\equiv_B\bar a'_0,\bar a'_1,\ldots,\bar a'_k$.
\end{fact}

\begin{proposition}\label{PropExMS}
	Suppose $\ind$ is an independence relation and $\bar a\ind_CB$.
	Then there is a Morley sequence in~$\tp(\bar a/BC)$ over~$C$.\\
	Actually, it is sufficient that $\ind$ satisfies
	the first five axioms and extension.
\end{proposition}

\begin{proof}
	Let $\bar a_0=\bar a$.
	We choose a very big cardinal $\kappa$ and use extension and transfinite induction
	to construct a sequence $(\bar a_i)_{i<\kappa}$ satisfying
	$\bar a_{\alpha}\equ_{BC}\bar a_0$ and
	$\bar a_{\alpha}\ind_C(\bar a_i)_{i<\alpha}$ for all $\alpha<\kappa$.
	If $\kappa$ has been chosen	sufficiently big, we can extract a sequence of
	$BC$-indiscernibles $(\bar a'_i)_{i<\omega}$ such that for every $n<\omega$ there are
	indices $\alpha_0$, \ldots, $\alpha_n$ such that
	$\bar a'_0\ldots\bar a'_n\equ_{BC}\bar a_{\alpha_0}\ldots\bar a_{\alpha_n}$.
	Note that $\bar a'_n\ind_C(\bar a'_i)_{i<n}$ by monotonicity and invariance.
	
	By compactness we can `invert' the sequence $(\bar a'_i)_{i<\omega}$, i.\;e.,
	find a new sequence $(\bar a''_i)_{i<\omega}$ such that
	$\bar a''_0\ldots\bar a''_n\equ_{BC}\bar a'_n\ldots\bar a'_0$ holds for every $n<\omega$.
	In particular, the new sequence is also indiscernible over $BC$.
	This new sequence satisfies $\bar a''_0\ind_C(\bar a''_i)_{0<i<n}$ for all $n<\omega$.
	Hence $(\bar a''_i)_{i<n}\ind_C\bar a_n$ for all $n<\omega$ by repeated use of
	Proposition~\ref{PropShifting}.
	
	Thus $(\bar a''_i)_{i<\omega}$ is a Morley sequence in $\tp(\bar a/BC)$ over~$C$.
\end{proof}

\begin{proposition}\label{PropConsMS}
	Suppose $\ind$ is an independence relation and there is a Morley sequence
	in~$\tp(\bar a/BC)$ over~$C$. Then $B\ind_C\bar a$.\\
	Actually, it is sufficient that $\ind$ satisfies
	the first five axioms, finite character and local character.
\end{proposition}

\begin{proof}
	Let $(\bar a_i)_{i<\omega}$ be a Morley sequence in~$\tp(\bar a/BC)$ over~$C$.
	
	Let $\kappa$ be a regular cardinal number that is greater than or equal to $\kappa(B)$
	as in the local character axiom. By compactness we can extend the sequence $(\bar a_i)_{i<\omega}$
	to obtain a $BC$-indiscernible sequence $(\bar a_i)_{i<\kappa}$.
	Using finite character, we see that
	$(\bar a_i)_{i<\alpha}\ind_C\bar a_{\alpha}$ for each $\alpha<\kappa$.

	By local character there is a set $D\subseteq C(\bar a_i)_{i<\kappa}$ of
	cardinality $|D|<\kappa$ such that $B\ind_DC(\bar a_i)_{i<\kappa}$.
	By regularity of $\kappa$ there is an index $\alpha<\kappa$ such that already
	$D\subseteq C(\bar a_i)_{i<\alpha}$. Thus, by base monotonicity and monotonicity,
	$B\ind_{C(\bar a_i)_{i<\alpha}}\bar a_{\alpha}$.
	
	Combining the results of the last two paragraphs using transitivity
	(actually, Remark~\ref{RemarkTransitivity} and monotonicity), we get $B\ind_C\bar a_{\alpha}$.
	Since $\bar a_{\alpha}\equ_{BC}\bar a$ this implies $B\ind_C\bar a$
	by invariance.
\end{proof}

For later use I have carefully noted which axioms were actually needed to prove
Propositions~\ref{PropExMS} and~\ref{PropConsMS}.
For our immediate use of them in this section it would not have been necessary:

\begin{theorem}\label{ThmSymmetry}
	Every independence relation $\ind$ is symmetric:\index{axiom!symmetry}\index{symmetry axiom}\ \\
	For any $A$, $B$ and $C$, $A\ind_CB$ iff $B\ind_CA$.
\end{theorem}

\begin{proof}
	If $\bar a\ind_CB$, then there is a Morley sequence in~$\tp(\bar a/BC)$ over~$C$
	by Proposition~\ref{PropExMS}.
	Hence $B\ind_C\bar a$ by Proposition~\ref{PropConsMS}.
\end{proof}

From now on we may use symmetry implicitly when working with an independence relation.
For the rest of this chapter, however, we focus on \emph{finding} independence relations.

\begin{example}\label{ExNoSir}
	\index{example!no strict indep.{} relation}
	(A theory with no \sir)\\
	Consider the following two-sorted theory $T_0$:
	$T_0$ has two sorts $P$ and $E$, the elements of which are
	called `points' and `equivalence relations', and a single
	ternary relation $\sim\;\subseteq P\times P\times E$ written as $p\sim_eq$.
	The axioms of $T_0$ say that $\sim_e$ is
	an equivalence relation on the points for every $e\in E$.
	
	Clearly every substructure of a model of $T_0$ is again a
	model of $T_0$. The signature of $T_0$ is finite and relational.
	Moreover, the class of finite models of $T_0$ satisfies the
	joint embedding property and the amalgamation property.
	So by~\cite[Theorem 7.4.1]{Wilfrid Hodges: Model Theory} $T_0$ has a Fra\"{\i}ss\'e limit~$T^*$
	which is $\omega$-categorical, has quantifier elimination,
	and whose finite submodels are precisely the
	finite models of $T_0^*$.
	
	Since $\acl A=A$ for all sets, $A\ind_CB \iff A\cap B\subseteq C$
	defines a \sir{} for~$T$.
	But there is no \sir{} for~$T\eq$:
	
	Suppose $\ind$ is an independence relation for~$T\eq$.
	Let $a_0\in P$ be a single point,
	and let $(a_i)_{i<\omega}$ be a Morley sequence for $\tp(a_0/\emptyset)$.
	Let $e\in E$ be an equivalence relation such that
	$a_i\sim_ea_j$ for any $i,j<\omega$.
	Then $(a_i)_{i<\omega}$ is indiscernible over~$e$.
	
	Note that $(a_{2i}a_{2i+1})_{i<\omega}$ is a Morley sequence
	for $\tp(a_0a_1/\emptyset)$ and is also indiscernible over~$e$.
	By Proposition~\ref{PropConsMS}, $e\ind_{\emptyset}a_0a_1$ holds,
	so by base monotonicity, $e\ind_{a_0}a_1$.
	On the other hand, $a_0\ind_{\emptyset}a_1$ also holds.
	Applying transitivity we get $a_0e\ind_{\emptyset}a_1$.
	Using symmetry and base monotonicity, we get $a_0\ind_ea_1$.
	
	But the equivalence class $c$ of $a_0$ and $a_1$ under
	$\sim_e$ is (an element of $T\eq$ and)
	definable both over $a_0e$ and over $a_1e$. 
	So $c\ind_ec$. Since $c$ is clearly not algebraic over~$e$
	this contradicts anti-reflexivity.
	So $\ind$ is not a \sir{} for $T\eq$.
\end{example}

\begin{notes}
	While the terms and techniques used in this section are not new,
	the specific treatment of abstract independence, and
	Theorem~\ref{ThmSymmetry} in particular, seems to be
	new.\index{independence relation!note on terminology}\index{notion of independence}
	Once you have the right set of axioms, it is somewhat implicit
	in the work of Byunghan Kim. I first met the technique of
	Proposition~\ref{PropConsMS} in~\cite{Byunghan Kim: Forking in simple unstable theories}.
	
	Theorem~\ref{ThmSymmetry} is made possible by the fact that transitivity is
	postulated on the left-hand side in Definition~\ref{DefLocalSir}.
	With the usual transitivity axiom\index{axiom!transitivity}\index{transitivity}
	(on the right-hand side, i.e.: $A\ind_CB$ and $A\ind_DC$ implies $A\ind_DB$
	for $D\subseteq C\subseteq B$),
	Theorem~\ref{ThmSymmetry} would not hold.
	That is why symmetry is traditionally included as an axiom for
	independence relations.
	This point is usually obscured by authors who combine monotonicity,
	base monotonicity and right-hand side transitivity into an axiom called
	`full transitivity'.
	I was set on the right track by~\cite[Corollary 1.9]{Itay Ben-Yaacov: Simplicity in compact abstract theories}.
	
	Example~\ref{ExNoSir} was suggested to me by Martin Ziegler.
\end{notes}


\section{A theorem on abstract forking}

We now show that, in a certain sense, the extension axiom comes free.

\begin{definition}\label{DefSind}
	For any invariant relation $\ind$ we define a new relation $\sind$ as follows:\index{012@$\sind$}
	\[
		A\sind_CB \iff \big( \textsl{ for all $\hat B\supseteq B$ there is
		$A'\equ_{BC}A$ s.t.\ } A'\ind_C\hat B \big).
	\]
\end{definition}

\noindent Note that $A\sind_CB$ implies $A\ind_CB$. Also, $\ind=\sind$ iff $\ind$ satisfies
the extension axiom. In analogy to the classical situation one might
call $\sind$ the notion of forking derived from the abstract notion
of dividing given by~$\ind$.

If $\ind$ already satisfies some of the axioms for independence relations,
then there cannot be much harm (possibly losing finite character and local character)
in passing from $\ind$ to $\sind$, but we get extension free:

\begin{lemma}\label{LemmaDF}\index{axiom!extension}\index{extension}\ \\
	If $\ind$ is a relation satisfying invariance and montonicity,
	then $\sind$ satisfies invariance, monotonicity and extension.
	If, moreover, $\ind$ satisfies one of the following axioms and properties,
	then $\sind$ also satisfies it:
	base monotonicity, transitivity, normality, anti-reflexivity, existence.
\end{lemma}

\begin{proof}
	Invariance of $\sind$ is obvious.
	
	\emph{Monotonicity:}
	Suppose $A\sind_CB$, $A_0\subseteq A$ and $B_0\subseteq B$.
	Then for all $\hat B\supseteq B$ there is $A'\equ_{BC}A$ such that
	$A'\ind_C\hat BC$. Let $A_0'\subseteq A'$ correspond to $A_0\subseteq A$.
	Then clearly $A_0'\equ_{B_0C}A_0$ and $A_0'\ind_C\hat B$.
	Thus $A_0\sind_CB_0$ holds.
	
	\emph{Extension:}
	Suppose $\bar a\sind_CB$, where $\bar a$ is a possibly infinite tuple,
	and let $\hat B\supseteq B$ be any superset of~$B$.
	
	We first claim that there is a type $p(\bar x)\in\S^*(\hat BC)$, extending
	$\tp(\bar a/BC)$, such that for all cardinals $\kappa$ there is a
	$\kappa$-saturated model $M\supseteq \hat BC$ and $\bar a'\models p(\bar x)$
	such that $\bar a'\ind_CM$.
	
	If not, then for each $p(\bar x)\in\S^*(\hat BC)$ extending $\tp(\bar a/BC)$
	there is a cardinal
	$\kappa(p)$ such that for no $\kappa(p)$-saturated model
	$M\supseteq\hat BC$ is there a tuple $\bar a'\models p$
	satisfying $\bar a'\ind_CM$.
	Let $\kappa$ be the supremum of the cardinals $\kappa(p)$,
	and let $M\supseteq\hat BC$ be $\kappa$-saturated.
	Then there is no $\bar a'\equ_{BC}\bar a$ such that $\bar a'\ind_CM$.
	So we have found a contradiction to the definition of $\sind$
	and thereby proved our claim.
	
	Now choose $\bar a'\models p(\bar x)$, where $p(\bar x)$ is as in the
	claim. Then clearly $\bar a'\equ_{BC}\bar a$, and
	$\bar a'\sind_C\hat B$ by monotonicity of~$\ind$.

	\emph{Base monotonicity:}
	Suppose $A\sind_CB$ and $C\subseteq C'\subseteq B$.
	So for any $\hat B\supseteq B$ there is $A'\equ_{BC}A$
	such that $A'\ind_C\hat BC$.
	Base monotonicity of $\ind$ yields $A'\ind_{C'}\hat BC$, so $A'\ind_{C'}\hat B$
	by monotonicity of~$\ind$.
	Thus $A\sind_{C'}B$.
	
	\emph{Transitivity:}
	Here we work with an alternative definition of $\sind$, which is equivalent
	to Definition~\ref{DefSind} by invariance of $\ind$:
	\[
		A\sind_CB \iff \big( \text{ for all $\hat B\supseteq B$ there is
		$\hat B'\equ_{BC}\hat B$ s.t.\ } A\ind_C\hat B' \big).
	\]
	
	So suppose $D\subseteq C\subseteq B$, $B\sind_{C}A$ and $C\sind_DA$ hold,
	and $\hat A\supseteq A$.
	We need to show that $B\ind_D\hat A^*$ for some $\hat A^*\equ_{AD}\hat A$.
	
	Let $\hat A'\equ_{AD}\hat A$ be such that $C\ind_D\hat A'$,
	and let $\hat A^*\equ_{AC}\hat A'$ be such that $B\ind_C\hat A^*$.
	Note that $\hat A^*\equ_{AD}\hat A$ and $C\ind_D\hat A^*$.
	By transitivity of $\ind$ we get $B\ind_D\hat A^*$.
	
	\emph{Normality:}
	Suppose $A\sind_CB$ and $\hat B\supseteq B$.
	Let $A'\equiv_{BC}A$ be such that $A'\ind_C\hat B$.
	Then also $A'C\ind_C\hat B$ by normality of~$\ind$.
	
	\emph{Anti-reflexivity:}
	Trivial, since $A\sind_CB$ implies $A\ind_CB$.
	
	\emph{Existence:}
	Suppose we are given $A,B,C$.
	Since $\ind$ satisfies existence by assumption, we have $A\sind_C\emptyset$.
	Since $\sind$ satisfies extension there is $A'\equiv_CA$ such that
	$A'\sind_CB$.
\end{proof}

\begin{theorem}\label{ThmDivFork}
	Suppose $\ind$ satisfies the first five axioms for independence relations
	and also finite character.
	Suppose $\sind$, derived from $\ind$ as in Definition~\ref{DefSind},
	has local character.
	Then $\sind$ is an independence relation.
\end{theorem}

\begin{proof}
	It follows from Lemma~\ref{LemmaDF} that $\sind$ satisfies
	the first five axioms and extension.
	As local character holds by assumption, we need only prove that
	$\sind$ satisfies finite character.
	We will prove some other facts on our way.
	
	First note that $\sind$ satisfies the conditions of
	Proposition~\ref{PropExMS}.
	
	Then note that $A\sind_CB$ implies $A\ind_CB$.
	Hence $\ind$ also has local character, and therefore $\ind$ satisfies
	the conditions of Proposition~\ref{PropConsMS}.
	
	Using the two propositions we can show that $A\sind_CB$ implies $B\ind_CA$:
	If $\bar a\sind_CB$, there is a $\sind$-Morley sequence in~$\tp(\bar a/BC)$
	over~$C$. This sequence is also a $\ind$-Morley sequence in~$\tp(\bar a/BC)$
	over~$C$, hence $B\ind_C\bar a$.
	
	It follows that $A\sind_CB$ implies $B\sind_CA$:
	Suppose $A\sind_CB$ and $\hat A\supseteq A$.
	Since $\hat A\sind_{AC}AC$ by local character and base monotonicity,
	we can use extension to find $\hat A'\equiv_{AC}\hat A$ such that
	$\hat A'\sind_{AC}ABC$, hence $\hat A'\sind_{AC}B$ by monotonicity.
	Combining this with $A\sind_CB$, we get $\hat A'\sind_CB$ by transitivity.
	This implies $B\ind_C\hat A$.
	Thus $B\sind_CA$.
	
	Now we can prove that $\sind$ has finite character.
	Suppose $\bar a\sind_CB$ holds for all finite $\bar a\in A$.
	We need to show that $A\sind_CB$.
	So suppose $\hat B\supseteq B$.
	Since $A\sind_{BC}BC$ by local character and base monotonicity, we can obtain
	$A'\equiv_{BC}A$ such that $A'\sind_{BC}\hat B$
	using existence and monotonicity.
	By invariance, there is also $\hat B'\equiv_{BC}\hat B$
	such that $A\sind_{BC}\hat B'$.
	It suffices to show that $A\ind_C\hat B'$.
	For every finite $\bar a\in A$ we have $\bar a\sind_CB$ by assumption,
	and $\bar a\sind_{BC}\hat B'$ by $A\sind_{BC}\hat B'$ and monotonicity.
	Since $\sind$ is symmetric we can combine these results using
	transitivity on the right-hand side.
	Thus we get $\bar a\sind_C\hat B'$ for all finite $\bar a\in A$.
	Hence $\bar a\ind_C\hat B'$ for all finite $\bar a\in A$.
	Since $\ind$ has finite character, this implies $A\ind_C\hat B'$.
\end{proof}

\begin{notes}\index{dividing}\index{forking}
	The traditional way to define independence in stability theory is by first defining
	a notion of `dividing' and then deriving a notion of `forking'. Both steps
	are usually done with reference to individual formulas.
	While I have not seen Definition~\ref{DefSind} in this form before (except for special
	cases in~\cite{Hans Scheuermann: Unabhaengigkeitsrelationen}), it merely makes the
	step from dividing to forking explicit, while expressing it in semantic
	rather than syntactic terms. Thus it generalises the relation
	between classical dividing and classical forking and between thorn-dividing and thorn-forking.
	
	Theorem~\ref{ThmDivFork} was probably not stated
	in this generality before. One reason is the fact that the local
	character axiom\index{axiom!local character}\index{local character}
	in its strong form (no finiteness condition on $A$) is needed to get $A\sind_BB$ for
	all $A$, $B$ from local character and base monotonicity.
	The proof of finite character\index{axiom!finite character}\index{finite character}
	in Theorem~\ref{ThmDivFork} is a bit contrived.
	Note that in Chapter~\ref{Chapter2} (Lemma~\ref{LemmaDivForkSFC}) we will prove a shortcut for it
	in case $\ind$ actually has \emph{strong} finite character.
\end{notes}


\section{A theorem on Shelah-forking}

\begin{definition}\label{DefShelahForking}
	\index{015@$\dind$}\index{015@$\find$}
	\index{Shelah-dividing independence}
	\index{Shelah-forking independence}
	\index{independence!Shelah-dividing}
	\index{independence!Shelah-forking}
	The relation $\dind$ \emph{(Shelah-dividing independence)} is defined by
	\[ A\dind_CB \iff
		\bigg(\;
			\parbox{10.2cm}{\slshape
				\noindent for any sequence of $C$-indiscernibles $(\bar b_i)_{i<\omega}$ s.\,t.{}
				$\bar b_0\in BC$:\\
				$\exists\;A'\equ_{BC}A$ s.\,t.{} the sequence is $A'C$-indiscernible
			}
		\;\bigg).
	\]
	
	\noindent The relation $\find$ \emph{(Shelah-forking independence)} is defined by $\find = \ind[d$\ast$]$, i.e.:
	\[ A\find_CB \iff
		\Big(
			\textsl{ for all $\hat B\supseteq B$ there is $A'\equ_{BC}A$ s.t.\ $A'\dind_C\hat B$ } 
		\Big).
	\]
\end{definition}

\setlength{\unitlength}{1mm}
\begin{figure}\begin{picture}(50,27)(-30,0)
\put(0,10){$A$}\put(5,10){\line(4,-1){32}}
\put(40,0){$C$}
\put(30,20){$\bar b_0$}\put(39,5){\line(-1,2){6.5}}
\put(40,20){$\bar b_1$}\put(41.5,5){\line(0,1){13}}
\put(50,20){$\bar b_2$}\put(44,5){\line(1,2){6.5}}
\put(47,5){\line(1,1){10}}
\put(60,20){$\ldots$}
\put(47,21.5){\oval(40,10)}
\end{picture}\caption{\footnotesize An attempt to illustrate the definition of $\dind$.}\end{figure}

\noindent These definitions are equivalent to dividing and forking as they were originally defined
by Saharon Shelah to study stable theories (cf.~Exercise~\ref{ExcDividingForking} below).
The definition of $\dind$ can be seen as motivated by the following remark:

\begin{remark}\label{RemarkDindFinest}\ \\
	If $\ind$ is any independence relation,	then $A\dind_CB$ implies $A\ind_CB$.
\end{remark}

\begin{proof}
	Suppose $A\dind_C\bar b$.
	Let $(\bar b_i)_{i<\omega}$ be a $\ind$-Morley sequence for $\tp(\bar b/C)$.
	This exists by Proposition~\ref{PropExMS}, since $\bar b\ind_CC$.
	We may assume that $\bar b_0=\bar b$.
	Since $A\dind_C\bar b$ there is $A'\equiv_{\bar bC}A$ such that
	the sequence $(\bar b_i)_{i<\omega}$ is $A'C$-in\-dis\-cer\-ni\-ble.
	By Proposition~\ref{PropConsMS}, it now follows that $A\ind_C\bar b$.
\end{proof}

Of course this remark implies that $\dind=\find$ whenever $\find$
is an independence relation. But we still need $\find$ for technical reasons.

\begin{definition}\index{theory!simple}\index{simple theory}
	A complete consistent first-order theory $T$ is \emph{simple}
	if $\find$ is an independence relation for $T$.
\end{definition}

We will see in Chapter~\ref{Chapter2} that $T$ is simple iff $T\eq$ is simple.

\begin{lemma}\label{LemmaBasicDividing}
	The relation $\dind$ of Shelah-dividing independence always satisfies
	the first five axioms for independence relations and finite character.
	It also satisfies anti-reflexivity.
\end{lemma}

\begin{proof}
	\emph{Invariance} and \emph{monotonicity} are obvious.
	
	\emph{Base monotonicity:} Suppose $A\dind_CB$ and $C\subseteq C'\subseteq B$.
	Let $(\bar b_i)_{i<\omega}$ be a sequence of $C'$-indiscernibles
	with $\bar b_0\in B=BC$. Let $\bar c'$ be an enumeration of $C'$.
	Then also $\bar b_0\bar c'\in BC$, and the sequence
	$(\bar b_i\bar c')_{i<\omega}$ is also $C$-indiscernible.
	Hence there is $A'\equ_{\bar b_0\bar c'}A$ such that
	$(\bar b_i\bar c')_{i<\omega}$ is $A'C$-indiscernible.
	Thus $(\bar b_i)_{i<\omega}$ is $A'C'$-indiscernible.
	
	\emph{Transitivity:} Suppose $D\subseteq C\subseteq B$, $B\dind_CA$ and $C\dind_DA$.
	Let $(\bar a_i)_{i<\omega}$ be any sequence of $D$-indiscernibles
	with $\bar a_0\in AD$.
	
	By $C\dind_DA$ there is $C'\equ_{AD}C$ such that the sequence
	$(\bar a_i)_{i<\omega}$ is indiscernible over $C'$.
	Choose any set $B'$ such that $(B',C')\equ_{AD}(B,C)$.
	Then $B'\dind_{C'}A$ holds by invariance.
	Hence there is $B''\equ_{AC'}B'$ such that the sequence
	is $B''$-indiscern\-ible.
	And really, $B''\equ_{AD}B$.
	
	\emph{Normality:}
	Suppose $A\dind_CB$.
	Let $(\bar b_i)_{i<\omega}$ be a sequence of $C$-indis\-cernibles
	such that $\bar b_0\in BC$.
	By definition there is $A'\equiv_{BC}A$ such that the sequence
	is $A'C$-indiscernible.
	But then also $A'C\equiv_{BC}AC$.
	
	\emph{Finite character:}
	Let $\bar a$ be a possibly infinite tuple s.\,t.{} $\bar a\ndind_CB$.
	Let $p(\bar x)=\tp(\bar a/BC)$.
	Then there is a sequence $(\bar b_i)_{i<\omega}$
	with $\bar b_0\in BC$ such that the type extending
	$p(\bar x)$ and the theory of the big model with constants for $BC(\bar b_i)_{i<\omega}$ and expressing that
	$(\bar b_i)_{i<\omega}$ is $\bar aC$-indiscernible is inconsistent.
	By compactness, a finite sub-tuple $\bar a_0$ of $\bar a$ is sufficient for
	this, so $\bar a_0\ndind_CB$.
	
	For \emph{anti-reflexivity} suppose $a\not\in\acl B$. Then there is a
	$B$-indiscernible sequence $(a_i)_{i<\omega}$ of distinct elements, with $a_0=a$.
	This sequence witnesses that $a\ndind_Ba$.
\end{proof}

\begin{theorem}\index{independence relation!finest}\label{ThmCharSimple}
	A theory $T$ is simple if and only if $\find$ has local character.
	If $T$ is simple, $\find=\dind$, and this is the finest independence relation for $T$.
\end{theorem}

\begin{proof}
	For the equivalence, just apply Theorem~\ref{ThmDivFork} to $\dind$.
	Now suppose $T$ is simple.
	While $A\find_CB$ always implies $A\dind_CB$, the converse is true by Remark~\ref{RemarkDindFinest}.
	Hence $\find=\dind$.
	Since $\find$ is an independence relation and $\dind$ is finer than every
	independence relation by Remark~\ref{RemarkDindFinest}, $\find=\dind$ is the finest.
\end{proof}

\begin{exercises}
	\begin{exercise}\label{ExcDividingForking}(dividing and forking of formulas)
	
		A formula $\phi(\bar x;\bar b)$ \emph{divides} over a set~$C$
		if there is a finite number $k<\omega$ and a sequence $(\bar b_i)_{i<\omega}$
		such that $\bar b_i\equiv_C\bar b$ holds for all $i<\omega$
		and $\{\phi(\bar x;\bar b_i)\mid i<\omega\}$ is $k$-inconsistent.
		A formula \emph{forks} over~$C$ if it implies a finite
		disjunction of formulas that divide over~$C$.
		
		(i) $\bar a\dind_CB$ iff there is a tuple $\bar b\in BC$ and a formula $\phi(\bar x;\bar y)$
		without parameters such that $\models\phi(\bar a;\bar b)$ holds and
		$\phi(\bar x;\bar b)$ divides over~$C$.
		
		(ii) $\bar a\find_CB$ iff there is a tuple $\bar b\in BC$ and a formula $\phi(\bar x;\bar y)$
		without parameters such that $\models\phi(\bar a;\bar b)$ holds and
		$\phi(\bar x;\bar b)$ forks over~$C$.
		
		(iii) For simple $T$, a formula $\phi(\bar x;\bar b)$ forks over a set
		$C$ if and only if $\phi(\bar x;\bar b)$ divides over~$C$.
	\end{exercise}
	
	\begin{exercise}\label{ExcPropertiesDind}(additional properties of $\dind$)
		
		(i) Every sequence of $B$-indiscernibles is also indiscernible over~$\acl B$.
		
		(ii) $A\dind_CB$ implies $\acl(AC)\cap B\subseteq\acl C$.
		
		(iii) $A\dind_CB$ implies $A\dind_C\acl(BC)$.
		So $\dind$ always satisfies a weak variant of the extension axiom.
		
		(iv) If $A\dind_CB$ and $C\subseteq C'\subseteq\acl(BC)$,
		then $\acl(AC')\cap\acl(BC)=\acl C'$.
		Hence $A\dind_CB$ implies $A\mind_CB$, where $\mind$ is as defined in the next section.
	\end{exercise}
\end{exercises}

\begin{notes}
	Most of Lemma~\ref{LemmaBasicDividing} can be found in~\cite{Saharon Shelah: Classification Theory. 2nd ed.}.
	Transitivity of $\dind$ in the general case is implicit in~\cite{Byunghan Kim: Forking in simple unstable theories},
	and most of Theorem~\ref{ThmCharSimple} is also due to Byunghan Kim~\cite{Byunghan Kim: Forking in simple unstable theories}.
	I have not found the fact that Shelah-forking in a simple (or stable) theory is the
	finest independence relation stated outside~\cite{Hans Scheuermann: Unabhaengigkeitsrelationen}.
	I think this is due to the fact that independence relations without
	any additional conditions are not usually an object of study.
\end{notes}


\section{A theorem on thorn-forking}\label{SectionThornForking}

\begin{definition}\label{DefMThorn}
	\index{M-dividing independence}\index{thorn-forking independence}
	\index{independence!M-dividing}\index{independence!thorn-forking}
	\index{017@$\mind$}\index{017@$\thind$}
	The relation $\mind$ \emph{(M-dividing independence)} is defined by
	\[ A\mind_CB \iff
		\bigg(\;
			\parbox{6.3cm}{
				\noindent\slshape
				for any $C'$ s.\,t.\ $C\subseteq C'\subseteq\acl(BC)$:\\
				$\acl(AC')\cap\acl(BC') = \acl C'$
			}
		\;\bigg).
	\]
	The relation $\thind$ \emph{(thorn-forking independence)} is defined by $\thind = \ind[M$\ast$]$\,, i.e.:
	\[ A\thind_CB \iff
		\Big(
			\textsl{ for all $\hat B\supseteq B$ there is $A'\equ_{BC}A$ s.t.\ $A'\mind_C\hat B$ }
		\Big).
	\]
\end{definition}

\setlength{\unitlength}{1mm}
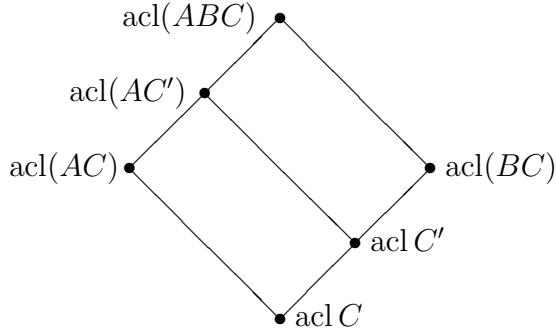
\begin{figure}
\begin{picture}(80,44)(-30,-1.5)
\newsavebox{\spot}
\savebox{\spot}(1,1)[bl]{\circle*{1.5}}
\put(15,20){\usebox{\spot}}		\put(-1,19){$\acl(AC)$}
\put(25,30){\usebox{\spot}}		\put(7,29){$\acl(AC')$}
\put(35,0){\usebox{\spot}}		\put(37,-1){$\acl C$}
\put(35,40){\usebox{\spot}}		\put(14,39){$\acl(ABC)$}
\put(45,10){\usebox{\spot}}		\put(47,9){$\acl C'$}
\put(55,20){\usebox{\spot}}		\put(57,19){$\acl(BC)$}
\put(15,20){\line(1,1){20}}
\put(35,0){\line(1,1){20}}
\put(35,40){\line(1,-1){20}}
\put(25,30){\line(1,-1){20}}
\put(15,20){\line(1,-1){20}}
\end{picture}\caption{\footnotesize A lattice diagram illustrating the definition of $\mind$
in the lattice of algebraically closed sets.
We have a map $\acl C'\mapsto \acl(AC')$ from the sublattice between
$\acl C$ and $\acl(BC)$ to the sublattice between $\acl(AC)$ and $\acl(ABC)$.
$A\mind_CB$ says that the map $D\mapsto D\cap\acl(BC)$ takes $D=\acl(AC')$
back to $\acl C'$.}\end{figure}

\noindent It is shown in Section~\ref{SectionThornDefinitions} that this
definition of $\thind$ agrees on $T\eq$ with thorn-forking independence as defined by
Alf Onshuus.
Here, however, we take the view that the definition of $\thind$ is motivated
by the following remark:

\begin{remark}\label{RemarkThindCoarsest}\ \\
	If $\ind$ is any \sir{}, then $A\ind_CB$ implies $A\thind_CB$.
\end{remark}

\begin{proof}
	Suppose $\ind$ is a \sir, $A\ind_CB$, and $\hat B\supseteq B$.
	We need to show that there is $A'\equ_{BC}A$ such that $A'\mind_C\hat B$.
	So choose $A'\equ_{BC}A$ such that $A'\ind_C\acl(\hat BC)$.
	For any $D$ satisfying $C\subseteq D\subseteq\acl(\hat BC)$ we get
	$A'\ind_D\acl(\hat BC)$ by base monotonicity of $\ind$.
	
	By extension and symmetry of $\ind$ there is a set
	$H\equ_{A'D}\acl(A'D)$ that satisfies $H\ind_D\acl(\hat BC)$.
	Clearly $H=\acl(A'D)$, so $\acl(A'D)\ind_D\acl(\hat BC)$.
	Now by anti-reflexivity of $\ind$,
	$\acl(A'D)\cap\acl(\hat BCD) \subseteq \acl D$,
	so $\acl(A'D)\cap\acl(\hat BCD) = \acl D$.
\end{proof}

By comparing Remarks~\ref{RemarkDindFinest} and~\ref{RemarkThindCoarsest}
one easily sees that $A\find_CB$ implies $A\thind_CB$, provided that a \sir{} exists.
Exercise~\ref{ExcPropertiesDind} showed that even without this
assumption a stronger statement is true:
$A\dind_CB$ always implies $A\mind_CB$,
hence $A\find_CB$ always implies $A\thind_CB$.

\begin{definition}
	\index{theory!rosy}\index{rosy theory}
	A complete consistent first-order theory $T$ is called \emph{rosy}
	if $\thind$ is an independence relation for $T\eq$.
\end{definition}

\begin{lemma}\label{LemmaBasicThorn}
	The relation $\mind$ of M-dividing independence always satisfies
	the first five axioms for independence relations and finite character.
	It also satisfies anti-reflexivity.
\end{lemma}

\begin{proof}
	\emph{Invariance, monotonicity, normality} and \emph{anti-reflexivity} are obvious.
	
	\emph{Base mo\-no\-to\-ni\-ci\-ty}:
	Suppose $A\mind_CB$ and $C\subseteq D\subseteq B$.
	Then for any $D'$ satisfying $D\subseteq D'\subseteq\acl(BD)$ we also have
	$C\subseteq D'\subseteq\acl(BC)$. So $A\mind_CB$ implies
	$\acl(AD')\cap\acl(BD') = \acl D'$. Hence $A\mind_{D'}B$.
	
	\emph{Transitivity:}
	Suppose $D\subseteq C\subseteq B$, $B\mind_{CD}A$ and $C\mind_DA$.
	Then for any $D'$ such that $D\subseteq D'\subseteq\acl(AD)$
	we can compute:
	\begin{align*}
		\acl(BD')\cap\acl(AD')&=\acl(BD')\cap\acl(ACD')\cap\acl(AD')\\
		                       &=\acl(CD')\cap\acl(AD')&&\!\!\big(\mbox{by $B\mind_CA$}\big)\\
		                       &=\acl D',&&\!\!\big(\mbox{by $C\mind_DA$}\big)
	\end{align*}
	so $B\mind_DA$ holds.
	
	\emph{Finite character:}
	Suppose $A\nmind_CB$. Let $C'$ be such that $C\subseteq C'\subseteq\acl(BC)$
	and $\acl(AC')\cap\acl(BC')\not\subseteq\acl C'$.
	Let $d\in(\acl(AC')\cap\acl(BC'))\setminus\acl C'$.
	Let $\bar a\in A$, finite, be such that $d\in\acl(\bar aC')$.
	Then clearly $\bar a\nmind_CB$.
\end{proof}

\begin{theorem}\label{ThmThind}\index{independence relation!coarsest strict}
	The relation $\thind$ of thorn-forking independence is a \sir{}
	if and only if it has local character, if and only if there is
	any \sir{} at all.
	If $\thind$ is a \sir{}, then it is the coarsest.
\end{theorem}

\begin{proof}
	To get the first equivalence, apply Theorem~\ref{ThmDivFork} to $\mind$.
	If $\ind$ is any \sir, then, since $\ind$ satisfies the local character axiom,
	so does $\thind$.
	If $\thind$ is a \sir{}, then it is the coarsest by Remark~\ref{RemarkThindCoarsest}.
\end{proof}

In particular, if there is any \sir{} for $T\eq$, then $T$ is rosy.
In a simple theory $T$, $\find$ is a \sir{} for $T\eq$ by Corollary~\ref{CorSimpleReduct}
below, so every simple theory is rosy.
Thus $\thind$ is the coarsest \sir{} on $T\eq$ while $\dind=\find$ is the finest.

Note that the assumptions of Theorem~\ref{ThmThind} do not imply $\mind=\thind$:

\begin{example}\label{ExForest2}\index{example!everywhere infinite forest}
	(Everywhere infinite forest, continued from Example~\ref{ExForest1})\\
	It follows from Theorem~\ref{ThmThind} that $\thind$ is also a \sir,
	and that $A\ind_CB\implies A\thind_CB$.
	It is straightforward to check that the converse is also true, so $\ind=\thind$.
	(Cf.~Exercise~\ref{ExcForest2}.)

	It is not hard to see that $\mind$ does not satisfy extension,
	so $\mind\neq\thind$:
	Let $a$ and $b$ be neighbours. Then $a\mind b$.
	However, there is no $c\equiv_ba$ such that $a\mind bc$:
	Either $c=a$, or $b$ lies between $a$ and $c$.
	In the first case, $a=c\in(\acl a\cap\acl(bc))\setminus\acl\emptyset$,
	so $a\nmind bc$.
	In the second case, $b\in(\acl(ac)\cap\acl(bc))\setminus\acl c$,
	so also $a\nmind bc$.
	Thus $\mind$ does not satisfy the extension axiom.
\end{example}

In some cases (most notably strongly minimal and o-minimal theories),
\th-forking as defined on the real elements of~$T$ is an important tool
for understanding the structure of models of~$T$. In these cases
\th-forking on~$T$ agrees with the restriction to~$T$ of \th-forking
in~$T\eq$. This is not the case in general, and the existence of a
\sir{} on the real elements of~$T$ \emph{per se} does not imply any
`structure' that is more than superficial:

\begin{example}\label{ExNotInteresting}
	\index{example!thorn-forking, in $T$ or in $T\eq$}
	(Thorn-forking must be computed in $T\eq$ in general)\\
	Let $T$ be any complete consistent theory in a relational language.
	Consider the following theory~$T'$:
	The language of~$T'$ is the language of~$T$ together with a new binary
	relation. The axioms of~$T'$ are the axioms of~$T$, but with equality
	replaced by the new relation, together with axioms saying that the
	new relation is an equivalence relation with infinite classes.
	Then $T'$ is a complete consistent theory satisfying $\acl A=A$
	for every set~$A$ of real elements.
	Hence the lattice of small algebraically closed sets is just the
	(modular) lattice of small subsets of the big model,
	and so the relation $A\ind_CB\iff A\cap B\subseteq C$ is a \sir{}
	for~$T$ (and agrees with $\thind$).
\end{example}

\begin{example}\label{ExTwoSirs}
	\index{example!two strict indep.{} relations}
	(A theory with two \sir s)\\
	Let $T$ be the theory of an equivalence relation with infinitely many
	infinite classes (i.\,e., \emph{all} classes are infinite),
	in the signature of a single relation~$E$.
	Then both $A\ind_CB\iff A\cap B\subseteq C$ (thorn-forking for $T$) and
	$A\eqind_CB\iff \acl\eq A\cap\acl\eq B\subseteq\acl\eq C$ (thorn-forking for $T\eq$)
	define \sir s on~$T$, but they are clearly not the same.
\end{example}

\begin{exercises}
	\begin{exercise}\label{ExcForest2}
		Check that $A\thind_CB$ implies $A\ind_CB$ in Example~\ref{ExForest2}.
	\end{exercise}
	\begin{exercise}\label{ExcNotInteresting}
		Check that $A\thind_CB\iff A\cap B\subseteq C$ is a \sir{} in Example~\ref{ExNotInteresting}.
	\end{exercise}
	\begin{exercise}\label{ExcTwoSirs}
		Check that $A\ind_CB\iff A\cap B\subseteq C$ and
		$A\eqind_CB\iff \acl\eq A\cap\acl\eq B\subseteq\acl\eq C$ are
		\sir s in Example~\ref{ExTwoSirs}.
	\end{exercise}
\end{exercises}

\begin{notes}
	The definition of thorn-forking independence $\thind$ via M-independence
	$\mind$ in Definition~\ref{DefMThorn} is new, but M-independence
	is from~\cite{Hans Scheuermann: Unabhaengigkeitsrelationen}.
	The original motiviation for the definition of $\mind$ was of course
	not thorn-forking. It was Exercise~\ref{ExcModular}, which shows
	that base monotonicity is the only problematic property for~$\aind$.
	If we try to force it, we get~$\mind$.
	The letter `M' was chosen because of the notation
	$\operatorname M(x,y)$ for modular pairs in lattices,\index{modular pair}\index{021@$\operatorname M(x,y)$}
	cf.{} Exercise~\ref{ExcMSymmetric}.
	
	Lemma~\ref{LemmaBasicThorn} is contained in~\cite[Lemmas 2.1.2 and 2.1.5]{Alf Onshuus: Properties and consequences of thorn-independence}.
	Alf Onshuus may have overlooked at first the fact
	(Theorem~\ref{ThmThind}, also Theorem~3 in~\cite{Hans Scheuermann: A note on rosy theories}) that
	thorn-forking is the coarsest \sir{} in every theory admitting one.
	
	Examples~\ref{ExNotInteresting} and~\ref{ExTwoSirs} are new.
\end{notes}


\section{Just a bunch of silly examples}
\label{SectionExamples}

In Sections 1 and 2 we found some relations that hold between the axioms
for independence relations and the existence and symmetry properties.
Our aim in this section is to show that we have actually found all of
them and, in particular, the axioms for independence relations are
independent. Most readers probably want to skip this section.

First we give two examples showing that invariance does not follow from the other axioms:

\begin{example}\textbf{(no invariance)}
	Let $\ind$ be any \sir. Let $F$ be a set such that
	$F\not\subseteq\acl\emptyset$.
	Define
	\[ A\pind_CB\iff A\ind_{CF}B. \]
	\noindent The relation $\pind$ satisfies all axioms for independence relations except invariance
	(but not anti-reflexivity).
	It also satisfies existence and symmetry.
\end{example}

\begin{example}\label{ExcNoInvariance}\textbf{(no invariance)}
	Consider the theory $T$ from Example~\ref{ExTwoSirs} with its two \sir s
	$\ind$ and $\eqind$\,. Now consider its reduct $T'$ that is just an infinite set
	without the equivalence relation. Take the big model of $T$ as a big model of $T'$.
	Then the relation $\eqind$\, satisfies all axioms for independence relations
	with respect to~$T'$ except invariance. It also satisfies existence and symmetry.
\end{example}

From this point on, we will only consider invariant relations in this section.

\begin{theorem}\label{ThmAxioms}
	Consider the following nine axioms that may hold for an invariant relation~$\ind$
	on the small sets of the big model of a complete theory:
	monotonicity, base monotonicity, transitivity, normality, extension, finite character,
	local character, existence, symmetry.
	The following relations hold between these axioms:
	\begin{enumerate}
		\item An invariant relation satisfying extension and symmetry also satisfies normality.
		\item An invariant relation satisfying extension and local character also satisfies existence.
		\item An invariant relation satisfying monotonicity, transitivity, normality, existence and symmetry
			also satisfies extension.
		\item An invariant relation satisfying monotonicity, base monotonicity, transitivity, normality,
			extension, finite character and local character also satisfies symmetry.
	\end{enumerate}
	This enumeration is complete: Every relation between these nine axioms that holds in general is a formal
	consequence of these four relations---with a grain of salt:
	The question whether monotonicity is needed in (3) is open.
\end{theorem}

\begin{proof}
	The relations (1)--(3) hold by Exercise~\ref{ExcExistence}.
	Relation (4) is Theorem~\ref{ThmSymmetry}.
	
	Completeness of the enumeration is proved by the following series of examples.
	Examples~\ref{ExNoMonotonicity}, \ref{ExNoBaseMonotonicity}, \ref{ExNoTransitivity}
	and~\ref{ExNoFC} show that monotonicity, base monotonicity, transitivity
	and finite character do not follow from any other axioms, respectively.
	Examples~\ref{ExNoNormalityNoSymmetry} and~\ref{ExNoNormalityNoExtension}
	show that normality does not follow from any set of other axioms that does
	not include at least extension and symmetry.
	Examples~\ref{ExNoExtensionNoExistence} and~\ref{ExNoLCNoExistence}
	show that existence does not follow from any set of other axioms that does
	not include at least extension and local character.
	Examples~\ref{ExNoExtensionNoTransitivity}, \ref{ExNoNormalityNoExtension},
	\ref{ExNoExtensionNoExistence} and~\ref{ExNoExtensionNoSymmetry}
	show that extension does not follow from any set of other axioms that does
	not include at least transitivity, normality, existence and symmetry.
	(Monotonicity is missing from this list.)
	Examples~\ref{ExNoMonotonicityNoSymmetry}, \ref{ExNoBaseMonotonicityNoSymmetry}, \ref{ExNoTransitivityNoSymmetry},
	\ref{ExNoNormalityNoSymmetry}, \ref{ExNoExtensionNoSymmetry}, \ref{ExNoFCNoSymmetry}
	and~\ref{ExNoLCNoSymmetry} show that symmetry does not follow from any set of other axioms that does
	not include at least monotonicity, base monotonicity, transitivity, normality,
	extension, finite character and local character.

	The author could not find an example satisfying all axioms except extension and monotonicity.
\end{proof}

\begin{example}\label{ExNoMonotonicity}\textbf{(no monotonicity)}\ \\
	Consider the theory from Example~\ref{ExTwoSirs} with its two \sir s
	$\ind$ and $\eqind$\,. Define
	\[ A\pind_CB\iff
		\bigg\{\;
			\parbox{5.5cm}{\sl
				$A\ind_CB$ if $A$ and $B$ are infinite\\[2mm]
				$A\eqind_CB$ if $A$ or $B$ is finite.
			}
	\]
	\noindent The relation $\pind$ satisfies all axioms for \sir s except monotonicity.
	It also satisfies existence and symmetry.
\end{example}

\begin{example}\label{ExNoBaseMonotonicity}\textbf{(no base monotonicity)}\ \\
	Consider the theory from Examples~\ref{ExForest1} and~\ref{ExForest2}.
	The relation $\aind$ from Exercise~\ref{ExcModular}
	satisfies all axioms for \sir s except base monotonicity.
	It also satisfies existence and symmetry.
\end{example}

\begin{example}\label{ExNoTransitivity}\textbf{(no transitivity)}\ \\
	Consider the theory from Example~\ref{ExTwoSirs} with its two \sir s
	$\ind$ and $\eqind$. Define
	\[ A\pind_CB\iff
		\bigg\{\;
			\parbox{4.5cm}{\sl
				$A\ind_CB$ if $C$ is infinite\\[2mm]
				$A\eqind_CB$ if $C$ is finite.
			}
	\]
	\noindent The relation $\pind$ satisfies all axioms for \sir s except transitivity.
	It also satisfies existence and symmetry.
	(By interchanging `finite' and `infinite' we would get another example without base monotonicity.)
\end{example}

\begin{example}\label{ExNoFC}\textbf{(no finite character)}\ \\
	In the theory of an infinite set with no structure,
	consider the relation
	\[ A\ind_CB \iff \big|(A\cap B)\setminus C\big|\leq\aleph_0. \]
	The relation $\ind$ satisfies all axioms for independence relations except finite character.
	It also satisfies existence and symmetry (but not anti-reflexivity).
	
	To get anti-reflexivity as well, consider the theory of an
	equivalence relation $E$ with infinitely many infinite classes.
	Let $\pi$ be the obvious projection from the standard sort
	to the imaginary sort of equivalence classes of~$E$.
	Define
	\[ A\pind_CB \iff A\cap B\subseteq C \textsl{ and }
		\big|\big(\pi(A)\cap\pi(B)\big)\setminus\pi(C)\big|\leq\aleph_0. \]
	\noindent The relation $\pind$ satisfies all axioms for \sir s except finite character.
	It also satisfies existence and symmetry.
\end{example}

\begin{example}\textbf{(no local character)}\label{ExNoLC}\ \\
	Consider the theory of the random graph, i.\,e.,
	the Fra\"issé limit of the finite undirected graphs,
	with the following relation:
	\[ A\ind_CB \iff A\cap B\subseteq C \textsl{ and there is no edge from $A\setminus C$ to $B\setminus C$.} \]
	\noindent The relation $\ind$ satisfies all axioms for \sir s except local character.
	It also satisfies existence and symmetry.
\end{example}

\begin{example}\label{ExNoNormalityNoSymmetry}\textbf{(no normality, no symmetry)}\ \\
	Consider the following relation:
	\[	A \ind_CB \iff \acl A\cap\acl(BC) \subseteq \acl C.
	\]
	It always satisfies all axioms for \sir s other than normality, and it also satisfies
	existence for every theory.
	But in the theory from Example~\ref{ExForest1} let $a\neq c$, and $b$ be points
	such that there is an edge from $a$ to $b$ and from $b$ to~$c$.
	Then $a\ind_cb$, $ac\nind_cb$ and $b\nind_ca$, so $\ind$ does not satisfy normality
	or symmetry.
\end{example}

\begin{example}\label{ExNoNormalityNoExtension}\textbf{(no normality, no extension)}\ \\
	Consider the following relation:
	\[	A \ind_CB \iff \acl A\cap\acl B \subseteq \acl C.
	\]
	It always satisfies all axioms for \sir s other than normality, and it also satisfies
	symmetry for every theory.
	But in the theory from Example~\ref{ExForest1} let $a\neq c$, and $b$ be points
	such that there is an edge from $a$ to $b$ and from $b$ to~$c$.
	Then $a\ind_cb$ but $ac\nind_cb$, so $\ind$ does not satisfy normality.
	It easily follows that $\ind$ does not satisfy extension either.
\end{example}

\begin{example}\label{ExNoLCNoExistence}\textbf{(no local character, no existence)}\ \\
	The empty ternary relation satisfies all axioms for \sir s except local character.
	It also satisfies symmetry, but not existence.
\end{example}

\begin{example}\label{ExNoExtensionNoExistence}\textbf{(no extension, no existence)}\ \\
	Let $\ind$ be a \sir{} for some theory $T$. Define $\ind[$\prime$]$ as follows:
	\[
		A\ind[$\prime$]_CB \quad\iff\quad
			\big( |C| \geq\aleph_0 \textsl{ and } A\ind_CB \big) \textsl{ or $A\subseteq C$ or $B\subseteq C$}.
	\]
	\noindent The relation $\ind[$\prime$]$ satisfies all axioms for \sir s except
	extension. It also satisfies symmetry, but not existence.
\end{example}

\begin{example}\label{ExNoExtensionNoTransitivity}\textbf{(no extension, no transitivity)}\ \\
	Given any \sir{} $\ind$, consider the following relation:
	\[
		A\pind_CB \iff A\ind_CB \textup{ or } \big(
			\card{A\setminus C} \leq 1\textsl{ and } \card{B\setminus C}\leq 1
		\big).
	\]
	It satisfies all axioms for \sir s except transitivity and extension,
	and it also satisfies existence and symmetry.
\end{example}

\begin{example}\label{ExNoExtensionNoSymmetry}\textbf{(no extension, no symmetry)}\ \\
	Consider the theory from Examples~\ref{ExForest1} and~\ref{ExForest2}.
	By Lemma~\ref{LemmaBasicThorn} the relation $\mind$ satisfies the axioms of \sir s except
	extension and local character. $\mind$ also satisfies local character
	and existence because $\thind$ does.
	We have already seen that $\mind$ does not satisfy extension.
	
	$\mind$ is not symmetric either: Suppose $b$ lies between $a$ and $c$.
	Then $a\nmind bc$ as we have just seen. But it is easy to see
	that $bc\mind a$.
\end{example}

\begin{example}\label{ExNoMonotonicityNoSymmetry}\textbf{(no monotonicity, no symmetry)}\ \\
	For any \sir{} $\ind$ consider the following relation:
	\[
		A\pind_CB \iff A\ind_CB \textsl{ or } \card{A\setminus\acl C}\geq 2.
	\]
	The relation $\pind$ satisfies all axioms for independence relations except monotonicity.
	It also satisfies extension, but not symmetry.
\end{example}

\begin{example}\label{ExNoBaseMonotonicityNoSymmetry}\textbf{(no base monotonicity, no symmetry)}\ \\
	Consider the theory of dense linear orders with the following relation:
	\[
		A\ind_CB\iff A\cap B\subseteq C \textsl{ or }
				\exists a\in A\exists c\in C: a<c.
	\]
	The relation $\ind$ satisfies all axioms for independence relations except
	base monotonicity. It also satisfies existence, but not symmetry.
\end{example}

\begin{example}\label{ExNoTransitivityNoSymmetry}\textbf{(no transitivity, no symmetry)}\ \\
	For any independence relation $\ind$ consider the following relation:
	\[
		A\pind_CB \iff A\ind_CB \textsl{ or } \card{A\setminus C}\leq 1.
	\]
	The relation $\pind$ satisfies all axioms for independence relations except
	transitivity. It also satisfies existence, but not symmetry.
\end{example}

\begin{example}\label{ExNoFCNoSymmetry}\textbf{(no finite character, no symmetry)}\ \\
	Given any \sir{} $\ind$, consider the following relation:
	\[
		A\pind_CB \iff \exists\textsl{ finite }B_0\subseteq B\textsl{ such that } A\ind_{B_0C}C.
	\]
	The relation $\pind$ always satisfies all axioms for \sir s except finite character,
	which it does not satisfy.
	It also satisfies existence, but not necessarily symmetry.
	
	Now let $T$ be the theory of $\omega$ cross-cutting equivalence relations $\epsilon_i$
	with infinitely many classes each. Let $\ind=\find$.
	Let $(a_i)_{i<\omega}$ and $b$ be such that $\models\epsilon_i(a_j,b)\iff i=j$.
	Then $(a_i)_{i<\omega}\pind_\emptyset b$ and $b\npind_\emptyset(a_i)_{i<\omega}$.
\end{example}

\begin{example}\label{ExNoLCNoSymmetry}\textbf{(no local character, no symmetry)}\ \\
	We will extend the theories $T_0$ and $T$ from Example~\ref{ExNoSir}.
	First we describe the signature of the respective extensions $T_0^*$ and $T^*$:
	It has the sorts $P$ (`points') and $E$ (`equivalence relations')
	as well as a new sort $\Gamma$ (`equivalence classes').
	The functions and relations of $T_0^*$ consist of the relation
	$p\sim_e q$ (for $p,q\in P$ and $e\in E$),
	a new relation written (slightly abusing notation)
	as $p/e=c$ for $p\in P$, $e\in E$ and $c\in\Gamma$,
	and a function $\epsilon:\Gamma\rightarrow E$.
	
	The axioms of $T_0^*$ include those of $T_0$, i.\,e.,
	$\sim_e$ is an equivalence relation for every $e\in E$.
	They also say that $\exists_{\leq 1}c(p/e=c)$,
	so it makes sense to regard $p/e$ as a partial function $P\times E\rightarrow\Gamma$
	which we will use informally in the following.
	The other axioms say $\epsilon(p/e)=e$ (if $p/e$ exists)
	and $p\sim_eq\leftrightarrow p/e=q/e$ (also if $p/e$ exists).
	
	Clearly every model of $T_0$ is also a model of $T_0^*$,
	and if we restrict a model of $T_0^*$ to the sorts $P$ and $E$
	we get a model of $T_0$.
	By the same arguments as for $T_0$ we can find an $\omega$-categorical
	theory $T^*$ with elimination of quantifiers which is the Fra\"{\i}ss\'e limit
	of the finite models of $T_0^*$. So $T^*$ extends both $T$ and $T_0^*$.
	
	For any subset $A$ of the big model of $T^*$ we write
	$P(A)=A\cap P$, $E(A)=(A\cap E)\cup\epsilon(A\cap\Gamma)$
	and $\Gamma(A)=(A\cap\Gamma)\cup\{p/e\mid p\in P(A), e\in E(A)\}$.
	It is not hard to check that
	$\acl A=\dcl A=P(A)\cup E(A)\cup \Gamma(A)$.
	It easily follows that
	\[
		A\mind_CB \iff 
			\left(\;
				\begin{aligned}
					P(A)\cap P(B)&\subseteq P(C)\textsl{ and }\\
					E(A)\cap E(B)&\subseteq E(C)\textsl{ and }\\
					P(A)/e \cap P(B)/e &\subseteq\Gamma(C) \textsl{ for all } e\in E(BC)
				\end{aligned}
			\;\right).
	\]
	Using this, it is not hard to check that $A\mind_CB\implies A\thind_CB$
	and that $A\mind_CB\implies A\dind_CB$, from which it easily follows
	that $\mind=\thind=\dind=\find$.
	Hence $\find$ satisfies all axioms for \sir s
	except local character (which would imply that there is
	a \sir{} for $T^*$). It also satisfies existence,
	though not symmetry.
\end{example}

\begin{exercises}
	\begin{exercise}\label{ExcNoExtensionNoExistence}\ \\
		Show that $\pind$ in Example~\ref{ExNoExtensionNoExistence}
		has the stated properties.
	\end{exercise}
\end{exercises}

\begin{notes}
	The author wishes to excuse for all the nonsensical examples in this section.
	Once he had found the first few, he could not resist the temptation to do a
	systematic search, the findings of which are now dumped on the reader.
\end{notes}
\chapter{Forking}
\label{Chapter2}

\noindent In this chapter we improve part of Chapter~\ref{Chapter1}
by exploring part of the local (i.\,e., relating to formulas),
or combinatorial, foundation of forking theory.
We will introduce the concept of inconsistency pairs. If $\Omega$
is a set of inconsistency pairs we get a relation $\omind$ such that
$\omsind$\; is a good candidate for being an independence relation.
In particular, $\omsind\;=\find$ or $\omsind\;=\thind$ for suitable
choices of $\Omega$. This will allow us to find out more about
$\find$ and $\thind$.

As in the previous chapter, the exposition is essentially self-contained.
I write $\bar a_{<k}$ for the tuple $\bar a_0\bar a_1\dots\bar a_{k-1}$
and $\bar a_{<\omega}$ for the sequence $(\bar a_i)_{i<\omega}$.
Tuples of variables or elements are often implicitly assumed to be
compatible: of the same lengths and with the same sorts at corresponding
positions.


\section{Strong finite character}

The purpose of this section is to give a foretaste of the improved results which
we will get in this chapter, while postponing the technicalities of the next two
sections as long as possible.
Here we show that we could have had slightly stronger results in Chapter~\ref{Chapter1} if we
had required the following stronger condition instead of finite character:

\begin{definition}\label{DefSFC}
	The \emph{strong finite character} condition is the following strong variant
	of the finite character axiom:
	\begin{description}
	\item[(strong finite character)]
		\index{axiom!strong finite character}
		\index{strong finite character}\ \\
		If $A\nind_CB$,
		then there are finite tuples $\bar a\in A$, $\bar b\in B$ and $\bar c\in C$
		and a formula $\phi(\bar x,\bar y,\bar z)$ without parameters
		such that
		\begin{itemize}
			\item $\models\phi(\bar a,\bar b,\bar c)$, and 
			\item $\bar a'\nind_C\bar b$ for all $\bar a'$ satisfying
				$\models\phi(\bar a',\bar b,\bar c)$.
		\end{itemize}
	\end{description}
\end{definition}

\noindent As we will see, this new condition is satisfied by all relations that are of interest
to us. It has two properties that make it more convenient than finite character.
The first is the following supplement to Lemma~\ref{LemmaDF}
that could have spared us the somewhat contrived proof of Theorem~\ref{ThmDivFork}.

\begin{lemma}\label{LemmaDivForkSFC}\ \\
	Let $\ind$ be a relation that satisfies invariance, monotonicity and the strong finite
	character condition. Then $\sind$ also satisfies the strong finite character condition.
\end{lemma}

\begin{proof}
	Suppose $\bar a\nsind_CB$
	($\bar a$ being a sequence of arbitrary length),
	and let this be witnessed by $\hat B\supseteq B$
	such that $\bar a'\nsind_C\hat B$ for all $\bar a'\equ_{BC}\bar a$.
	Let $\bar x$ be a sequence of the same length as $\bar a$, and let $p(\bar x)$ be the
	set of formulas over $\hat BC$ consisting of the negations of all those formulas
	$\phi_i(\bar x,\bar b_i,\bar c_i)$
	with parameters $\bar b_i\in\hat B$ and $\bar c_i\in C$ that have the property that
	$\bar a'\nind_C\bar b_i$ for all $\bar a'$ satisfying $\models\phi_i(\bar a',\bar b_i,\bar c_i)$.
	By choice of $\hat B$ and strong finite character of
	$\ind$, $p(\bar x)\cup\tp(\bar a/BC)$ is inconsistent.
	So by compactness there is a formula $\psi(\bar x,\bar b,\bar c)\in\tp(\bar a/BC)$
	such that $p(\bar x)\cup\{\psi(\bar x,\bar b,\bar c)\}$
	is inconsistent.
	
	Now suppose $\bar a'$ satisfies $\models\psi(\bar a',\bar b,\bar c)$.
	To finish our proof we claim that $\bar a'\nsind_C\bar b$.
	Otherwise there would be $\bar a^*\equ_{C\bar b}\bar a'$ such that $\bar a^*\ind_C\hat B$.
	But then $\models\psi(\bar a^*,\bar b,\bar c)$ would also hold.
	On the other hand, $\bar a^*$ would realise $p(\bar x)$, in contradiction to inconsistency
	of $p(\bar x)\cup\{\psi(\bar x,\bar b,\bar c)\}$.
\end{proof}

For the second advantage of strong finite character
recall that a type $p(\bar x)$ is called finitely satisfied\index{finitely satisfied} in a set~$C$ if
for every formula $\phi(\bar x,\bar b)\in p$ there is a tuple $\bar c\in C$
such that $\models\phi(\bar c,\bar b)$.

\begin{remark}\label{RemarkFSInd}
	Suppose $\ind$ satisfies monotonicity and strong finite character, and
	$\bar a,B,C$ are such that $C\ind_CB$ holds and $\tp(\bar a/BC)$ is finitely satisfied in~$C$.
	Then $\bar a\ind_CB$.
\end{remark}

\begin{proof}
	Suppose $\bar a\nind_CB$.
	Let $\phi(\bar x_0,\bar y,\bar z)$ and $\bar a_0\subseteq\bar a$, $\bar b\in B$, $\bar c\in C$ be as in the
	strong finite character condition.
	Since $\tp(\bar a/BC)$ is finitely satisfied in~$C$
	there is $\bar a'\in C$ such that $\models\phi(\bar a',\bar b,\bar c)$ holds.
	Hence $\bar a'\nind_CB$, hence $C\nind_CB$ by monotonicity.
\end{proof}

This is quite useful because of the following well-known fact:

\begin{remark}\label{RemarkFS}
	For any $\bar a$, $B$ there is a subset $C\subseteq\bar a$
	of size $\card C\leq\card T+\card B$ such that
	$\tp(\bar a/BC)$ is finitely satisfied in~$C$.
\end{remark}

\begin{proof}
	Let $C_0=\emptyset$. Given any set $C_n$ we construct $C_{n+1}\supseteq C_n$ as follows:
	For every formula $\phi(\bar x_0,\bar b)$, $\bar b\in BC_n$,
	that is satisfied by a finite subtuple $\bar a_0\subseteq\bar a$,
	we make sure that $C_{n+1}$ contains one such tuple $\bar a_0$.
	Clearly we can make sure that $\card{C_{n+1}}\leq\card T+\card{C_n}$.
	Now we can just take $C=\bigcup_{n<\omega}C_n$.
\end{proof}

Putting both results together it is easy to get the dual (left and right sides reversed)
of local character. Therefore we have:

\begin{theorem}\label{TheoremSFCExSymm}
	Suppose $\ind$ satisfies the first five axioms for independence relations as well
	as the strong finite character condition.
	Then $\ind$ is an independence relation if and only if $\ind$ satisfies
	existence and symmetry.
\end{theorem}

\begin{proof}
	First note that strong finite character implies finite character.
	We already know the forward direction, so we only need to prove
	extension and local character from existence and symmetry.
	
	Extension easily follows from transitivity, normality, existence and symmetry.
	For local character we can take $\kappa(B)=(\card T+\card B)^+$:
	Given $\bar a$ and $B$ there is $C\subseteq\bar a$ such that $\card C<\kappa(B)$
	and $\tp(\bar a/BC)$ is finitely satisfied in~$C$.
	Now $C\ind_CB$ holds by existence, so $\bar a\ind_CB$ by monotonicity
	and strong finite character.
\end{proof}

It follows that the relation in Example~\ref{ExNoLC} does not have strong finite character.

\setlength{\unitlength}{1mm}
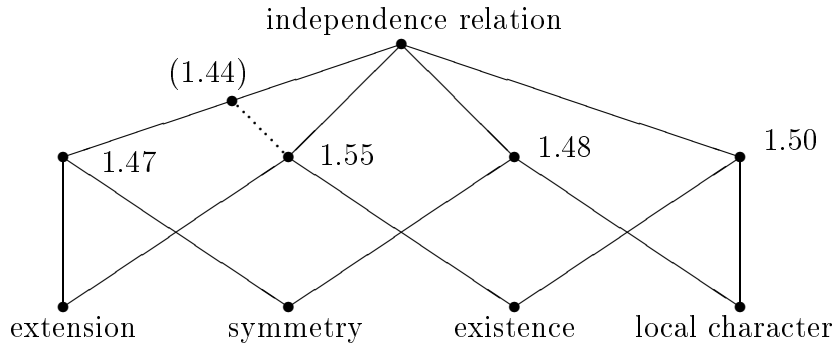
\begin{figure}\begin{picture}(50,46)(-23,0)
\savebox{\spot}(1,1)[bl]{\circle*{1.5}}
\put(0,5){\usebox{\spot}}	\put(-7,1){extension}
\put(30,5){\usebox{\spot}}	\put(22,1){symmetry}
\put(60,5){\usebox{\spot}}	\put(52,1){existence}
\put(90,5){\usebox{\spot}}	\put(76,1){local character}
\put(0,25){\usebox{\spot}}	\put(5,23){\ref{ExNoLCNoExistence}}
\put(30,25){\usebox{\spot}}	\put(34,24){\ref{ExNoLCNoSymmetry}}
\put(60,25){\usebox{\spot}}	\put(63,25){\ref{ExNoExtensionNoExistence}}
\put(90,25){\usebox{\spot}}	\put(93,26){\ref{ExNoExtensionNoSymmetry}}
\put(45,40){\usebox{\spot}}	\put(27,42){independence relation}
\put(0,5){\line(0,1){20}}
\put(90,5){\line(0,1){20}}
\put(0,5){\line(3,2){30}}
\put(30,5){\line(3,2){30}}
\put(60,5){\line(3,2){30}}
\put(0,25){\line(3,-2){30}}
\put(30,25){\line(3,-2){30}}
\put(60,25){\line(3,-2){30}}
\put(45,40){\line(-3,-1){45}}
\put(45,40){\line(3,-1){45}}
\put(45,40){\line(-1,-1){15}}
\put(45,40){\line(1,-1){15}}
\put(22.5,32.5){\usebox{\spot}}	\put(14,35){(\ref{ExNoLC})}
\multiput(22.5,32.5)(0.75,-0.75){10}{\circle*{0.5}}
\end{picture}\caption{\footnotesize Classification of relations satisfying the
first 5 axioms of independence relations and finite character,
according to which of 4 remaining properties hold.
For each point in the middle row of this lattice diagram there is an
example in Section~\ref{SectionExamples}.
These can be used to assemble examples for the bottom row.
If we also require strong finite character, the dotted line and
the point represented by Example~\ref{ExNoLC} disappear.}\label{FigureClassification}\end{figure}

\begin{exercises}
	\begin{exercise}\label{ExcNoExtensionNoSymmetry} ($\mind$ has strong finite character)
		
		The relation $\mind$ always satisfies the strong finite character condition.
	\end{exercise}
	
	\begin{exercise}(Figure \ref{FigureClassification})\label{ExcFigureClassification}
		
		Check that the examples mentioned in the middle row of Figure~\ref{FigureClassification}
		satisfy the strong finite character condition.
		For Example~\ref{ExNoExtensionNoExistence} assume that
		$\ind$ satisfies strong finite character.
	\end{exercise}

	\begin{exercise}(alternative definition for strong finite character)\label{ExcStrongFiniteCharacter}
	
		Suppose the relation $\ind$ satisfies invariance, monotonicity and extension.
		Prove that it satisfies the strong finite character condition 
		if and only if it satisfies the following condition:
		For any sequence of variables $\bar x$ and any sets $B$, $C$,
		the set $\big\{\tp(\bar a/BC) \;\big|\; \bar a\ind_CB\big\}$ is a closed
		subset of $\S^{\bar x}(BC)$.
	\end{exercise}
\end{exercises}

\begin{notes}
	The term `strong finite character' is probably new, but the property itself is
	essentially the anonymous axiom A.7 in~\cite{Michael Makkai: A survey of basic stability theory}.
	Exercise~\ref{ExcStrongFiniteCharacter}
	(read in conjunction with Theorem~\ref{TheoremSFCExSymm})
	shows that independence relations satisfying the strong finite character condition
	are precisely the relations considered in~\cite[Section~3]{Enrique Casanovas: Some remarks on indiscernible sequences}.
	None of the arguments in this section is new.
	
	I do not believe that \emph{every} independence relation has
	strong finite character, but I do not have a counter-example
	(cf.~Question~\ref{QEveryIndSFC}).
\end{notes}


\section{Local dividing}\index{dividing!local}\index{local dividing}

\begin{definition}\label{DefXi}\index{inconsistency witness}\index{inconsistency pair}\index{031@$\Xi$}
	The formula $\psi(\bar y_{<k})$ is called a \emph{$k$-inconsistency witness}
	for $\phi(\bar x;\bar y)$
	if the formula $\big(\bigwedge_{i<k}\phi(\bar x;\bar y_i)\big) \wedge \psi(\bar y_{<k})$
	is inconsistent.
	When the precise value of $k$ is immaterial we will omit it.
	We write 
	\[
		\Xi = \Big\{\big(\phi(\bar x;\bar y),\psi(\bar y_{<k})\big)
		\;\Big|\; \textsl{$\psi$ is a $k$-inconsistency witness for $\phi$; $k<\omega$} \Big\}
	\]
	for the set of all \emph{inconsistency pairs}.
\end{definition}

Note that in the preceding definition the free variables of $\phi(\bar x;\bar y)$
are partitioned in two blocks. The definition depends crucially on this partition.

A $k$-inconsistency witness $\psi(\bar y_{<k})$ for $\phi(\bar x;\bar y)$
`witnesses' $k$-inconsistency in the following way:
Suppose $(\bar b_i)_{i<\omega}$ is a sequence such that
$\models\psi(\bar b_{i_0},\dots,\bar b_{i_{k-1}})$ for any
$i_0<\dots<i_{k-1}<\omega$. Then the set $\{\phi(\bar x;\bar b_i)\mid i<\omega\}$
is $k$-inconsistent, i.\,e., there is no tuple $\bar a$ satisfying $k$
formulas from the set simultaneously.

\begin{definition}\index{dividing $(\phi,\psi)@$(\phi,\psi)$-dividing}
	A formula $\phi(\bar x;\bar b)$ \  \emph{$(\phi,\psi)$-divides} over a set $C$
	if $(\phi,\psi)\in\Xi$ and there is a sequence $\bar b_{<\omega}$
	such that
	\begin{itemize}
		\item each $\bar b_i$ realises $\tp(\bar b/C)$, and
		\item $\models\psi(\bar b_{i_0},\ldots,\bar b_{i_{k-1}})$ holds for all
			$i_0<\ldots <i_{k-1}<\omega$.
	\end{itemize}
	We say that $\bar b_{<\omega}$ \emph{witnesses} that $\phi(\bar x;\bar b)$ \ 
	$(\phi,\psi)$-divides over~$C$.
	
	A partial type $p(\bar x)$ \emph{$(\phi,\psi)$-divides} over a set $C$
	if it contains a formula $\phi(\bar x;\bar b)\in p(\bar x)$
	that $(\phi,\psi)$-divides over~$C$.
\end{definition}

Note that when $\phi(\bar x;\bar b)$ \  $(\phi,\psi)$-divides over a set~$C$,
then there is a sequence $\bar b_{<\omega}$ witnessing this with $\bar b_0=\bar b$.
Also note that $\phi(\bar x;\bar b)$ also $(\phi,\psi)$-divides over every
subset of~$C$.

\begin{definition}\label{DefOmega}\index{dividing $\Omega$@$\Omega$-dividing}\index{032@$\omind$}
	Let $\Omega\subseteq\Xi$ be a subset of $\Xi$ that is \emph{closed under
	variable substitution} in the following sense:
	If $\big(\phi(\bar x;\bar y),(\psi(\bar y_{<k})\big)\in\Omega$
	and $\bar u$, $\bar v$, $\bar v_{<k}$ are appropriate tuples of
	variables---possibly with repetitions,
	but $\bar u$, $\bar v$ and the tuples $\bar v_i$ being pairwise disjoint from each other---then
	$\big(\phi(\bar u,\bar v),\psi(\bar v_{<k})\big)\in\Omega$.
	
	We say that a partial type $p(\bar x)$ \  \emph{$\Omega$-divides} over a
	set $C$ if it $(\phi,\psi)$-divides over $C$ for some $(\phi,\psi)\in\Omega$.
	We define a relation $\omind$ as follows:
	\[ A\omind_CB \iff
		\textsl{there is no $\bar a\in A$ such that $\tp(\bar a/BC)$ $\Omega$-divides over $C$.}
	\]
\end{definition}

Note that $\Xi$ itself is closed under variable substitution.

\begin{proposition}\label{PropBasicOmind}
	Suppose $\Omega\subseteq\Xi$ is closed under variable substitution.
	Then $\omind$ satisfies the following axioms for independence relations:
	invariance, monotonicity, base monotonicity and finite character.
	In fact, $\omind$ has strong finite character.
	Moreover, $A\omind_BB$ and $A\omind_AB$ for any sets $A$ and $B$.
\end{proposition}

\begin{proof}
	Invariance and monotonicity are obvious.
	
	\emph{Base monotonicity:}
	Suppose $A\nomind_CB$ and $D\subseteq C\subseteq B$. It suffices to show that $A\nomind_DB$.
	There is $\bar a\in A$, $(\phi,\psi)\in\Omega$ and $\bar b\in B$
	such that $\phi(\bar x;\bar b)\in\tp(\bar a/B)$ and
	$\phi(\bar x;\bar b)$ \  $(\phi,\psi)$-divides over~$C$.
	It is immediate from the definition of $(\phi,\psi)$-dividing that
	$\phi(\bar x;\bar b)$ also $(\phi,\psi)$-divides over~$D$.
	So $A\nomind_DB$ does in fact hold.
	
	\emph{Strong finite character:}
	Suppose $A\nomind_CB$.
	Let $\bar a\in A$ be such that $\tp(\bar a/BC)$ \  $\Omega$-divides over~$C$.
	So there is $(\phi,\psi)\in\Omega$ and $\bar b\in BC$
	such that $\phi(\bar x;\bar b)\in\tp(\bar a/BC)$ and
	$\phi(\bar x;\bar b)$ \  $(\phi,\psi)$-divides over~$C$.
	Hence for every $\bar a'$ satisfying $\models\phi(\bar a';\bar b)$,
	$\tp(\bar a'/\bar bC)$ also $(\phi,\psi)$-divides over~$C$,
	so $\bar a'\nomind_C\bar b$.	
	
	For the first `moreover' statement, suppose $(\phi,\psi)\in\Xi$,
	$\models\phi(\bar a;\bar b)$ for some tuples $\bar a\in A$ and $\bar b\in B$,
	and $\phi(\bar x;\bar b)$ \  $(\phi,\psi)$-divides over~$B$.
	This would be witnessed by a sequence $\bar b_{<\omega}$ of
	tuples realising $\tp(\bar b/B)$, so $\bar b_i=\bar b$. But then
	$\models\big(\bigwedge_{i<k}\phi(\bar a;\bar b_i)\big)\wedge\psi(\bar b_{<k})$,
	contradicting the assumption that $(\phi,\psi)$ is a $k$-inconsistency
	witness.
	
	For the second `moreover' statement, suppose $(\phi,\psi)\in\Xi$,
	$\models\phi(\bar a;\bar b)$ for some tuples $\bar a\in A$ and $\bar b\in B$,
	and $\phi(\bar x;\bar b)$ \  $(\phi,\psi)$-divides over~$A$.
	This would be witnessed by a sequence $\bar b_{<\omega}$ of
	tuples realising $\tp(\bar b/A)$. But then again
	$\models\big(\bigwedge_{i<k}\psi(\bar a;\bar b_i)\big)\wedge\psi(\bar b_{<k})$,
	contradicting the assumption that $(\phi,\psi)$ is a $k$-inconsistency
	witness.
\end{proof}

So the missing axioms are transitivity, normality, extension and local character.
Heuristically speaking,
from our experience in Chapter~\ref{Chapter1} we can say that extension is probably no problem
since we can fix it by passing to $\omsind$\,, while we could not expect
local character to hold in general. Only the fact that we cannot prove
transitivity and normality is a bit annoying (since both are among the
first five axioms), so let us check that the problem is real:

\begin{example}\label{ExOmindNotTransitiveNotNormal}
	Let $T$ be the theory of an infinite set in the empty signature.
	Let $\Omega$ consist of all inconsistency pairs of the form
	$\big(\phi(xx';yy'),\psi(y_0y'_0,y_1y'_1)\big)$, where
	$\phi(xx';yy')\equiv x\neq x'\wedge x=y\wedge x'=y'$
	and $\psi(y_0y'_0,y_1y'_1)\equiv y_0\neq y_1 \wedge y'_0\neq y'_1$.
	
	It is not hard to see that
	$A\omind_CB\iff \big|(A\cap B)\setminus C\big| \leq 1$,
	that $\omind$ satisfies extension, local character,
	existence and symmetry, and that $\omind=\omsind$\,.
	But $\omind$ does not satisfy transitivity:
	Suppose $b\neq c$.
	Then $bc\omind_cbc$ and $c\omind_{\emptyset}bc$,
	but $bc\nomind_{\emptyset}bc$.
	Since $b\omind_\emptyset bc$ and $bc\nomind_\emptyset bc$,
	$\omind$ does not satisfy normality, either.
\end{example}

\begin{example}\label{ExTypeForkBase}
	Let $T$ be a theory in which there is a type that forks over its domain
	in the sense of Shelah-forking. Two examples of this phenomenon were given
	by Saharon Shelah in~\cite[Exercise III.1.3]{Saharon Shelah: Classification Theory. 2nd ed.}.
	In Proposition~\ref{PropDindXind} below we will see that $\xind=\dind$.
	From this it easily follows that $\xind$ does
	not satisfy extension or existence.
	Moreover, it follows from Theorem~\ref{ThmCharacterisationsSimple} below
	that $\xind=\dind$ does not satisfy local character or symmetry, either.
\end{example}

The choice of $\Omega$ in Example~\ref{ExOmindNotTransitiveNotNormal}
was of course perverse. Here are two natural conditions that we may require so
that $\omind$ makes sense:

\begin{definition} Suppose $\Omega\subseteq\Xi$ is closed under variable substitution.

	We say that $\Omega$ is \emph{transitive}\index{transitive $\Omega$}\index{transitivity}\index{axiom!transitivity}
	if the following holds:\\
	Suppose $(\phi(\bar y;\bar x),\psi(\bar x_{<k})\in\Omega$ and $C\omind_D\bar a$, where $D\subseteq C$.
	If $\phi(\bar y;\bar a)$ \  $(\phi,\psi)$-divides over~$D$
	then $\phi(\bar y;\bar a)$ \  $(\phi,\psi)$-divides over~$C$.
	
	We say that $\Omega$ is \emph{normal}\index{normal $\Omega$}\index{normality}
	if the following holds:\\
	If $\big(\phi(\bar x,\bar z;\bar y), \psi(\bar y_{<k})\big)\in\Omega$,
	then also $\big(\phi(\bar x;\bar z,\bar y), \psi'(\bar y_{<k},\bar z_{<k})\big)\in\Omega$,
	where $\psi'(\bar y_{<k},\bar z_{<k}) \equiv \psi(\bar y_{<k})\wedge (z_0=z_1=\dots =z_{k-1})$.
\end{definition}

\begin{proposition}\ \\[-7mm]\mbox{}\hfill\mbox{}
	\begin{enumerate}
	\item If $\Omega$ is transitive, then $\omind$ satisfies the transitivity axiom.
	\item If $\Omega$ is normal, then $\omind$ satisfies the normality axiom.
	\end{enumerate}
\end{proposition}

\begin{proof}
	(1) Suppose $\Omega$ is transitive, $D\subseteq C\subseteq B$, $C\omind_DA$
	and $B\nomind_DA$. Then there is $(\phi,\psi)\in\Omega$, $\bar b\in B$ and $\bar a\in AD$
	such that $\models\phi(\bar b;\bar a)$ and $\phi(\bar y;\bar a)$ \  
	$(\phi,\psi)$-divides over~$D$. Since $C\ind_D\bar a$ it follows that
	$\phi(\bar y;\bar a)$ also $(\phi,\psi)$-divides over~$C$.
	Hence $B\nomind_CA$.
	
	(2) Now suppose instead that $\Omega$ is normal and $AC\nomind_CB$.
	Then there is $\big(\phi(\bar x,\bar z;\bar y),\psi(\bar x_{<k},\bar z_{<k})\big)\in\Omega$,
	$\bar a\in A$, $\bar b\in B$ and $\bar c\in C$
	such that $\models\phi(\bar a,\bar c;\bar b)$
	and $\phi(\bar x,\bar z;\bar b)$ \  $(\phi,\psi)$-divides over~$C$.
	Let $\psi'$ be as in the definition of normality for $\Omega$.
	Then $\big(\phi(\bar x;\bar z,\bar b),\psi'\big)\in\Omega$
	and $\phi(\bar x;\bar c,\bar b)$ clearly
	$(\phi,\psi')$-divides over~$C$.
	Hence $A\nomind_CB$.
\end{proof}

\begin{exercises}
	\begin{exercise}\label{ExcForkStar}\index{forking $\Delta$@$\Delta$-forking}($\Delta$-forking)\\
		For $\Omega\subseteq\Xi$ let
		$\Omega\restrict\bar x =
			\big\{ \big(\phi(\bar x;\bar y),\psi(\bar y_{<k})\big) \;\big|\; k<\omega,\;(\phi,\psi)\in\Omega \big\}$
		be the set of those tuples $(\phi,\psi)\in\Omega$ for which the left block
		of variables of $\phi$ is~$\bar x$.
		For any $\Delta\subseteq\Xi\restrict\bar x$, we say that a partial type $p(\bar x)$ \emph{$\Delta$-forks}
		over a set~$C$ if there are $n<\omega$,
		$\big(\phi^i(\bar x;\bar y^i),\psi^i(\bar y^i_{<k_i})\big)\in\Delta$ for $i<n$,
		and tuples $\bar b^0,\ldots,\bar b^{n-1}$ such that
		$p(\bar x)\vdash\bigvee_{i<n}\phi^i(\bar x;\bar b^i)$
		and $\phi^i(\bar x;\bar b^i)$ $(\phi^i,\psi^i)$-divides over~$C$ for each $i<n$.
		
		Suppose $\Omega\subseteq\Xi$ is closed under variable substitution.
		Show that $\bar a\omsind_CB$ iff $\tp(\bar a/BC)$ does not $\Delta$-fork
		over~$C$ for any finite $\Delta\subseteq\Omega$.
	\end{exercise}
	
	\begin{exercise}\label{ExcDividing}(more on dividing and forking of formulas, cf.~Exercise~\ref{ExcDividingForking})\\
		(i) Given any formula $\phi(\bar x;\bar b)$, show that
		$\phi(\bar x;\bar b)$ divides over~$C$ iff 
		$\phi(\bar x;\bar b)$ \ $(\phi,\psi)$-divides over~$C$
		for some formula $\psi(\bar x_{<k})$ that is a
		$k$-inconsistency witness for~$\phi(\bar x;\bar y)$.\\
		(ii) Show that $\phi(\bar x;\bar b)$ forks over~$C$ iff 
		$\phi(\bar x;\bar b)$ \  $\Delta$-forks over~$C$
		for some set $\Delta\subseteq\Xi\restrict\bar x$.
	\end{exercise}
	
	\begin{exercise}\label{ExcOmindNotTransitiveNotNormal}\ \\
		Check the claim that $A\omind_CB\iff \big|(A\cap B)\setminus C\big| \leq 1$
		in Example~\ref{ExOmindNotTransitiveNotNormal}.
	\end{exercise}
\end{exercises}

\begin{notes}
	This section was derived from a small part of~\cite{Itay Ben-Yaacov: Simplicity in compact abstract theories}
	by localising and simplifying it. (Note that there is no array-dividing here.)
	Example~\ref{ExOmindNotTransitiveNotNormal} is new, Example~\ref{ExTypeForkBase} is from Saharon Shelah.
	The definitions of transitivity and normality of $\Omega$ seem to be new.
	They were found by the author when he looked for a condition that makes
	Lemma~\ref{LemmaRankSymmetry2} true and holds for both $\Xi$ and $\Xim$.
	($\Xim$ is defined below in Definition~\ref{DefinitionXim}.)
	
	It should be noted that both $\Xi$ and $\Xim$ are of the form
	$\Omega(\Psi)=\{(\phi,\psi)\in\Xi\mid\psi\in\Psi\}$.
	Yet it seems to be necessary to localise dividing in $\phi$
	as well as in $\psi$ in order to get a good theory of local rank.
\end{notes}


\section{Dividing patterns}

For a tuple $\bar x$ of variables we write $\Xi\restrict\bar x$
for the set of inconsistency pairs $(\phi,\psi)\in\Xi$ such that
$\phi$ has the form $\phi(\bar x;\bar y)$, with arbitrary $\bar y$.

\begin{definition}\label{DefDivPat}\index{dividing pattern}
	Let $p(\bar x)$ be a partial type over~$C$ and $I$ a linearly ordered set.
	
	An $I$-sequence $\xi=((\phi^i,\psi^i))_{i\in I}\in\Xi^I$ is a \emph{dividing pattern
	for~$p(\bar x)$ (over~$C$)} if there is
	an $I$-sequence $(\bar b^i)_{i\in I}$ that \emph{realises~$\xi$ over~$C$,}
	i.\,e.:
	\begin{itemize}
	\item $p(\bar x)\cup\{\phi^i(\bar x;\bar b^i)\mid i\in I\}$ (makes sense and) is consistent, and
	\item each formula $\phi^i(\bar x;\bar b^i)$ \  $(\phi^i,\psi^i)$-divides
		over~$C\bar b^{<i}$.
	\end{itemize}
	If $\Delta\subseteq\Xi\restrict\bar x$ and $\xi\in\Delta^I$ we may call $\xi$ a $\Delta$-dividing pattern.
\end{definition}

Vaguely speaking, dividing patterns measure how many dividing extensions
a type has. Under certain conditions an extension of a type that admits exactly the same
dividing patterns will be shown not to divide.

If $I$ is a
linearly ordered set and $i\in I$ we will temporarily write $<i$ and $\leq i$
for the initial sequences $\{j\in I\mid j<i\}$ and
$\{j\in I\mid j\leq i\}$, respectively.

\newcommand{\divseq}{\operatorname{divpat}}
\begin{theorem}\label{TheoremTypeDivpat}\index{035@divpat}
	Let $p(\bar x)$ be a partial type, definable over a set~$C$.
	An $I$-sequence $\xi=\big(\big(\phi^i(\bar x;\bar y^i), \psi^i(\bar y^i_{<k_i})\big)\big)_{i\in I}\in\Xi^I$
	is a dividing pattern for~$p(\bar x)$ over~$C$ iff the following type
	$\divseq_p^\xi\big((\bar x_\alpha)_{\alpha\in\omega^I},(\bar y_\alpha)_{\alpha\in\omega^{\leq i}, i\in I}\big)$
	is consistent:
	\begin{multline*}
		\bigcup_{\alpha\in\omega^I}p(\bar x_\alpha)
		\quad\cup\quad \big\{\phi^i(\bar x_\alpha;\bar y_{\alpha\restrict\leq i})
			\;\big|\; i\in I,\alpha\in\omega^I\big\}\\
		\cup\quad \big\{\;\psi^i(\bar y_{\alpha_0},\dots,\bar y_{\alpha_{k_i-1}})
			\quad\big|\quad i\in I,\; \alpha_0,\dots,\alpha_{k_i-1}\in\omega^{\leq i},\\
			(\alpha_0\restrict<i)=\dots=(\alpha_{k_i-1}\restrict<i), 
			\textsl{ and }\alpha_0(i)<\dots<\alpha_{k_i-1}(i) \;\big\}.
	\end{multline*}
\end{theorem}

\noindent Before proving this theorem let us try to understand what it says.
Without understanding the structure of the type $\divseq_p^\xi$
it is at least easy to see that it does not mention the set~$C$.
That's why the qualification `over~$C$' is in parentheses in
Definition~\ref{DefDivPat}. The next easy observation is that
the surrounding theory is not involved in the definition of $\divseq_p^\xi$.
Hence if $p$ and $\xi$ make sense in a reduct $T'$ of $T$,
then $\xi$ is a dividing pattern for $p$ in the context of~$T$
iff it is one in the context of~$T'$. We will use this to
prove that simplicity and rosiness are preserved in reducts.

For understanding the structure of $\divseq_p^\xi$ it is perhaps
best to imagine this type partially realised by tuples
$(\bar b_\alpha)_{\alpha\in\omega^{\leq i}, i\in I}$.
These tuples form a non-standard tree, and the last part of the conjunction
requires that the tuples $\bar b_\alpha$
of level~$i$ (i.e.: $\alpha\in\omega^{\leq i}$)
that define the same non-standard path $\alpha\restrict<i$
through the tree are related by the inconsistency witness $\psi^i$.
The type
$\divseq_p^\xi\big((\bar x_\alpha)_{\alpha\in\omega^I},(\bar b_\alpha)_{\alpha\in\omega^{\leq i}, i\in I}\big)$
then merely expresses that for every branch $\alpha\in\omega^I$ of this tree the set
$\big\{\phi^i(\bar x,\bar b_{\alpha\restrict\leq i}) \;\big|\; i\in I\big\}$
is consistent with~$p(\bar x)$.

With this tree structure in mind it is easy to see that, by
compactness, the property of being a dividing pattern has finite character:
$\divseq_p^\xi$ is consistent iff $\divseq_p^{\xi\restrict J}$
is consistent for every finite $J\subseteq I$.

The tree structure of $\divseq_p^\xi$ already suggests a proof strategy.

\begin{proof}
	We will prove the equivalence of the following statements:
	\begin{enumerate}
		\item $\xi$ is a dividing pattern for $p$ over $C$.
		\item $\divseq_p^\xi$ is consistent.
		\item The type
			\begin{multline*}
				\divseq'{}_p^\xi\big((\bar x_\alpha)_{\alpha\in\omega^I},(\bar y_\alpha)_{\alpha\in\omega^{\leq i}, i\in I}\big)\\
				\begin{aligned}
				= \;&\divseq_p^\xi\big((\bar x_\alpha)_{\alpha\in\omega^I},(\bar y_\alpha)_{\alpha\in\omega^{\leq i}, i\in I}\big)\\
				\cup \;&\big\{ \bar y_\alpha\equiv_{C\{\bar y_{\alpha\restrict\leq j} \mid j < i\}}\bar y_{\alpha'} 
				\;\big|\;
					i\in I,\;
					\alpha,\alpha'\in\omega^{\leq i},\;
					(\alpha\restrict<i)=(\alpha'\restrict<i)
				\}
				\end{aligned}
			\end{multline*}
			is consistent.
	\end{enumerate}
	
	\noindent We first prove that (3) implies (1):
	Let the tuples $(\bar b_\alpha)_{\alpha\in\omega^{\leq i}, i\in I}$
	be a partial realisation of $\divseq'{}_p^\xi$.
	For $i\in I$ write $\zeta^i$ for the unique function $\zeta^i\in\{0\}^{<i}$,
	and for $m<\omega$ write $\zeta^i{}\concat(m)$ for the extension of $\zeta^i$
	that maps $i$ to~$m$.
	Then for every $i\in I$ the sequence $(\bar b_{\zeta^i{}\concat(m)})_{m<\omega}$
	witnesses that $\bar b_{\zeta^i{}\concat(0)}$ \  $(\phi^i,\psi^i)$-divides
	over~$C\{\bar b_{\zeta^j{}\concat(0)}\mid j<i\}$.
	Hence the $I$-sequence $(\bar b_{\zeta^i{}\concat(0)})_{i\in I}$
	realises $\xi$ over~$C$.
	
	Next we observe that we need only prove that (1) implies (2)
	and that (2) implies (3) in case $I$ is finite. The general
	case then follows by compactness. Thus we can use induction
	on the size of~$I$.
	
	The case $I=\emptyset$ is trivial:
	$()\in\Xi^0$ is a dividing pattern for~$p$ over~$C$
	iff $p$ is consistent,
	and we have $\divseq_p^{()}=\divseq'{}_p^{()}=p(\bar x_{()})$.
	
	Now suppose the implications $(1)\implies (2)\implies (3)$ hold for $I$,
	and we are given
	$((\phi^s(\bar x;\bar y),\psi^s(\bar y_{<k})))\concat\xi\in\Xi^{\{s\}\cup I}$,
	where $s\not\in I$ is less than every element of~$I$.
	It is not hard to see that $(1)\implies (2)\implies (3)$
	for $((\phi^s,\psi^s))\concat\xi$,
	using the following three easy facts:
	
	(i) \  $((\phi^s,\psi^s))\concat\xi$ is a dividing pattern for $p$ over $C$
	iff there is a tuple $\bar b$ such that
	$\phi^s(\bar x;\bar b)$ \  $(\phi^s,\psi^s)$-divides over~$C$,
	and $\xi$ is a dividing pattern for $p(\bar x)\cup\phi^s(\bar x;\bar b)$.
	
	(ii) $\divseq_p^{((\phi^s,\psi^s))\concat\xi}$ is consistent iff
	there is a sequence $(\bar b_m)_{m<\omega}$ such that
	$\models\psi^s(\bar b_{m_0},\dots,\bar b_{m_{k-1}})$
	for any $m_0<\dots<m_{k-1}<\omega$ and
	the type $\divseq_{p(\bar x)\cup\phi^s(\bar x;\bar b_m)}^\xi$
	is consistent for every $m<\omega$.
	
	(iii) $\divseq'{}_p^{((\phi^s,\psi^s))\concat\xi}$ is consistent iff
	there is a sequence $(\bar b_m)_{m<\omega}$ such that
	$\models\psi^s(\bar b_{m_0},\dots,\bar b_{m_{k-1}})$
	for any $m_0<\dots<m_{k-1}<\omega$,
	$\bar b_m\equiv_C\bar b_0$ for all $m<\omega$, and
	the type $\divseq'{}_{p(\bar x)\cup\phi^s(\bar x;\bar b_m)}^\xi$
	is consistent for every $m<\omega$.
	(Thus the sequence $(\bar b_m)_{m<\omega}$ witnesses that
	$\bar b_0$ \  $(\phi^s,\psi^s)$-divides over~$C$.)
\end{proof}

\begin{exercises}
	\begin{exercise}(tree property)\label{ExcTreeProperty}\index{tree property}
	
		A formula $\phi(\bar x;\bar y)$ \emph{has the tree-property (of order~$k$)}
		if there is a tree of tuples $(\bar b_\alpha)_{\alpha\in\omega^{<\omega}}$
		such that for every limit point $\alpha\in\omega^\omega$ the branch
		$\{\phi(\bar x;\bar b_{\alpha\restrict n})\mid n<\omega\}$,
		is consistent, and at every node $\alpha\in\omega^{<\omega}$ the set of
		successors $\{\phi(\bar x;\bar b_{\alpha\concat(i)})\mid i<\omega\}$
		is $k$-inconsistent (i.\,e., every subset with $k$ elements is inconsistent).
		
		Show that formula $\phi(\bar x;\bar y)$ has the tree-property of order $k$
		if and only if there is a $k$-inconsistency witness $\psi(\bar y_{<k})$
		for $\phi$ such that $\D_{\phi,\psi}(\emptyset) = \infty$.
		Here $\D_{\phi,\psi}$ is $\D_\Delta$ as defined in Definition~\ref{DefinitionDRank} below
		for the case $\Delta=\{(\phi,\psi)\}$.
	\end{exercise}
\end{exercises}

\begin{notes}
	The only new things in this section are the \emph{term}
	`dividing pattern' and the idea of admitting
	arbitrary linear orders in order to allow a uniform treatment
	of dividing patterns and the tree property.
	
	Dividing patterns appear in~\cite{Itay Ben-Yaacov: Thickness and a categoric view of type-space functors} in the
	following guise:
	Let $\alpha$ be an ordinal and $I=\alpha^{\operatorname{opp}}$, i.\,e.,
	$\alpha$ with the opposite order.
	Then $\xi\in\Xi^I$ is a dividing pattern for $p$ iff
	$\xi\in\D(p,\Xi)$ in the notation of~\cite[Definition 1.8]{Itay Ben-Yaacov: Thickness and a categoric view of type-space functors}.
	The idea that $\xi$ being a dividing pattern can be expressed
	by a partial type is also from~\cite{Itay Ben-Yaacov: Thickness and a categoric view of type-space functors}.
	
	A realisation of a dividing pattern $\xi\in\Xi^\alpha$ is also the same
	thing as a dividing chain as defined in~\cite{Enrique Casanovas: The number of types in simple theories}.
\end{notes}


\section{Local rank and symmetry}

Suppose $\Delta\subseteq\Xi\restrict\bar x$ is finite.
If there are arbitrarily long finite
$\Delta$-dividing patterns for~$p$, then there is an inconsistency
pair $(\phi,\psi)\in\Delta$ such that there are arbitrarily long
finite $(\phi,\psi)$-dividing patterns for~$p$. It follows that
$(\phi,\psi)^I$ is a $\Delta$-dividing pattern for~$p$ for every
linearly ordered set~$I$. Therefore the following definition
makes sense:

\begin{definition}\label{DefinitionDRank}\index{037@$\D_\Delta$}\index{local rank}\index{rank!local}
	Let $p(\bar x)$ be a partial type and $\Delta\subseteq\Xi\restrict\bar x$
	a finite set of inconsistency pairs	$(\phi(\bar x;\bar y_{\phi}),\psi)$.
	Then $\D_\Delta(p)\in\omega\cup\{\infty\}$ is $\infty$ if $p$ has
	$\Delta$-dividing patterns of arbitrary order type,
	or otherwise the greatest number $n<\omega$ such that $\Delta$-dividing patterns
	of length~$n$ exist for~$p$.
\end{definition}

If $p=\tp(\bar a/B)$ we abbreviate $\D_\Delta(\bar a/B)=\D_\Delta(p)$.

\begin{remark}
	For any $\bar a$, $B$, $C$ and finite $\Delta\subseteq\Xi\restrict\bar x$:
	\[ \D_\Delta(\bar a/BC)\leq\D_\Delta(\bar a/C). \]
\end{remark}

\begin{proof}
	Let $p=\tp(\bar a/BC)$.
	If $\xi\in\Delta^n$ is a dividing pattern for~$p$, then
	$\divseq_p^\xi$ is consistent. Hence $\divseq_{p\restrict C}^\xi$
	is consistent, so $\xi$ is a dividing pattern for~$p\restrict C$.
\end{proof}

For the rest of this section we fix a set $\Omega\subset\Xi$
that is closed under variable substitution.

\begin{lemma}\label{LemmaRankSymmetry1}
	Suppose $\D_\Delta(\bar a/BC)=\D_\Delta(\bar a/C)<\infty$
	for all finite $\Delta\subseteq\Omega\restrict\bar x$.
	Then $\bar a\omsind_CB$.
\end{lemma}

\begin{proof}
	Towards a contradiction, suppose $\bar a\nomsind_CB$.
	Then there is a set $\hat B\supset B$ such that
	$\bar a'\nomind_C\hat B$ holds for every $\bar a'$ realising $p(\bar x)=\tp(\bar a/BC)$.
	Hence the set \[ p(\bar x) \cup \big\{\neg\phi(\bar x;\bar b) \;\big|\;
	\textsl{$(\phi,\psi)\in\Delta$, and $\phi(\bar x;\bar b)$ \  $(\phi,\psi)$-divides over~$C$} \big\}\]
	is inconsistent.
	By compactness there are inconsistency pairs $(\phi^i,\psi^i)\in\Delta$
	and tuples $\bar b^i$ such that
	$p(\bar x)\vdash\bigvee_{i<k}\phi^i(\bar x;\bar b^i)$
	and $\phi^i(\bar x;\bar b^i)$ \  $(\phi^i,\psi^i)$-divides over~$C$.
	
	Let $\xi$ be a $\Delta$-dividing pattern for~$p$ of maximal length
	$\card\xi=\D_\Delta(p)$, realised over $C\bar b^0\bar b^1\dots\bar b^{k-1}$
	by, say, $(\bar b_j)_{j<\card\xi}$.
	Let $\bar a'$ realise $p(\bar x)\cup\{\phi_j(\bar x;\bar b_j)\mid j<\card\xi\}$.
	Then $\models\phi^i(\bar a';\bar b^i)$ for an index $i<k$.
	Hence $(\phi^i,\psi^i)\concat\xi$ is a $\Delta$-dividing pattern for~$p$,
	realised by $\bar b^i{}\concat(\bar b_j)_{j<\card\xi}$.
	This contradicts maximality of $\card\xi$.
\end{proof}

Having shown a connection between $\omind$ and the $\D_\Delta$-ranks under
a combinatorial condition, we now show a sort of converse under
a geometric condition.

\begin{lemma}\label{LemmaRankSymmetry2}
	Suppose $\Omega$ is transitive and normal and $B\omsind_C\bar a$.\\
	Then for every finite $\Delta\subseteq\Omega\restrict\bar x$
	we have $\D_\Delta(\bar a/BC)=\D_\Delta(\bar a/C)$.
\end{lemma}

\begin{proof}
	Since $\D_\Delta(\bar a/BC)\leq\D_\Delta(\bar a/C)$
	holds anyway we need only prove that $\D_\Delta(\bar a/C)\geq n$
	implies $\D_\Delta(\bar a/BC)\geq n$.
	By definition of $\D_\Delta$ there is a $\Delta$-dividing pattern
	$\xi=((\phi^i,\psi^i))_{i<n}\in\Delta^n$ for $\tp(\bar a/C)$,
	and this is witnessed by tuples $(\bar b^i)_{i<n}$
	such that $\models\phi^i(\bar a;\bar b^i)$
	and $\phi^i(\bar x;\bar b^i)$ \  $(\phi^i,\psi^i)$-divides
	over $C\bar b^{<i}$ for all $i<n$.
	Since $B\omsind_C\bar a$
	we may assume that $B\omind_C\bar a\bar b^{<n}$.
	Hence $BC\bar b^{<i}\omind_{C\bar b^{<i}}\bar b_i$ by base monotonicity and normality
	for all $i<n$.
	Now since $\phi^i(\bar x;\bar b^i)$ \  $(\phi^i,\psi^i)$-divides
	over~$C\bar b^{<i}$ we get by transitivity that
	$\phi^i(\bar x;\bar b^i)$ \  $(\phi^i,\psi^i)$-divides
	over~$BC\bar b^{<i}$ as well.
	Therefore $\bar b^{<n}$ also witnesses that $\xi$ is a
	dividing pattern for $\tp(\bar a/BC)$,
	so $\D(\bar a/BC)\geq n$.
\end{proof}

\begin{theorem}\label{TheoremRankSymmetry}
	Suppose $\Omega$ is transitive and normal, and
	$\D_\Delta(\emptyset)<\infty$ for all finite $\Delta\subseteq\Omega\restrict\bar x$.
	Then the following conditions are equivalent:
	\begin{enumerate}
		\item $\bar a\omsind_CB$.
		\item $\D_\Delta(\bar a/BC)=\D_\Delta(\bar a/C)$ for all finite
			$\Delta\subseteq\Omega\restrict\bar x$.
		\item $B\omsind_C\bar a$.
	\end{enumerate}
\end{theorem}

\begin{proof}
	(3) implies (2) by Lemma~\ref{LemmaRankSymmetry2},
	and (2) implies (1) by Lemma~\ref{LemmaRankSymmetry1}.
	Hence $\omsind$\,\, is symmetric, so (1) implies (3).
\end{proof}

\begin{notes}
	Theorem~\ref{TheoremRankSymmetry} was proved by Byunghan Kim in
	the case $\Omega=\Xi$, cf.~\cite[Theorem~5.1]{Byunghan Kim: Forking in simple unstable theories},
	and by Alf Onshuus in the case $\Omega=\Xim$, cf.~\cite[Theorem~4.1.3]{Alf Onshuus: Properties and consequences of thorn-independence}.
	(To be pedantic, both Kim and Onshuus prove their results for ranks that are very similar
	to, but not exactly the same as, our ranks $\D_\Delta$.)
\end{notes}


\section{A better theorem on local forking}

Like in the previous section we fix a set $\Omega\subseteq\Xi$
which is closed under variable substitution.

\begin{lemma}\label{LemOmindStarRaw}
	The following statements are equivalent:
	\begin{enumerate}
		\item
			$\omsind$\, satisfies the local character axiom.
		\item
			$\omind$ satisfies the local character axiom.
		\item
			$\D_{\phi,\psi}(\emptyset)<\infty$ for every $(\phi,\psi)\in\Omega$.
	\end{enumerate}
\end{lemma}

\begin{proof}
	(1) implies (2):
	This follows from $A\omsind_CB\implies A\omind_CB$.
	
	(2) implies (3):
	Suppose $\omind$ has local character with a constant~$\kappa$,
	but $\D_{\phi,\psi}(\emptyset)=\infty$.
	We may assume that $\kappa$ is regular.
	$(\phi,\psi)^\kappa$ is a dividing pattern for the empty type,
	so it has a realisation $(\bar b_i)_{i<\kappa}$ over~$\emptyset$.
	Let $\bar a$ realise $\{\phi(\bar x;\bar b_i)\mid i<\kappa\}$.
	By local character there is a subset $C\subseteq\bar b_{<\kappa}$
	such that $\bar a\omind_C\bar b_{<\kappa}$ and $\card C<\kappa$.
	Since $\kappa$ is regular, there is $\alpha<\kappa$ such that
	$C\subseteq\bar b_{<\alpha}$.
	Hence $\bar a\omind_{\bar b_{<\alpha}}\bar b_{<\kappa}$,
	a contradiction to the fact that $\models\phi(\bar a;\bar b_\alpha)$ holds
	and $\phi(\bar x;\bar b_\alpha)$ \  $(\phi,\psi)$-divides over $\bar b_{<\alpha}$.
	
	(3) implies (1):
	Note that $\D_\Delta(p(\bar x))<\omega$ for all finite $\Delta\subseteq\Omega\restrict\bar x$
	and partial types $p(\bar x)$.
	We will prove local character for $\omsind$\; with $\kappa=\card T^+$.
	So suppose we have a type $p(\bar x)=\tp(\bar a/B)$ with finite~$\bar a$.
	
	For every finite $\Delta\subseteq\Omega\restrict\bar x$ we can find a finite subset
	$C_\Delta\subseteq B$ such that $\D_\Delta(p\restrict C_\Delta)=\D_\Delta(p)$:
	For each $\Delta$-dividing pattern $\xi$ of length $\card\xi = \D_\Delta(p)+1$
	(there are only finitely many) the type $\divseq^\xi_p$
	is inconsistent, so there is a finite subset $C_\xi\subseteq B$ such that
	$\divseq^\xi_{p\restrict C_\xi}$ is still inconsistent.
	If $C_\Delta$ is the union of these sets $C_\xi$,
	then clearly $C_\Delta$ is a finite set such that
	$\D_\Delta(p\restrict C_\Delta)=\D_\Delta(p)$.
	
	Now let $C$ be the union of these sets $C_\Delta$ for all finite $\Delta\subseteq\Omega\restrict\bar x$.
	Then $\card C\leq\card T$, so $\card C<\kappa$.
	Moreover, $\D_\Delta(p\restrict C)=\D_\Delta(p\restrict C_\Delta)=\D_\Delta(p)$
	for all finite $\Delta\subseteq\Omega\restrict\bar x$.
	Hence $\bar a\omsind_CB$ by Lemma~\ref{LemmaRankSymmetry1}.
\end{proof}

Now we are prepared for the following improved version of Theorem~\ref{ThmDivFork}
for the case $\ind=\omind$:

\begin{theorem}\label{TheoremOmindStar}
	Suppose $\omind$ satisfies the transitivity and normality axioms.\\
	Then $\omsind$\, is an independence relation if and only if
	the following, equivalent, conditions are satisfied:
	\begin{enumerate}
		\item
			$\omsind$\, satisfies the local character axiom.
		\item
			$\omind$ satisfies the local character axiom.
		\item
			$\D_{\phi,\psi}(\emptyset)<\infty$ for every $(\phi,\psi)\in\Omega$.
		\item
			$A\omsind_CB$ implies $B\omsind_CA$.
		\item
			$A\omsind_CB$ implies $B\omind_CA$.
	\end{enumerate}
\end{theorem}

\begin{proof}
	By Proposition~\ref{PropBasicOmind} (and since we have assumed transitivity
	and normality), $\omind$ satisfies the first five axioms for independence
	relations and strong finite character.
	Hence by Lemmas~\ref{LemmaDF} and~\ref{LemmaDivForkSFC},
	the relation $\omsind$\, satisfies all axioms except local character.
	Therefore it is an independence relation if and only if (1) holds.
	
	Conditions (1) to (3) are equivalent by Lemma~\ref{LemOmindStarRaw}.
	If $\omsind$\, is an independence relation, then (4) holds
	by Theorem~\ref{TheoremRankSymmetry}.
	(4) implies (5) because $B\omsind_CA$ implies $B\omind_CA$.
	
	(5) implies (2):
	We choose $\kappa(A)=(\card T+\card A)^+$ for every set~$A$.
	Given sets $A$ and~$B$, let $\bar b$ be an enumeration of~$B$.
	By Remark~\ref{RemarkFS} there is a subset $C\subseteq B$
	such that $\card C<\kappa(A)$ and $\tp(\bar b/AC)$ is finitely
	realised in~$C$.
	Hence by Remark~\ref{RemarkFSInd} we have $B\find_CA$,
	so $A\dind_CB$ by~(5).
\end{proof}

Note that if $\Omega$ satisfies the slightly stronger conditions
of transitivity and normality (and it will do so in our applications),
then Theorem~\ref{TheoremRankSymmetry}
gives us additional information on $\omsind$.

\begin{exercises}
	\begin{exercise}\label{ExcOmindStar} (symmetry of $\omind$)
	
		Suppose $\omind$ satisfies transitivity, normality and symmetry.
		Then $\omind$ is an independence relation, and $\omind=\omsind$\,.
	\end{exercise}
\end{exercises}

\begin{notes}
	This section extends some standard results to a more general context.
\end{notes}


\section{A better theorem on Shelah-forking}

In order to avoid the use of exercises in the main text
(Exercises~\ref{ExcDividingForking} and~\ref{ExcDividing} in this case),
we give a direct proof of the following result:

\begin{proposition}\label{PropDindXind}
	$\dind=\xind$.
\end{proposition}

\begin{proof}
	First suppose $A\nind[$\Xi$]_CB$, so there is $\bar a\in A$ and
	$\big(\phi(\bar x;\bar y),\psi(\bar y_{<\omega})\big)\in\Xi$
	such that $\tp(\bar a/BC)$ \  $(\phi,\psi)$-divides over~$C$.
	Hence there is $\bar b\in BC$ such that $\phi(\bar x;\bar b)\in\tp(\bar a/BC)$
	and $\phi(\bar x;\bar b)$ \  $(\phi,\psi)$-divides over~$C$.
	Let $\bar b_{<\omega}$ be a sequence witnessing this, i.e.,
	each $\bar b_i$ realises $\tp(\bar b/C)$ and
	$\models\psi(\bar b_{i_0},\ldots\bar b_{i_{k-1}})$
	holds for all $i_0<\ldots<i_{k-1}<\omega$.
	We can extend this sequence and extract a sequence of
	$C$-indiscernibles from it, so we may assume that
	$\bar b_{<\omega}$ is $C$-indiscernible.
	Moreover, we may assume that $\bar b_0=\bar b$.
	
	Towards a contradiction, suppose $A\dind_CB$. Then there would be
	$\bar a'\equiv_{BC}\bar a$ such that the sequence is
	$\bar a'C$-indiscernible. But then
	$\models\big(\bigwedge_{i<k}\phi(\bar a';\bar b_i)\big)\wedge\psi(\bar b_{<k})$,
	contradicting the fact that $(\phi,\psi)$ is a $k$-inconsistency witness.
	
	For the converse, suppose $A\ndind_CB$, so there is a sequence
	$\bar b_{<\omega}$ of $C$-in\-dis\-cer\-nibles with $\bar b_0\in BC$
	that cannot be $A'C$-indiscernible for any $A'\equiv_{BC}A$.
	By compactness this must be due to a formula $\phi(\bar a;\bar b_0)$
	($\bar a\in A)$ and a formula $\psi(\bar b_{<k})$ such that
	$\big(\bigwedge_{i<k}\phi(\bar x;\bar b_0)\big)\wedge\psi(\bar b_{<k})$
	is inconsistent.
	So $\big(\phi(\bar x;\bar y),\psi(\bar y_{<k})\big)$
	is a $k$-inconsistency witness and $\tp(\bar a/BC)$
	\  $(\phi,\psi)$-divides over~$C$.
	Hence $A\nind[$\Xi$]_CB$ as well.
\end{proof}

\begin{theorem}\label{ThmCharacterisationsSimple}\index{theory!simple}\index{simple theory}
	A complete consistent theory $T$ is simple if and only if the following, equivalent,
	conditions are satisfied:
	\begin{enumerate}
		\item $\find$ satisfies the local character axiom.
		\item $\dind$ satisfies the local character axiom.
		\item $\D_{\phi,\psi}(\emptyset)<\omega$ for each $(\phi,\psi)\in\Xi$.
		\item $\find$ is symmetric.
		\item $A\find_CB$ implies $B\dind_CA$.
		\item $\dind$ is symmetric.
	\end{enumerate}
	Moreover, in a simple theory $\find=\dind$ is the finest
	independence relation.
\end{theorem}

\begin{proof}
	Simplicity is equivalent to (1) by Theorem~\ref{ThmCharSimple}.
	For the equivalence of (1)--(5) note that
	$\find=\indb[${\scriptscriptstyle\textnormal d}{\scriptstyle*}$]=\indb[${\scriptscriptstyle\Xi}{\scriptstyle*}$]$.
	Then apply Theorem~\ref{TheoremOmindStar}.
	(6) clearly implies (5). The converse holds since (6) implies $\find=\dind$
	by Theorem~\ref{ThmCharSimple}. The `moreover' statement is also
	from Theorem~\ref{ThmCharSimple}.
\end{proof}

\begin{corollary}\label{CorSimpleReduct}
	Every reduct of a simple theory is simple.\\
	$T$ is simple if and only if $T\eq$ is simple.
\end{corollary}

\begin{proof}
	We first show that a reduct of a simple theory is simple.
	We already know that a theory $T$ is simple iff
	$\D_{\phi,\psi}(\emptyset)<\infty$ holds for every
	inconsistency pair $(\phi,\psi)$.
	Let $T'$ be a reduct of $T$.
	For formulas $\phi$ and $\psi$ in the signature of $T'$,
	$(\phi,\psi)$ is an inconsistency pair for $T'$ if and only if
	it is an inconsistency pair for $T$.
	Moreover, $\D_{\phi,\psi}(\emptyset)$ is the maximal $n<\omega$
	such that for the unique $\xi\in\{(\phi,\psi)\}^n$ the type
	$\divseq_p^\xi$ from Theorem~\ref{TheoremTypeDivpat} is consistent.
	Since this type is independent of the ambient theory,
	it does not matter whether we evaluate $\D_{\phi,\psi}$
	in $T$ or in $T'$.
	Thus if $T$ is simple then so is $T'$.
	
	One consequence is that if $T\eq$ is simple, then so is its
	reduct~$T$.
	We now show the converse.
	So suppose $T$ is simples and $(\phi,\psi)$ is a $k$-inconsistency pair for~$T\eq$.
	We may assume that as much as possible is coded in a single
	imaginary variable, so $\phi\equiv\phi(x;y)$ and $\psi\equiv\psi(y_{<k})$.
	The sorts of $x$ and $y$ correspond to definable equivalence
	relations $\epsilon_x$ and $\epsilon_y$.
	Now consider $\phi'(\bar x;\bar y)\equiv\phi(\bar x/\epsilon_x;\bar y/\epsilon_y)$
	and $\psi'(\bar y_{<k})\equiv\psi(\bar y_0/\epsilon_y,\ldots,\bar y_{k-1}/\epsilon_y)$.
	$\phi'$~and $\psi'$ can be expressed in $T$, and
	$(\phi',\psi')$ is a $k$-inconsistency pair for~$T$.
	Clearly $\D_{\phi,\psi}(\emptyset) = \D_{\phi',\psi'}(\emptyset)$,
	so $T\eq$ also satisfies condition (3) of Theorem~\ref{ThmCharacterisationsSimple}.
\end{proof}

The following remark shows that we can also apply Theorem~\ref{TheoremRankSymmetry} to get an alternative
characterisation of~$\find$ in a simple theory.

\begin{remark}
	$\Xi$ is transitive and normal.
\end{remark}

\begin{proof}
	For transitivity of $\Xi$ suppose $C\xind_D\bar a_0$, $D\subseteq C$ and
	$\phi(\bar y;\bar a_0)$ \  $(\phi,\psi)$-divides over~$D$,
	witnessed by $\bar a_{<\omega}$.
	By compactness and Fact~\ref{FactExtraction}
	we may assume that $\bar a_{<\omega}$ is $D$-indiscernible.
	Since $C\dind_D\bar a_0$ we may assume that $\bar a_{<\omega}$ is in fact $C$-indiscernible.
	Thus $\bar a_{<\omega}$ witnesses that $\phi(\bar y;\bar a_0)$ \ 
	$(\phi,\psi)$-divides over~$C$.
	
	For normality of $\Xi$ just observe that
	$\psi(\bar y_{<k})$ is a $k$-inconsistency witness for
	$\phi(\bar x,\bar z;\bar y)$ iff
	$\psi(\bar y_{<k})\wedge z_0=z_1=\dots z_{k-1}$ is
	a $k$-inconsistency witness for $\phi(\bar x;\bar z,\bar y)$.
\end{proof}

\begin{notes}
	Apart from the style of presentation, nothing is new in this section.
\end{notes}


\section{A better theorem on thorn-forking}\label{SectThorn2}

\begin{definition}\index{043@$\Xim$}\label{DefinitionXim}
	We define the following subset of $\Xi$:
	\begin{multline*}
		\Xim = \Big\{
			\big(
				\phi(\bar x;u\bar v),\psi((u\bar v)_{<k})
			\big) \in\Xi
		\;\Big|\\
				\psi((u\bar v)_{<k})
					\quad\equiv\quad
						{\bigwedge}_{i<j<k} (u_i\neq u_j \wedge \bar v_i=\bar v_j)
		\Big\}.
	\end{multline*}
\end{definition}

\noindent Note that if $\psi((u\bar v)_{<k})$ (as in the definition) is a $k$-inconsistency witness for
$\phi(\bar x;u\bar v)$, then whenever $\phi(\bar a;g\bar h)$ holds, $g$ must be algebraic over~$\bar a\bar h$.

\begin{proposition}\label{PropXim}
	Some properties of $\ind[$\Xim$]$\;:
	\begin{enumerate}
		\item
			$\ind[$\Xim$]$\; has the following characterisation:
			\[ A\ind[$\Xim$]_C\;B \iff
				\bigg(\;
					\parbox{8cm}{
						\noindent $\acl(AD)\cap B\subseteq\acl D$\\
						\mbox{}\ \ \ \ \ \textsl{for every set $D$ such that $C\subseteq D\subseteq BC$}
					}
				\;\bigg).
			\]
		\item
			$\Xim$ is transitive and normal.
		\item
			$A\mind_CB$ implies $A\ind[$\Xim$]_CB$.
		\item
			$\thind = \indb[${\scriptscriptstyle\textnormal M}\scriptstyle*$]\; = \indb[${\scriptscriptstyle\Xim}\scriptstyle*$]$\ \;.
	\end{enumerate}
\end{proposition}

\begin{proof}
	(1) Suppose there is a set $D$ such that $C\subseteq D\subseteq BC$ and
	$\acl(\bar aD)\cap B\nsubseteq\acl D$.
	So there is an element $e\in\acl(\bar aD)\cap B\setminus\acl D$.
	Let $\alpha(u,\bar a,\bar d)$ with $\bar d\in D$ be an algebraic formula
	realised by~$e$.
	
	Then for some $k<\omega$, $\models\phi(\bar a;e\bar d)$ holds,
	where $\phi(\bar x;u\bar v)
		\equiv\alpha(u,\bar x,\bar v)\wedge\exists_{<k}u'\alpha(u',\bar x,\bar v)$.
	We set $\psi((u\bar v)_{<k}) \equiv \bigwedge_{i<j<k}(u_i\neq u_j \wedge \bar v_i=\bar v_j)$.
	Clearly $(\phi,\psi)\in\Xim$.
	
	Let $e_{<\omega}$ be a sequence of distinct realisations of the (non-algebraic)
	type $\tp(e/D)$.
	Then the sequence $(e_i\bar d)_{i<\omega}$ witnesses that
	$\tp(\bar a/BC)$ \  $(\phi,\psi)$-divides over~$C$.
	Hence $\bar a\nind[$\Xim$]_CB$.
	
	Conversely, suppose $\bar a\nind[$\Xim$]_CB$.
	So $\tp(\bar a/BC)$ \  $(\phi,\psi)$-divides over~$C$
	for some $(\phi,\psi)\in\Xim$.
	Let this be witnessed by $(e_i\bar d)_{i<\omega}$.
	We may assume that $e_0\bar d\in BC$.
	
	Since $e_i\bar d\equiv_Ce_j\bar d$ for $i<j<\omega$,
	$e_i\equiv_{C\bar d}e_j$ holds as well, so the sequence $e_{<\omega}$
	witnesses that $e_0\not\in\acl(C\bar d)$.
	In particular, $e_0\in BC\setminus C$, so $e_0\in B$.
	Moreover, $\models\phi(\bar a;e_0\bar d)$ implies that
	$e_0\in\acl(\bar a\bar d)\subseteq\acl(\bar aC\bar d)$.
	So choosing $D=C\bar d$ we get
	$e_0\in\acl(\bar aD)\cap B\setminus\acl D$.
	
	(2) For transitivity suppose $D\subseteq C$, $C\ind[$\Xim$]_Dg\bar h$
	and $\phi(\bar y;\bar a)$ \  $(\phi,\psi)$-divides over~$D$,
	where $(\phi,\psi)\in\Xim$.
	Note that $\bar a=g\bar h$ and $g$ is not algebraic over~$D\bar h$.
	By~(1) we have $\acl(C\bar h)\cap\{g\}\subseteq\acl(D\bar h)$,
	so $g$ is not algebraic over $C\bar h$.
	Hence there is a sequence $(\bar g_i)_{i<\omega}$
	of distinct elements $g_i\equiv_{C\bar h}g$.
	It is easy to see that the sequence $(g_i\bar h)_{i<\omega}$
	witnesses the fact that $\phi(\bar y;\bar a)$ \  $(\phi,\psi)$-divides over~$C$.
	
	For normality just observe that $\phi(\bar x,\bar z;\bar y)$
	and $\psi(\bar y_{<k})$ are of the form necessary for $(\phi,\psi)\in\Xim$
	by Definition~\ref{DefinitionXim} iff
	$\phi(\bar x;\bar z,\bar y)$ and
	$\psi'(\bar y_{<k},\bar z_{<k})\equiv\psi(\bar y)\wedge z_0=z_1=\dots=z_{k-1}$
	are of this form.

	(3) Suppose $A\mind_CB$.
	So $\acl(AD)\cap\acl(BD)\subseteq\acl D$ for every set $D$ such that
	$C\subseteq D\subseteq\acl(BC)$.
	Hence $\acl(AD)\cap B\subseteq\acl D$ for every set $D$ such that
	$C\subseteq D\subseteq BC$.
	
	(4) $\thind=\indb[${\scriptscriptstyle\textnormal M}\scriptstyle*$]$\;\, by definition.
	$A\indb[${\scriptscriptstyle\textnormal M}\scriptstyle*$]_CB$ implies
	$A\indb[${\scriptscriptstyle\Xim}\scriptstyle*$]_C\;B$ by (2).
	
	For the converse suppose $A\indb[${\scriptscriptstyle\Xim}\scriptstyle*$]_C\;B$
	holds and we are given $\hat B\supseteq B$.
	Let $A'\equiv_{BC}A$ be such that $A\ind[$\Xim$]_C\acl(\hat BC)$.
	Then $\acl(A'D)\cap\acl(\hat BC)\subseteq\acl D$ for every set $D$
	such that $C\subseteq D\subseteq\acl(\hat BC)$.
	
	Since $\acl(\hat BC)=\acl(\hat BD)$, and since
	$\acl D\subseteq\acl(AD)\cap\acl(\hat BD)$ holds trivially,
	$A'\mind_C\acl(\hat BC)$ follows.
\end{proof}

It is not true in general that $\ind[$\Xim$]\,=\mind$:
Let $T$ be the theory of an everywhere infinite forest,
as in Example~\ref{ExForest2}.
Let $a$, $b_0$ and $b_1$ be nodes such that
$\models Rab_0$, $\models Rab_1$ and $\models b_0\neq b_1$.
Then $a\nmind_{\emptyset}b_0b_1$ because
$a\in\acl(a)\cap\acl(b_0b_1)\setminus\acl\emptyset$.
But $a\ind[$\Xim$]_{\emptyset}b_0b_1$ holds.
This can be verified by checking
$\acl(aD)\cap \{b_0,b_1\}\subseteq\acl D$ for the four possible values
of $D$ such that $\emptyset\subseteq D\subseteq\{b_0,b_1\}$.

\begin{theorem}\label{ThmRosy}\index{theory!rosy}\index{rosy theory}
	$\thind$ is an independence relation for $T$ if and only if the following, equivalent, conditions
	are satisfied:
	\begin{enumerate}
		\item
			$\thind$ satisfies the local character axiom.
		\item
			$\mind$ satisfies the local character axiom.
		\item
			$\D_{\phi,\psi}(\emptyset)<\infty$ for every $(\phi,\psi)\in\Xim$.
		\item
			$A\thind_CB$ implies $B\thind_CA$.
		\item
			$A\thind_CB$ implies $B\mind_CA$.
		\item
			$T$ admits a \sir{}.
	\end{enumerate}
	Moreover, in a theory $T$ satisfying these conditions, $\thind$ is the coarsest \sir.
	
	In particular, $T$ is rosy iff \,$T\eq$ satisfies the equivalent conditions above.
\end{theorem}

\begin{proof}
	First note that $\thind=\indb[${\scriptscriptstyle\Xim}\scriptstyle*$]$\;\;\;
	by Proposition~\ref{PropXim}.
	Therefore we can apply Theorem~\ref{TheoremOmindStar}:
	(1), (3) and (4) are equivalent, and they hold if and only if $\thind$ is an
	independence relation.
	Moreover, they are equivalent to
	(2') $\ind[$\Xim$]$\;\,  satisfies the local character axiom,
	and to
	(5') $A\thind_CB$ implies $B\ind[$\Xim$]_CA$.
	
	The `moreover' statement and the equivalence of (6) with the other
	conditions are by Theorem~\ref{ThmThind}.
	
	Finally, (1) $\implies$ (2) $\implies$ (2') and
	(4) $\implies$ (5) $\implies$ (5')
	since $A\thind_CB\implies A\mind_CB\implies A\ind[$\Xim$]_CB$
	by Proposition~\ref{PropXim} (3) and (4).
\end{proof}

Moreover, by Proposition~\ref{PropXim} (2), if $\thind$ is an independence relation
then we also have an alternative characterisation of $\thind$ by
Theorem~\ref{TheoremRankSymmetry}.

\begin{corollary}\label{CorollaryReductRosy}
	Every reduct of a rosy theory is rosy.
\end{corollary}

\begin{proof}
	Use condition (3) as in Corollary~\ref{CorSimpleReduct}.
\end{proof}

\begin{theorem}\label{ThmMsymmetricSir}\index{theory!M-symmetric}\index{M-symmetric theory}
	The relation $\mind$ is a (strict) independence relation iff
	it is symmetric.
\end{theorem}

\begin{proof}
	If $\mind$ is symmetric,
	condition (5) of Theorem~\ref{ThmRosy} is satisfied, and
	so $\thind$ is an independence relation and therefore satisfies existence.
	Hence $\mind$ satisfies existence as well.
	Using symmetry and transitivity of $\mind$ it is easy to see that
	$\mind$ also satisfies extension, so $\mind=\thind$.
	Conversely, if $\mind$ is an independence relation, then $\mind$ is symmetric by
	Theorem~\ref{ThmSymmetry}.
\end{proof}

The property of $\mind$ being symmetric is not stable under taking reducts:

\begin{example}\label{ExReductNotMSymmetric}\index{example!no preservation of M-symmetry}
	(Symmetry of $\mind$ is not preserved under reducts)\\
	Let $T_0$ be the theory, in the signature of one unary function~$f$,
	which states the following: there is at least one element;
	for every element $b$ there are infinitely many elements $a$ such that $f(a)=b$;
	and $f$ has no periodic points. Note that $T$ is complete.
	
	Let $T$ be the theory extending $T_0$, in the signature consisting of $f$ and a binary
	relation $E$, stating that $\forall xy(Exy\leftrightarrow f(x)=y\vee f(y)=x)$.
	Then $A\mind_CB\iff \acl(AC)\cap\acl(BC)=\acl C$, so $\mind$ is clearly symmetric.
	Yet the theory of an everywhere infinite forest from Example~\ref{ExForest2},
	for which $\mind$ is not symmetric, is a reduct of~$T$.
\end{example}

\begin{exercises}
	\begin{exercise}(M-symmetry)\label{ExcMSymmetric}\\
		(i) Two algebraically closed sets $A$ and $B$ form a \emph{modular pair}\index{modular pair}\index{021@$\operatorname M(x,y)$} in the
		lattice of algebraically closed sets, written $\operatorname M(A,B)$,
		if the following rule holds:
		For any algebraically closed set $C\subseteq B$,
		$\acl(AC)\cap B = \acl(C(A\cap B))$.
		Show that $\operatorname M(A,B)\iff A\mind_{A\cap B}B$.\\
		(ii) A lattice is called \emph{M-symmetric}\index{M-symmetric lattice} if $\operatorname M(A,B)\implies \operatorname M(B,A)$.
		Show that $\mind$ is symmetric iff the lattice of algebraically closed sets
		is M-symmetric.
	\end{exercise}
\end{exercises}

\begin{notes}
	The core of the results presented in this section is, of course,
	from Alf Onshuus.
	The entire development of the theory as presented here, however, is new.
	This is true, in particular, for the use of inconsistency pairs for thorn-forking,
	the definition of $\Xim$, and Theorem~\ref{ThmMsymmetricSir}.
	
	Exercise~\ref{ExcMSymmetric} is from~\cite{Hans Scheuermann: Unabhaengigkeitsrelationen}.
	M-symmetric lattices are studied in a general context
	in~\cite{Manfred Stern: Semimodular Lattices. Theory and Applications}. Note that a finite lattice is
	M-symmetric iff it is semimodular.
	Theorem~\ref{ThmMsymmetricSir} is a generalisation
	of~\cite[Theorems 6.2.8 and 6.2.10]{Hans Scheuermann: Unabhaengigkeitsrelationen},
	which state that $\mind$ is an independence relation
	if the lattice of algebraically closed subsets of the big model is
	M-symmetric and one of two additional conditions is satisfied.
	One of the additional conditions is simplicity of $T$. The other
	is strong atomicity of the lattice of algebraically closed sets
	(a condition that, in conjunction with M-symmetry, is roughly
	equivalent to the Steinitz exchange property for $\acl$).

	Example~\ref{ExReductNotMSymmetric} is due to Wilfrid Hodges and
	first appeared in~\cite{Evans + Pillay + Poizat: Le groupe dans le groupe}.
	It has become the standard example of a 1-based stable
	theory with a reduct that is not 1-based.
	It can also be found in~\cite[Chapter 4, Example 6.1]{Anand Pillay: Geometric Stability Theory}.
\end{notes}

\chapter{Thorn-forking}
\label{Chapter3}

This chapter covers some topics around the concepts of
canonical bases in $T\eq$, both of types and of sequences of
indiscernibles. It is \emph{not} intended to be self-contained.
It is based on Chapter~\ref{Chapter1} and a previous familiarity with simple theories,
but it is completely independent of Chapter~\ref{Chapter2}.
The title is sort of justified by the fact that, as we will see,
only thorn-forking independence can have canonical bases.

If $I$ is a linearly ordered set, we occasionally write
$\bar a_{\in I}$ for the $I$-sequence $(\bar a_i)_{i\in I}$.
Sequences of indiscernibles can be indexed by an arbitrary \emph{infinite}
linearly ordered set~$I$.


\section{Canonical independence relations}

Independence relations (on $T\eq$) that appear in the real world often have the following
property:

\begin{definition}
	A relation $\ind$ has the intersection property
	if it satisfies the following condition:
	\begin{description}\index{axiom!intersection}\index{intersection property}
		\item[(intersection)]\ \\
			Suppose $C_1\subseteq B$ and $C_2\subseteq B$ are such that
			$A\ind_{C_1}B$ and $A\ind_{C_2}B$.
			Then $A\ind_{\acl C_1\cap\acl C_2}B$.
	\end{description}
	
	\noindent An independence relation is \emph{canonical}\index{independence relation!canonical}
	if it is strict and has the intersection property.\index{canonical independence relation}
\end{definition}

Note that if $T$ is a simple theory with elimination of hyperimaginaries
(such as a stable or supersimple theory), then
$\find$ is a canonical independence relation for $T\eq$, because
$\bar a\find_CB$ holds for $C\subseteq B$ iff $\cb(\stp(\bar a/B))\subseteq\acl\eq C$.
It is surprisingly easy to see that there can be at most one \sir{}
with this property:

\begin{lemma}\label{LemAclMorley}
	Suppose $\ind$ is a canonical independence relation.\\
	Then $\bar a\ind_CB$ if and only if there is a sequence
	$(\bar a_i)_{i<\omega}$ of $BC$-indiscernibles that realise $\tp(\bar a/BC)$
	and such that $\acl(C\bar a_{<k})\cap\acl(C\bar a_{\geq k})=\acl C$
	for all $k<\omega$.
\end{lemma}

\begin{proof}
	First suppose that $\bar a\ind_CB$.
	By Proposition~\ref{PropExMS} there is
	a $\ind$-Morley sequence in $\tp(\bar a/BC)$ over~$C$.
	By Proposition~\ref{PropShifting} and finite character,
	$\bar a_{<k}\ind_C\bar a_{\geq k}$ for all $k<\omega$.
	Hence $C\bar a_{<k}\ind_CC\bar a_{\geq k}$.
	By anti-reflexivity, $\acl(C\bar a_{<k})\cap\acl(C\bar a_{\geq k})=\acl C$.
	Note that for this direction it is only necessary that $\ind$ is a strict
	independence relation.
	
	Conversely, suppose there is a sequence
	$(\bar a_i)_{i<\omega}$ of $BC$-indiscernibles realising $\tp(\bar a/BC)$
	and such that $\acl(C\bar a_{<k})\cap\acl(C\bar a_{\geq k})=\acl C$
	for all $k<\omega$.
	Let $\kappa$ be a regular cardinal number sufficiently big for the
	local character axiom.
	Define a totally ordered set $I=\kappa+\{\ast\}+\kappa'$, where $\kappa'$
	is a disjoint copy of $\kappa$ having the opposite order, and $\ast$ is
	a new element greater than any element of $\kappa$ and smaller than any
	element of $\kappa'$.
	Let $(\bar a_i)_{i\in I}$ be a $BC$-indiscernible extension of $(\bar a_i)_{i<\omega}$.
	Note that $\acl(C\bar a_{\in\kappa})\cap\acl(C\bar a_{\in\kappa'})=\acl C$.
	
	By finite character and symmetry it is sufficient to
	prove that $\bar a\ind_C\bar b$ for every finite tuple $\bar b\in B$.
	So let $\bar b\in B$ be a finite tuple.
	
	By local character there is a subset $D\subseteq \bar a_{<\kappa}C$
	such that $|D|<\kappa$ and $\bar b\ind_D\bar a_{<\kappa}C$.
	By regularity of $\kappa$ there is $\lambda<\kappa$ such that
	$D\subseteq\bar a_{<\lambda}C$. Therefore, by base monotonicity,
	$\bar b\ind_{\bar a_{<\lambda}C}\bar a_{<\kappa}C$.
	Now using finite character it is easy to see that
	$\bar b\ind_{\bar a_{<\lambda}C}\bar a_{\in I}C$.
	Hence, using base monotonicity again (and monotonicity),
	$\bar b\ind_{\bar a_{\in\kappa}C}\bar a_{\ast}$.
	
	Since the setup is symmetric with respect to reversing the order of~$I$,
	$\bar b\ind_{\bar a_{\in\kappa'}C}\bar a_{\ast}$ holds as well.
	
	Applying the intersection property we get $\bar b\ind_{\acl C}\bar a_{\ast}$.
	On the other hand $\bar b\ind_C\acl C$ by extension.
	Applying symmetry to the last two statements, and then transitivity,
	we get $\bar a_{\ast}\ind_C\bar b$, hence $\bar a\ind_C\bar b$.
\end{proof}

\begin{theorem}\label{ThmCanonicalThind}
	If $\ind$ is a canonical independence relation, then $\ind=\thind$.
\end{theorem}

\begin{proof}
	Suppose $\ind$ is a canonical independence relation.
	By Theorem~\ref{ThmThind} it follows that $\thind$ is a (strict) independence relation,
	and it is the coarsest, so $\bar a\ind_CB\implies\bar a\thind_CB$.
	Therefore the only thing left to show is $\bar a\thind_CB\implies\bar a\ind_CB$.
	So suppose $\bar a\thind_CB$. As in the first part of the proof of Lemma~\ref{LemAclMorley},
	there is a sequence
	$(\bar a_i)_{i<\omega}$ of $BC$-indiscernibles realising $\tp(\bar a/BC)$
	and such that $\acl(C\bar a_{<k})\cap\acl(C\bar a_{\geq k})=\acl C$
	for all $k<\omega$.
	Hence by Lemma~\ref{LemAclMorley}, $\bar a\ind_CB$.
\end{proof}

There are rosy theories for which $\thind$ is not canonical;
cf.~Section~\ref{SectionQuestions} in the appendix.

\begin{corollary}\label{CorSimpleForkingEHI}
	Suppose $T$ is simple and has elimination of hyperimaginaries
	(e.\,g., $T$ is stable or supersimple).
	Then $\find=\thind$, and this is the only \sir{} for $T\eq$.
\end{corollary}

\begin{proof}
	If $T$ is simple, $\find$ is a \sir{} on $T\eq$, and types
	have canonical bases as hyperimaginary elements.
	If $T$ also has elimination of hyperimaginaries, then
	types have canonical bases as sequences of imaginary elements.
	Therefore $\find$ satisfies the intersection property.
	Hence by the theorem, $\find=\thind$.
	Since $\find$ is the finest \sir{} for $T\eq$ and
	$\thind$ is the coarsest, $\find=\thind$ is the only one.
\end{proof}

\begin{exercises}
	\begin{exercise}(independence in a reduct)\label{ExcCanonicalReduct}
	
		(i) Suppose $T'$ is a reduct of $T$, $\ind$ is a \sir{} for $T$,
		and $\pind$ is a canonical independence relation for~$T'$.
		Let $A,B,C$ be subsets of the big model of $T$ such that $C=\acl C$
		and $A\ind_CB$ in~$T$.
		Then $A\pind_CB$.
		
		(ii) Suppose $T$ is simple and $T'$ is a reduct of $T$ that has
		elimination of hyperimaginaries.
		Let $A,B,C$ be subsets of the big model of $T$ such that $C=\acl\eq C$
		and $A\find_CB$ in~$T$.
		Then $A\find_CB$ holds in $T'$ as well.
	\end{exercise}
\end{exercises}

\begin{notes}
	Lemma~\ref{LemAclMorley} is from~\cite[Theorem~1.7.3]{Hans Scheuermann: Unabhaengigkeitsrelationen}.
	Theorem~\ref{ThmCanonicalThind} is the result of
	combining~\cite[Theorem~1.7.4]{Hans Scheuermann: Unabhaengigkeitsrelationen}
	(which states that any canonical independence
	relation is the coarsest \sir{}) with Theorem~\ref{ThmThind}.
	The definition of the term `canonical independence relation' in~\cite{Hans Scheuermann: Unabhaengigkeitsrelationen}
	was via existence of weak canonical bases.
	Both definitions are in fact equivalent,
	as will be shown in Theorem~\ref{ThmCanonicalSirWCB}.
	
	$\find=\thind$, the first statement of Corollary~\ref{CorSimpleForkingEHI},
	is proved for stable theories in~\cite{Alf Onshuus: Properties and consequences of thorn-independence}, and stated
	for simple theories with stable forking in~\cite{Alf Onshuus: Th-forking algebraic independence and examples of rosy theories}.
	The generalisation to simple theories with elimination of hyperimaginaries was
	first proved, in a different way, by Clifton Ealy.
\end{notes}


\section{Kernels of indiscernibles}

Recall the convention that sequences of indiscernibles are infinite by definition.

\begin{definition}\index{kernel}\index{049@ker}\index{049@aker}
	The \emph{kernel} of a sequence of indiscernibles $(\bar a_i)_{i\in I}$
	is the set
	\[
		\ker(\bar a_i)_{i\in I} = \big\{
			d \;\big|\;
				d\in\dcl(\bar a_{\in I}) \textsl{ and } (\bar a_i)_{i\in I}
				\textsl{ is indiscernible over }d
		\big\}.
	\]
	The \emph{algebraic kernel} of a sequence of indiscernibles $(\bar a_i)_{i\in I}$
	is the set
	\[
		\aker(\bar a_i)_{i\in I} = \big\{
			d \;\big|\;
				d\in\acl(\bar a_{\in I}) \textsl{ and } (\bar a_i)_{i\in I}
				\textsl{ is indiscernible over }d
		\big\}.
	\]
	For the (algebraic) kernel computed in $T\eq$ we write $\ker\eq$ and $\aker\eq$.
\end{definition}

We will see below that $(\bar a_i)_{i\in I}$ is actually indiscernible over
$\ker(\bar a_i)_{i\in I}$, so the kernel is the greatest set definable over the
sequence over which it is indiscernible. Similarly, the algebraic kernel is the
greatest set algebraic over the sequence over which it is indiscernible.

\begin{definition}
	Two sequences $(\bar a_i)_{i\in I}$, $(\bar a_j)_{j\in J}$
	of indiscernibles are
	\emph{cleanly collinear}\index{cleanly collinear}\index{collinear!cleanly} if one of their concatenations
	$(\bar a_i)_{i\in I+J}$ and $(\bar a_i)_{i\in J+I}$ is indiscernible.
	($I+J$ denotes the disjoint union ordered in such a way that $I<J$.)
	We write $\approx$ for the transitive closure of clean collinearity.\index{050@$\approx$}
\end{definition}

\begin{theorem}\label{ThmKernel}
	Properties of the kernel:
	\begin{enumerate}
		\item The kernel is a $\approx$-invariant:\\
			$(\bar a_i)_{i\in I} \approx (\bar b_j)_{j\in J}
				\implies \ker(\bar a_i)_{i\in I} = \ker(\bar b_j)_{j\in J}$.
		\item If $(\bar a_i)_{i\in I}$ and $(\bar b_j)_{j\in J}$
			are cleanly collinear,\\
			then $\ker(\bar a_i)_{i\in I} =
				\dcl(\bar a_{\in I}) \cap \dcl(\bar b_{\in J})$.
		\item $\ker(\bar a_i)_{i\in I} = \dcl\ker(\bar a_i)_{i\in I}$.
		\item $(\bar a_i)_{i\in I}$ is indiscernible over
			$\ker(\bar a_i)_{i\in I}$.
	\end{enumerate}
\end{theorem}

\begin{proof}
	(1) We may assume that $(\bar a_i)_{i\in I}$ and $(\bar b_j)_{j\in J}$
	are cleanly collinear.
	Moreover, it suffices to prove
	$\ker(\bar a_i)_{i\in I} \subseteq \ker (\bar b_j)_{j\in J}$.
	So suppose $d\in\ker(\bar a_i)_{i\in I}$.
	Since $(\bar b_j)_{j\in J}$ is indiscernible over $\bar a_{\in I}$
	and $d\in\dcl\big(\bar a_{\in I}\big)$, $(\bar b_j)_{j\in J}$ is indiscernible over $d$.
	
	Since $d\in\dcl(\bar a_{\in I})$, there are $i_0<\ldots<i_{k-1}$ in $I$ and
	a formula $\phi(\bar x_{<k},\bar y)$ such that
	$\phi(\bar a_{i_0},\ldots,\bar a_{i_{k-1}},\bar y)$
	is a definition of~$d$.
	Choose any $j_0<\ldots<j_{k-1}$ in $J$.
	Then by indiscernibility of the concatenated sequence,
	$\phi(\bar b_{j_0},\ldots,\bar b_{j_{k-1}},\bar y)$
	is also a definition of~$d$.
	Therefore $d\in\dcl(\bar b_{\in J})$.
	
	(2) $\ker(\bar a_i)_{i\in I} \subseteq
	\dcl(\bar b_{\in J})$ follows from (1)
	and the definition of the kernel.
	The opposite inclusion holds because $d\in\dcl(b_{\in J})$
	implies that $(\bar a_i)_{i\in I}$ is indiscernible over $d$.
	
	(3) It follows from (2) that the kernel can be written as the intersection
	of two $\dcl$-closed sets, so it is itself closed under $\dcl$.
	
	(4) Let $(\bar b_j)_{j\in J}$ be cleanly collinear with $(\bar a_i)_{i\in I}$.
	Then $(\bar a_i)_{i\in I}$ is indiscernible over $(\bar b_j)_{j\in J}$,
	hence also over $\dcl\big((\bar b_j)_{j\in J}\big)$,
	hence also over
	$\ker(\bar a_i)_{i\in I}\subseteq\dcl\big((\bar b_j)_{j\in J}\big)$.
\end{proof}

\begin{remark}
	If a sequence $(\bar a_i)_{i\in I}$ is indiscernible over a set $B$,
	then it is also indiscernible over $\acl B$.
\end{remark}

\begin{proof}
	Suppose not. We may assume that $I$ is such that we can extract a
	sequence from $(\bar a_i)_{i\in I}$ that is indiscernible over~$\acl B$.
	Since the fact that the extracted sequence is indiscernible over~$\acl B$
	can be expressed by the type of the sequence over~$B$, the original
	sequence must also be indiscernible over~$\acl B$.
\end{proof}

\begin{lemma}\label{LemKerAker}
	In $T\eq$, the algebraic kernel is just the algebraic closure of the kernel:
		\[
			\aker\eq(\bar a_i)_{i\in I}=
			\acl\eq\ker\eq(\bar a_i)_{i\in I}
		\]
\end{lemma}

\begin{proof}
	Suppose $d\in\acl\eq\ker\eq(\bar a_i)_{i\in I}$.
	Let $\bar d\in\ker\eq(\bar a_i)_{i\in I}$ be a finite tuple such that
	$d\in\acl\eq(\bar d)$.
	Then clearly $d\in\acl\eq(\bar d)\subseteq\acl\eq(\bar a_{\in I})$,
	and $(\bar a_i)_{i\in I}$ is indiscernible over $\bar d$ by
	Theorem~\ref{ThmKernel} (4). Hence $(\bar a_i)_{i\in I}$ is indiscernible
	over $d\in\acl\eq(\bar d)$.
	
	Conversely, suppose $d\in\acl\eq(\bar a_{\in I})$
	and $(\bar a_i)_{i\in I}$ is indiscernible over~$d$.
	Let $D$ be the set of conjugates of $d$ over~$\acl\eq(\bar a_{\in I})$.
	Since $D$ is finite, there is an element $e$ coding it.
	(This is where we need $T\eq$.)
	Clearly $e\in\dcl\eq(\bar a_{\in I})$.
	It is easy to see that $(\bar a_i)_{i\in I}$ is indiscernible over~$d'$
	for every $d'\in D$. Hence $D\subseteq\ker\eq\big((\bar a_i)_{i\in I}\big)$.
	By Theorem~\ref{ThmKernel} (4),
	$(\bar a_i)_{i\in I}$ is indiscernible over~$D$,
	hence also over~$e$.
	Thus $e\in\ker\eq(\bar a_i)_{i\in I})$
	and $d\in\acl\eq(e)\subseteq\acl\eq\ker\eq(\bar a_i)_{i\in I}$.
\end{proof}

\begin{theorem}\label{ThmAKernel}
	Properties of the algebraic kernel:
	\begin{enumerate}
		\item The algebraic kernel is a $\approx$-invariant:\\
			$(\bar a_i)_{i\in I} \approx (\bar b_j)_{j\in J}
				\implies \aker(\bar a_i)_{i\in I} = \aker(\bar b_j)_{j\in J}$.
		\item If $(\bar a_i)_{i\in I}$ and $(\bar b_j)_{j\in J}$
			are cleanly collinear,\\
			then $\aker(\bar a_i)_{i\in I} =
				\acl(\bar a_{\in I}) \cap \acl(\bar b_{\in J})$.
		\item $\aker\bar a_i)_{i\in I} = \acl\aker(\bar a_i)_{i\in I}$.
		\item $(\bar a_i)_{i\in I}$ is indiscernible over
			$\aker(\bar a_i)_{i\in I}$.
	\end{enumerate}
\end{theorem}

\begin{proof}
	(1) It would be straightforward to prove this directly, in the same way
	as Theorem~\ref{ThmAKernel} (1).
	Instead we note the following equation:
	\[
			\aker(\bar a_i)_{i\in I}
			=\monster\cap\aker\eq(\bar a_i)_{i\in I}
			=\monster \cap \acl\eq\ker\eq(\bar a_i)_{i\in I}.
			\tag{$*$}
	\]
	By means this equation the statement is immediate from Theorem~\ref{ThmKernel} (1).

	(2) $\aker(\bar a_i)_{i\in I} \subseteq
	\acl(\bar a_{\in I})$ follows from the definition of the kernel.
	Analogously, $\aker(\bar a_i)_{i\in I}
	=\aker(\bar b_j)_{j\in J} \subseteq
	\acl(\bar b_{\in J})$.
	The opposite inclusion holds because $d\in\acl(b_{\in J})$
	implies that $(\bar a_i)_{i\in I}$ is indiscernible over~$d$.
	
	(3) also follows immediately from ($*$).
	(It is also straightforward to infer (3) from (2) as in Theorem~\ref{ThmKernel} (3).)
	
	(4) $(\bar a_i)_{i\in I}$ is indiscernible over
	$\ker\eq(\bar a_i)_{i\in I}$ by Theorem~\ref{ThmKernel}.
	Hence $(\bar a_i)_{i\in I}$ is indiscernible over
	$\aker(\bar a_i)_{i\in I}$ by ($*$).
\end{proof}

\begin{definition}\label{DefCbIndiscernibles}
	A set of indiscernibles $(\bar a_i)_{i\in I}$ is \emph{based} on
	a set $B$ if every automorphism $\sigma$ of the big model that fixes
	$B$ pointwise satisfies $\big(\sigma(\bar a_i)\big)_{i\in I}\approx(\bar a_i)_{i\in I}$.
	The set $B$ is a \emph{canonical base}\index{canonical base of a sequence of indiscernibles}
	for $(\bar a_i)_{i\in I}$ if the converse is also true.
\end{definition}

If a sequence of indiscernibles has a canonical base at all,
then it is unique up to interdefinability.
But the canonical base of a sequence of indiscernibles need not exist:

\begin{example} (A sequence without a canonical base)\label{ExSequenceNoCb}\index{example!no canonical base for sequence of indiscernibles}\\
	Let $T$ be the theory of a dense linear order without endpoints.
	Consider the type $p\big((x_i)_{i\in\Q}\big) = \{ x_i<x_j \mid i < j\}$.
	Now let $\big((a_i^n)_{i\in\Q}\big)_{n<\omega}$ be a sequence
	of sequences $(a_i^n)_{i\in\Q}$ realising this type,
	and more precisely $a_i^n < a_j^m$ iff either $i<j$, or $i=j$ and $n<m$.
	The sequence is indiscernible, and any two elements are equivalent under
	the type-definable equivalence relation
	$E\big((x_i)_{i\in\Q},(y_i)_{i\in\Q}\big) = \{x_i<y_j \wedge y_i<x_j\mid i<j\}$.
	Therefore every automorphism of the big model that fixes
	the $\approx$-class of the sequence also fixes the $E$-class of the sequence.
	Yet $\ker\big((a_i^n)_{i\in\Q}\big)_{n<\omega}=\emptyset$.
	(Note that $E$ also witnesses that $T$ does not eliminate hyperimaginaries.)
\end{example}

\begin{remark}\label{RemKernelCanonicalParameter}
	Let $(\bar a_i)_{i\in I}$ be an indiscernible sequence.
	Every automorphism~$\sigma$ of the big model that satisfies
	$\big(\sigma(\bar a_i)\big)_{i\in I}\approx(\bar a_i)_{i\in I}$
	fixes $\ker\big((\bar a_i)_{i\in I}\big)$ pointwise.
\end{remark}

\begin{proof}
	Suppose $\sigma$ satisfies $\big(\sigma(\bar a_i)\big)_{i\in I}\approx(\bar a_i)_{i\in I}$.
	It is easy to see that there are sequences $(\bar a_i)_{i\in I}=(\bar a_i^{(0)})_{i\in I},
	(\bar a_i^{(1)})_{i\in I},\dots,(\bar a_i^{(n)})_{i\in I}=\big(\sigma(\bar a_i)\big)_{i\in I}$
	such that $(\bar a_i^{(j)})_{i\in I}$ and $(\bar a_i^{(j+1)})_{i\in I}$ are cleanly collinear.
	Let $\sigma^{(0)}$ be the identity automorphism of the big model, $\sigma^{(n)}=\sigma$
	and $\sigma^{(j)}$ for $0<j<n$ an arbitrary automorphism such that
	$\big(\sigma^{(j)}(\bar a_i)\big)_{i\in I}=(\bar a_i^{(j)})_{i\in I}$.
	It is clearly sufficient to prove the remark for the automorphisms
	$\sigma^{(j+1)}\circ(\sigma^{(j)})^{-1}$.
	
	Therefore we may assume that $(\bar a_i)_{i\in I}$ and
	$\big(\sigma(\bar a_i)\big)_{i\in I}$ are cleanly collinear.
	Since the concatenated sequence is indiscernible over
	$\ker(\bar a_i)_{i\in I}$, there is an automorphism
	$\tau$ of the big model that fixes $\ker(\bar a_i)_{i\in I}$
	pointwise and such that
	$\big(\tau(\bar a_i)\big)_{i\in I}=\big(\sigma(\bar a_i)\big)_{i\in I}$.
	Since $\sigma$ and $\tau$ agree on $\bar a_{\in I}$, they also
	agree on $\ker(\bar a_i)_{i\in I}\subseteq\dcl(\bar a_{\in I})$.
	Thus $\sigma$ also fixes $\ker(\bar a_i)_{i\in I}$ pointwise.
\end{proof}

\begin{corollary}\label{CorCanonicalBaseIndiscernibles}
	If an indiscernible sequence $(\bar a_i)_{i\in I}$ has a canonical base,
	then the canonical base is interdefinable with $\ker(\bar a_i)_{i\in I}$.
\end{corollary}

\begin{proof}
	Let $(\bar a_i)_{i\in I}$ be an indiscernible sequence that has a
	canonical base~$B$.
	By Remark~\ref{RemKernelCanonicalParameter}, every automorphism of the
	big model that fixes $B$ pointwise also fixes $\ker(\bar a_i)_{i\in I}$
	pointwise, so $\ker(\bar a_i)_{i\in I}\subseteq\dcl B$.
	
	For the converse, first note that an automorphism fixing $\bar a_{\in I}$
	also fixes $\bar a_{\in I}/\approx$, hence fixes $B$ pointwise.
	Now consider a sequence $(\bar b_j)_{j\in J}$ that is
	cleanly collinear with $(\bar a_i)_{i\in I}$.
	Then $B\subseteq\dcl(\bar a_{\in I}) \cap \dcl(\bar b_{\in J})$,
	so $B\subseteq\ker(\bar a_i)_{i\in I}$.
	Therefore $\dcl B=\ker(\bar a_i)_{i\in I}$.
\end{proof}

\begin{notes}
	The kernel was first defined by the author of~\cite{Hans Scheuermann: Unabhaengigkeitsrelationen},
	who has changed his surname and his treatment of kernels since then.
	He now calls `algebraic kernel' what he used to call `kernel'.
	The canonical base of a sequence of indiscernibles was defined by
	Steven Buechler in~\cite{Steven Buechler: Canonical bases in some supersimple theories}.
	The relation of this notion to the kernel as stated in
	Corollary~\ref{CorCanonicalBaseIndiscernibles} is new.
	
	Generalisation of this section (and, in fact, the entire chapter)
	to hyperimaginaries is essentially straightforward.
	A paper carrying this out is in preparation.
\end{notes}


\section{Weak canonical bases}

In the following variant of Lemma~\ref{LemAclMorley} the hard direction is
strengthened:

\begin{lemma}\label{LemAclMorleyImproved}
	Suppose $\ind$ is a canonical independence relation.\\
	Then $\bar a\ind_CB$ if and only if there is a sequence
	$(\bar a_i)_{i<\omega}$ of $BC$-indiscernibles that realise $\tp(\bar a/BC)$
	and such that $\acl(\bar a_{<k})\cap\acl(\bar a_{\geq k})\subseteq\acl C$
	for all $k<\omega$.
\end{lemma}

\begin{proof}
	The `only if' part is a weakening of the `only if' part of
	Lemma~\ref{LemAclMorleyImproved}.
	
	Conversely, suppose there is a sequence
	$(\bar a_i)_{i<\omega}$ of $BC$-indiscernibles realising $\tp(\bar a/BC)$
	and such that $\acl(\bar a_{<k})\cap\acl(\bar a_{\geq k})\subseteq\acl C$
	for all $k<\omega$.
	Let $(\kappa=|T|+|BC|)^+$.
	Define a totally ordered set $I=\kappa+\{\ast\}+\kappa'$, where $\kappa'$
	is a disjoint copy of $\kappa$ having the opposite order, and $\ast$ is
	a new element greater than any element of $\kappa$ and smaller than any 
	element of $\kappa'$.
	Let $(\bar a_i)_{i\in I}$ be a $BC$-indiscernible extension of $(\bar a_i)_{i<\omega}$.
	Note that $\acl(\bar a_{\in\kappa})\cap\acl(\bar a_{\in\kappa'})\subseteq\acl C$.
	
	By symmetry it is sufficient to
	prove that $\bar a\ind_CB$.
	
	By extended local character 
	(see Exercise~\ref{ExcLocalCharacter} (ii))
	there is a subset $D\subseteq \bar a_{<\kappa}$
	such that $|D|<\kappa$ and $BC\ind_D\bar a_{<\kappa}$.
	By regularity of $\kappa$ there is $\lambda<\kappa$ such that
	$D\subseteq\bar a_{<\lambda}$. Therefore, by base monotonicity,
	$BC\ind_{\bar a_{<\lambda}}\bar a_{<\kappa}$.
	Now using finite character it is easy to see that
	$BC\ind_{\bar a_{<\lambda}}\bar a_{\in I}$.
	Hence, using base monotonicity again (and monotonicity),
	$BC\ind_{\bar a_{\in\kappa}}\bar a_{\ast}$.
	
	Since the setup is symmetric with respect to reversing the order of~$I$,
	$BC\ind_{\bar a_{\in\kappa'}}\bar a_{\ast}$ holds as well.
	
	Applying the intersection property we get $BC\ind_{D}\bar a_{\ast}$,
	where $D$ is defined as $D=\acl(\bar a_{\in\kappa})\cap\acl(\bar a_{\in\kappa'})\subseteq\acl C$.
	By symmetry, $\bar a_{\ast}\ind_DBC$, and by extension,
	$\bar a_{\ast}\ind_DB\acl C$.
	Hence $\bar a_{\ast}\ind_{\acl C} B$ by base monotonicity and monotonicity.
	
	On the other hand $\bar a_{\ast}\ind_C\acl C$ by extension.
	Applying symmetry to the last two statements, and then transitivity,
	we get $B\ind_C\bar a_{\ast}$, hence $B\ind_C\bar a$, hence $\bar a\ind_CB$.
\end{proof}

The following theorem is essentially just a more conceptual re-formulation
of Lemma~\ref{LemAclMorleyImproved}:

\begin{theorem}\label{ThmCanonicalAker}
	Let $\ind$ be a canonical independence relation.\\
	Then $\bar a\ind_CB$ holds if and only if there is a sequence $(\bar a_i)_{i<\omega}$
	of indiscernibles realising $\tp(\bar a/BC)$ such that
	$\aker(\bar a_i)_{i<\omega}\subseteq\acl C$.
\end{theorem}

\begin{proof}
	To reduce the theorem to Lemma~\ref{LemAclMorleyImproved}, it suffices to prove
	that for a $BC$-indiscernible sequence $(\bar a_i)_{i<\omega}$,
	$\aker(\bar a_i)_{i<\omega}\subseteq\acl C$ holds if and only if
	$\acl(\bar a_{<k})\cap\acl(\bar a_{\geq k})\subseteq\acl C$ holds
	for all $k<\omega$.
	
	For the `if' part of the claim, suppose
	$\aker\big((\bar a_i)_{i<\omega}\big)\subseteq\acl C$ holds.
	Let $I$ be a totally ordered set and $J=I+\omega$.
	Let $(\bar a_j)_{j\in J}$ be a $BC$-indiscernible extension
	of $(\bar a_i)_{i<\omega}$.
	Then for every $k<\omega$, $(\bar a_i)_{i\in J, i<k}$
	and $(\bar a_i)_{i\in J, i\geq k}$ are cleanly collinear sequences
	of $BC$-indiscernibles. Hence
	$\acl(\bar a_{<k})\cap\acl(\bar a_{\geq k})
		\subseteq\aker(\bar a_i)_{i\in J}\subseteq\acl C$.
		
	For the `only if' part of the claim,
	suppose $\acl(\bar a_{<k})\cap\acl(\bar a_{\geq k})\subseteq\acl C$ holds
	for all $k<\omega$ and consider an arbitrary element
	$d\in\aker\bar a_i)_{i<\omega}$ of the algebraic kernel.
	Let $\phi(\bar x_{<k},y)$ be such that $\phi(\bar a_{<k},y)$ is an
	algebraic formula realised by~$d$.
	Since $(\bar a_i)_{i<\omega}$ is indiscernible over~$d$,
	$\phi(\bar a_{k+1},\ldots,\bar a_{2k},y)$ is also an
	algebraic formula realised by~$d$.
	Hence $d\in\acl(\bar a_{<k})\cap\acl(\bar a_{\geq k})\subseteq\acl C$.
\end{proof}

\begin{corollary}\label{CorAkerMorley}
	Let $\ind$ be a canonical independence relation.\\
	A sequence $(\bar a_i)_{i<\omega}$ of $C$-indiscernibles
	is a $\ind$-Morley sequence for $\tp(\bar a_0/C)$ if and only if
	$\aker(\bar a_i)_{i<\omega}\subseteq\acl C$.
\end{corollary}

\begin{proof}
	First suppose $\aker(\bar a_i)_{i<\omega}\subseteq\acl C$.
	Then $\bar a_0\ind_C\{\bar a_i\mid 0<i<\omega\}$
	by Theorem~\ref{ThmCanonicalAker},
	so $(\bar a_i)_{i<\omega}$ is $\ind$-independent over~$C$, hence a
	$\ind$-Morley sequence over~$C$.
	Conversely, suppose $(\bar a_i)_{i<\omega}$ is a $\ind$-Morley
	sequence over~$C$.
	Consider a $C$-indiscernible extension
	$(\bar a_i)_{i<\omega+\omega}$ of $(\bar a_i)_{i<\omega}$.
	Then $\bar a_{<\omega}\ind_C\{\bar a_i\mid\omega\leq i<\omega+\omega\}$,
	so $\aker(\bar a_i)_{i<\omega}\subseteq\acl C$ by
	anti-reflexivity of~$\ind$.
\end{proof}

\begin{definition}
	Let $\ind$ be an independence relation, $\bar a$ a tuple and $B=\acl B$ an
	algebraically closed set.
	An algebraically closed set $C\subseteq B$ is called the
	\emph{weak canonical base}\index{weak canonical base} (with respect to $\ind$) of $\tp(\bar a/B)$
	if $C$ is the smallest algebraically closed subset of~$B$
	such that $\bar a\ind_CB$.
	
	We say that $\ind$ \emph{has weak canonical bases} if $\tp(\bar a/B)$ has
	a weak canonical base for every tuple
	$\bar a$ and every algebraically closed set~$B$.
\end{definition}

Thus $C=\acl C\subseteq B=\acl B$ is the weak canonical base of
$\tp(\bar a/B)$ iff $\bar a\ind_CB$ holds, and for every $D=\acl D\subseteq B$
such that $\bar a\ind_DB$ we have $C\subseteq D$.
If a type has a weak canonical base, then it is of course unique.

\begin{theorem}\label{ThmCanonicalSirWCB}
	A strict independence relation has weak canonical bases if and only if
	it is canonical.
	
	Moreover, if $\ind$ is a canonical independence relation and $B$ is an
	algebraically closed set,
	then the weak canonical base of a type $\tp(\bar a/B)$
	is $\aker(\bar a_i)_{i<\omega}$, where $(\bar a_i)_{i<\omega}$
	is an arbitrary $\ind$-Morley sequence for $\tp(\bar a/B)$.
\end{theorem}

\begin{proof}
	First suppose $\ind$ has weak canonical bases
	and $C_1\subseteq B$ and $C_2\subseteq B$ are such that
	$\bar a\ind_{C_1}B$ and $\bar a\ind_{C_1}B$.
	It easily follows that $\bar a\ind_{\acl C_1}\acl B$
	and $\bar a\ind_{\acl C_1}\acl B$ also hold.
	Thus if $C$ is the weak canonical base of $\tp(\bar a/\acl B)$,
	then $C\subseteq\acl C_1$ and $C\subseteq\acl C_2$.
	Hence $C\subseteq\acl C_1\cap\acl C_2$, and so
	$\bar a\ind_{\acl C_1\cap\acl C_2}\acl B$ by base monotonicity.
	
	Conversely, suppose $\ind$ is a canonical independence relation.
	Let $B$ be an algebraically closed set and
	$(\bar a_i)_{i<\omega}$ a $\ind$-Morley sequence
	for $\tp(\bar a/B)$.
	Then $C=\aker(\bar a_i)_{i<\omega}$ is a weak
	canonical base for $\tp(\bar a/B)$ (note that this also
	proves the `moreover' part):
	
	First we observe that in fact $C\subseteq B$ by Corollary~\ref{CorAkerMorley}.
	It follows from Theorem~\ref{ThmCanonicalAker} that $\bar a\ind_CB$.
	Finally we need to prove minimality of~$C$.
	So suppose $C'\subseteq B$ and $\bar a\ind_{C'}B$.
	Then $(\bar a_i)_{i<\omega}$ is a $\ind$-Morley sequence over~$C'$.
	Hence by Corollary~\ref{CorAkerMorley},
	$C=\aker(\bar a_i)_{i<\omega}\subseteq\acl C'$ holds.
\end{proof}

\begin{exercises}
	\begin{exercise}(weak canonical bases)\label{ExcWcb}
		
		Suppose $\ind$ is a canonical independence relation.
		Given a tuple $\bar a$ and an algebraically closed set $B$,
		let $\wcb(\bar a/B)$\index{055@wcb} denote the weak canonical base of $\tp(\bar a/B)$.
		
		(i) $\bar a\equiv_B\bar a'$ implies $\wcb(\bar a/\acl B)=\wcb(\bar a'/\acl B)$.
		Hence it is in no way misleading to define $\wcb(\bar a/B)=\wcb(\bar a/\acl B)$
		if $B$ is not algebraically closed.
		
		(ii) If $\bar a'$ is a subtuple of $\bar a$, then
		$\wcb(\bar a'/B)\subseteq\wcb(\bar a/B)$.
		
		(iii) The following conditions are equivalent for a tuple $\bar a$ and sets $B,C$:\\
		(1) $\bar a\ind_CB$,\quad
		(2) $\wcb(\bar a/BC)\subseteq\acl C$,\quad
		(3) $\wcb(\bar a/BC)=\wcb(\bar a/C)$.
	\end{exercise}
\end{exercises}

\begin{notes}
	Theorem~\ref{ThmCanonicalSirWCB} is new.
	Weak canonical bases in the sense of this section were defined in~\cite{Hans Scheuermann: Unabhaengigkeitsrelationen}.
	They first appeared, with a different definition, in~\cite{Evans + Pillay + Poizat: Le groupe dans le groupe}.
	The definition in~\cite{Evans + Pillay + Poizat: Le groupe dans le groupe} implies stability,
	and both definitions are equivalent for stable theories.
\end{notes}


\section{Canonical bases in simple theories}

If $T$ is a stable theory or, more generally, a simple theory with elimination of
hyperimaginaries, then forking independence is a canonical independence relation for $T\eq$
because an amalgamation base $\tp(\bar a/B)$ does not fork over a subset $C\subseteq B$
if and only if $\cb(\bar a/B)\subseteq\acl\eq C$, and so $\acl\eq\cb(\bar a/B)$
is the weak canonical base of~$\tp(\bar a/B)$.
Hence $\cb(\bar a/B)$ is characterised by Theorem~\ref{ThmCanonicalSirWCB}
up to interalgebraicity.
(Note that $\cb(\bar a/B)$ is only defined up to interdefinability. So it makes
sense to regard $\cb(\bar a/B)$ as a definably closed set.)
In this section we will try to improve this result.
We will also see that the phenomenon exhibited by Example~\ref{ExSequenceNoCb}
cannot occur in this context.

We will need the well-known fact that the type of an indiscernible sequence over a
cleanly collinear sequence is an amalgamation base:

\begin{remark}\label{RemarkCollStp}
	Suppose $\bar a_{\in I}$ and $\bar b_{\in J}$ are cleanly collinear sequences
	of indiscernibles. Then $\tp(\bar a_{\in I}/\bar b_{\in J})$ is a strong type.
\end{remark}

\begin{proof}
	We may assume that the concatenation $\bar a_{\in I}{}\concat\bar b_{\in J}$ is indiscernible
	(otherwise exchange $\bar a_{\in I}$ and $\bar b_{\in J}$).
	Suppose $\bar a'_{\in I}$ realises $\tp(\bar a_{\in I}/\bar b_{\in J})$.
	Consider a finite equivalence relation $\phi(\bar x_0\dots\bar x_{n-1},\bar x'_0\dots\bar x'_{n-1})$
	definable over~$\bar b_{\in J}$.
	We need to show that
	$\models\phi(\bar a_{i_0}\dots\bar a_{i_{n-1}},\bar a'_{i_0}\dots\bar a'_{i_{n-1}})$
	for any $i_0<\dots<i_{n-1}$ in $I$.
	
	Note that by indiscernibility all tuples $\bar a_{i_0}\dots\bar a_{n-1}$
	($i_0<\dots<i_{n-1}$ in $I$) must be equivalent, since otherwise all would
	have to be inequivalent, which is impossible because $\phi$ has only finitely
	many classes. By compactness there is a sequence $\bar c_{<\omega}$
	such that $\bar a_{\in I}{}\concat \bar c_{<\omega}{}\concat \bar b_{\in J}$
	and $\bar a_{\in I}{}\concat \bar c_{<\omega}{}\concat \bar b_{\in J}$
	are both indiscernible.
	Now clearly $\bar a_{i_0}\dots\bar a_{i_{n-1}}$,
	$\bar c_0\dots\bar c_{n-1}$ and $\bar a'_{i_0}\dots\bar a'_{i_{n-1}}$
	are $\phi$-equivalent.
\end{proof}

For the rest of this section we work in $T\eq$ for a fixed simple theory $T$
with elimination of hyperimaginaries (EHI),\index{EHI, elimination of hyperimaginaries}
and we freely use well-known facts about simple theories
(cf.~\cite{Byunghan Kim: Forking in simple unstable theories}, \cite{Kim + Pillay: Simple theories} and~\cite{Hart + Kim + Pillay: Coordinatisation and canonical bases in simple theories}).

\begin{lemma}($T$ simple with EHI)\label{LemmaCBCPSimple}\\
	Let $\bar a_{\in I}$ and $\bar b_{\in J}$ be cleanly collinear sequences
	of indiscernibles and $C=\cb(\bar a_{\in I}/\bar b_{\in J})$.
	If $\bar a'_{\in I}\equiv_C\bar a_{\in I}$,
	then $\bar a'_{\in I}\approx\bar a_{\in I}$.
\end{lemma}

\begin{proof}
	Let $\bar b'_{\in J}$ be such that
	$\bar a_{\in I}\bar b_{\in J}\equiv_C\bar a'_{\in I}\bar b'_{\in J}$.
	Without loss of generality $\bar b_{\in J}\ind_C\bar b'_{\in J}$
	(otherwise we can find
	$\bar a''_{\in I}\bar b''_{\in J}\equiv_C\bar a_{\in I}\bar b_{\in J}$
	independent from $\bar b_{\in J}\bar b'_{\in J}$ over~$C$).
	Note that $\tp(\bar a_{\in I}/\bar b_{\in J})$
	and $\tp(\bar a'_{\in I}/\bar b'_{\in J})$ are non-forking extensions
	of the same amalgamation base $\tp(\bar a_{\in I}/\bar b_{\in J})$.
	By the amalgamation theorem for amalgamation bases there is a type
	$\tp(\bar c_{\in I}/\bar b_{\in J}\bar b'_{\in J})$ which is 
	a common non-forking extension of both.
	Since $\bar c_{\in I}$ is cleanly collinear with both $\bar b_{\in J}$
	and $\bar b'_{\in J}$, we have
	$\bar a_{\in I}\approx\bar b_{\in J}\approx\bar c_{\in I}\approx\bar b'_{\in J}\approx\bar a'_{\in I}$.
\end{proof}

\begin{remark}($T$ simple with EHI)\label{RemarkCBKerSimple}\\
	If $\bar b_{\in J}$ is indiscernible over~$\bar a$,
	then the relation $\cb(\bar a/\bar b_{\in J})\subseteq\ker\bar b_{\in J}$ holds.
\end{remark}

\begin{proof}
	Let $\bar c_{<\omega}$ be cleanly collinear with $\bar b_{\in J}$
	over~$\bar a$.
	Then $\cb(\bar a/\bar b_{\in J}) = \cb(\bar a/\bar c_{<\omega})$
	is definable over $\bar b_{\in J}$ and over $\bar c_{\in K}$,
	hence
	$\cb(\bar a/\bar b_{\in J})\subseteq\dcl\bar b_{\in J}\cap\dcl\bar c_{<\omega}=\ker\bar b_{\in J}$.
\end{proof}

Putting together the last two results and Remark~\ref{RemKernelCanonicalParameter},
we get something quite close to the desired refinement of Theorem~\ref{ThmCanonicalSirWCB}:

\begin{theorem}($T$ simple with EHI)\\
	If $\bar a_{\in I}$ and $\bar b_{\in J}$ are cleanly collinear sequences
	of indiscernibles, then
	\[
		\cb(\bar a_{\in I} / \bar b_{\in J}) = \ker\bar a_{\in I},
	\]
	and this is a canonical base for $\bar a_{\in I}$ in the sense of
	Definition~\ref{DefCbIndiscernibles}.
\end{theorem}

\begin{proof}
	An automorphism $\sigma$ of the big model that fixes pointwise $\ker a_{\in I}$
	fixes $\cb(\bar a_{\in I}/\bar b_{\in J})$ pointwise by Remark~\ref{RemarkCBKerSimple}.
	An automorphism $\sigma$ that fixes pointwise $\cb(\bar a_{\in I}/\bar b_{\in J})$
	satisfies $\sigma(\bar a_{\in I})\approx\bar a_{\in I}$ by
	Lemma~\ref{LemmaCBCPSimple}.
	An automorphism $\sigma$ that satisfies $\sigma(\bar a_{\in I})\approx\bar a_{\in I}$
	fixes $\ker a_{\in I}$ pointwise by Remark~\ref{RemKernelCanonicalParameter}.
\end{proof}

\begin{corollary}($T$ simple with EHI)\label{CorollarySimpleCb}\\
	Let $p(\bar x)$ be an amalgamation base.
	If $\bar a_{<\omega}$ is a Morley sequence for $p(\bar x)$, then
	\[
		\cb(p) \subseteq \ker(\bar a_{<\omega}) \subseteq \acl\cb(p).
	\]
\end{corollary}

\begin{proof}
	By compactness there is a sequence $\bar b_{<\omega}$ such that
	the concatenation $\bar b_{<\omega}{}\concat\bar a_{<\omega}$
	is still a Morley sequence for~$p$.
	Note that $\tp(\bar b_{<\omega}/\bar a_{<\omega})$ is an
	amalgamation base by Remark~\ref{RemarkCollStp}.
	Hence $\cb(p) = \cb(\bar b_0/\bar a_{<\omega}) \subseteq
	\cb(\bar b_{<\omega}/\bar a_{<\omega}) = \ker(\bar a_{\in I})$.
	
	For the inclusion on the right-hand side just observe that
	$\ker(\bar a_{<\omega}) \subseteq \aker(\bar a_{<\omega}) = \acl\cb(p)$
	because $\aker(\bar a_{<\omega})$ is the weak canonical base
	of~$p$ by Theorem~\ref{ThmCanonicalSirWCB}.
\end{proof}

I would have liked to show that $\cb\big(\stp(\bar a/B)\big)=\ker\bar a_{<\omega}$,
but here is a (perfectly trivial 1-based supersimple) counter-example:

\begin{example} (Alzheimer's random graph)\\
	Consider the following theory $T$:
	There are two sorts $N$ and $C$ and a partial function
	$f:N\times N\rightarrow C$. $C$ has precisely 2 elements
	(`edge' and `no edge'). $f(a,b)$ is defined iff $a\neq b$.
	If we fix an element $e\in C$, the relation $f(x,y)=e$
	defines a random graph on~$N$.
	Note that $\dcl\emptyset=\emptyset$, while $\acl\emptyset=C$.
	Like the theory of the random graph $T$ is supersimple.
	
	We have $\tp(a/\emptyset)\vdash\stp(a/\emptyset)$ for every
	single element $a\in N$. For distinct $a,b\in N$, however,
	we have $\tp(ab/\emptyset)\not\vdash\stp(ab/\emptyset)$
	because $\stp(ab/\emptyset)$ fixes $f(a,b)$ while
	$\tp(ab/\emptyset)$ does not.
	
	Hence for any Morley sequence $(a_i)_{i<\omega+2}$
	in $\tp(a/\emptyset)$ we have $\cb(a_\omega/a_{<\omega}) = \emptyset$
	while $\cb(a_\omega a_{\omega+1}/a_{<\omega}) = C = \ker(a_{<\omega})$.
\end{example}

Yet it turns out that the stronger statement does hold
for stable theories:

\begin{corollary}($T$ stable)\\
	Let $p(\bar x)$ be a stationary type.
	If $\bar a_{<\omega}$ is a Morley sequence for $p(\bar x)$, then
	\[
		\cb(p) = \ker\bar a_{<\omega}.
	\]
\end{corollary}

\begin{proof}
	It is sufficient to show that $\cb(\bar b_0/\bar a_{<\omega}) =
	\cb(\bar b_{<\omega}/\bar a_{<\omega})$ holds in the proof of
	Corollary~\ref{CorollarySimpleCb}.
	Let $C=\cb(\bar b_0/\bar a_{<\omega})$.
	Since clearly $C\subseteq\cb(\bar b_{<\omega}/\bar a_{<\omega})$
	we need only prove $\cb(\bar b_{<\omega}/\bar a_{<\omega})\subseteq C$.
	$\bar b_{<\omega}\ind_C\bar a_{<\omega}$ holds by a standard
	forking calculation, so the only thing left to show is that
	$\tp(\bar b_{<\omega}/C)$ is stationary.
	But this is true because $\tp(\bar a_i/C)$ is stationary for every
	$i<\omega$ and $\bar a_{<\omega}$ is independent over~$C$.
\end{proof}

\begin{exercises}
	\begin{exercise}\label{ExcOneBased} (1-based theories, cf.~Exercise~\ref{ExcModular})\\
		Let $T$ be a simple theory with elimination of hyperimaginaries.\\
		(i) $T$ is called 1-based if $A\find_{\acl\eq A\cap\acl\eq B}B$ holds for
		all $A$, $B$.\index{theory!1-based}
		Show that $T$ is 1-based iff the lattice of algebraically closed sets
		in the big model of $T\eq$ is modular.\index{lattice!modular}\index{modular lattice}\\
		(ii) Show that $T$ is 1-based and perfectly trivial iff the lattice of algebraically
		closed sets in the big model of $T\eq$ is distributive.\index{lattice!distributive}\index{distributive lattice}
	\end{exercise}
\end{exercises}

\begin{notes}
	Exercise~\ref{ExcOneBased} is from~\cite{Hans Scheuermann: Unabhaengigkeitsrelationen},
	the rest is new.
	The entire chapter can be generalised to hyperimaginaries---a paper
	containing the details is in preparation.
\end{notes}


\appendix
\chapter{Appendix}

\section{Thorn-forking---the official definition}\label{SectionThornDefinitions}

Thorn-forking (or \th-forking) is a notion of independence first defined by Alf Onshuus
as the notion of independence corresponding to certain ranks (the thorn-ranks
or \th-ranks) suggested by Thomas Scanlon.
In this section I present the original definition and prove that it is equivalent
to the definition in Section~\ref{SectionThornForking} if it is read in $T\eq$ (as is the original definition).

Therefore we work in $T\eq$ throughout.
Here is the definition from~\cite[Definition 2.1]{Alf Onshuus: Properties and consequences of thorn-independence}:

\begin{definition}\label{DefOnshuus}
	\index{dividing!strong}
	\index{strong dividing}
	\index{059@\th|see{thorn}}
	\index{thorn-dividing}
	\index{thorn-forking}
	Let $\phi(x,y)$ be an $\mathcal L$-formula without parameters, let $b$ be an element,
	and let $C$ be a set.
	\begin{itemize}
		\item The formula $\phi(x,b)$ is said to \emph{strongly divide} over~$C$ if\\
			$\tp(b/C)$ is not algebraic
			and $\big\{\phi(x,b') \;\big|\; {b'\models\tp(b/C)}\big\}$
			is $k$-inconsistent for some natural number $k < \omega$.
		\item The formula $\phi(x,b)$ is said to \emph{\th-divide} over~$C$ if\\
			there is a tuple $c$ such that
			$\phi(x,b)$ strongly divides over $Cc$.
		\item The formula $\phi(x,b)$ is said to \emph{\th-fork} over~$C$ if\\
			it implies a (finite) disjunction
			of formulas (with arbitrary parameters), each of which \th-divides over~$C$.
	\end{itemize}
\end{definition}

\noindent Let us begin with some easy observations about these definitions:

Following Onshuus, we never mentioned identifying logically equivalent formulas (so we will not do it)
or required that $b$ be a canonical parameter of $\phi(x,b)$. Thus, if $\phi(x,b)$ strongly divides
over~$C$, $z$ is a variable in a non-algebraic sort, $\psi(x,y,z)$ is the
formula $\phi(x,y)\wedge z=z$, and $c$ is an element of the same sort as $z$, then even though
$\phi(x,b)$ and $\psi(x,bc)$ are equivalent, $\psi(x,bc)$ does not strongly divide over $C$.
More generally, a formula containing an unused parameter of a non-algebraic sort can never
\th-divide over any
set.

Using compactness, one can prove that a formula $\phi(x,b)$ strongly divides over~$C$ if and only if
the set $\{\phi(x,b')\mid b'\models\tp(b/C)\}$ of its $C$-conjugates is infinite and has no infinite
consistent subset. (This is~\cite[Remark 2.1.1]{Alf Onshuus: Properties and consequences of thorn-independence}.)

Again using compactness, it is easy to see that in the definition of \th-dividing it does not matter
whether we demand that $c$ be a \emph{finite} tuple. In fact, a formula $\phi(x,b)$ \th-divides over
a set~$C$ if and only if it strongly divides over a superset $C'\supseteq C$ of~$C$.

\begin{proposition}\label{propaltdef}
    Our new definition of \th-forking ($\thind$ as defined in Definition~\ref{DefMThorn})
    agrees with the original definition by Alf Onshuus:
    
	$A\thind_CB$ if and only if for every tuple $a\in A$ and
	every tuple $b\in BC$ and $\mathcal L$-formula $\phi(x,y)$
	such that $\models\phi(a,b)$, the formula $\phi(x,b)$ does not \th-fork over~$C$.
\end{proposition}

\begin{proof}
	We first prove the implication from left to right, assuming that for some
	$a\in A$ a formula in the type $\tp(a/BC)$ \th-forks over~$C$.
	So we have $\tp(a/BC)\models\phi_1(x,b_1)\vee\ldots\vee\phi_n(x,b_n)$, where
	each $\phi_i(x,b_i)$ strongly divides over some superset $C_i\supseteq C$.
	Now choosing $\hat B=BC_1\ldots C_nb_1\ldots b_n$ we can demonstrate
	that $A'\nmind_C\hat B$ for all $A'\equ_{BC}A$:
	
	For any $A'\equ_{BC}A$ there is $a'\equ_{BC}a$, $a'\in A'$,
	and $i$ such that $\models\phi_i(a',b_i)$ holds. Since $\phi_i(x,b_i)$ strongly divides over~$C_i$,
	we have the following: $b_i\not\in\acl C_i$,
	but only finitely many realisations $b_i'\models\tp(b_i/C_i)$ can simultaneously
	satisfy $\phi_i(a',b_i')$.
	Thus $b_i\in(\acl(C_iA') \cap\acl\hat B)\setminus\acl C_i$, so $A'\nmind_C\hat B$.
	
	We will now prove the implication from right to left in two steps.
	
	First we claim the following:
	Suppose $B=\acl(BC)$.
	If there is no $a\in A$ s.t.\ a formula in $\tp(a/BC)$ \th-forks over~$C$,
	then $A\mind_CB$.
	
	Suppose $A\nmind_CB$.
	Then there is a set $C'$ satisfying $C\subseteq C'\subseteq B$ such that
	$\acl(AC')\cap B\supsetneq\acl C'$.
	Let $b\in\acl(AC')\cap B\setminus\acl C'$ witness this.
	Let $a\in A$, $c\in C'$ and $\phi(x,y,z)$ be such that the formula $\phi(a,y,c)$ is algebraic
	(has at most $k$ realisations, say) and realised by~$b$.
	Define $\phi'(x,y,z)$ as follows: $\phi(x,y,z) \wedge \exists_{\leq k}y\phi(x,y,z)$.
	Then $\models\phi'(a,b,c)$ also holds.
	The formula $\phi'(x,b,c)$ strongly divides over $C'$ since $bc\not\in\acl C'$ (because $b\not\in\acl C'$)
	and every consistent subset of
	$\big\{\phi'(x,b',c') \;\big|\; b'c'\models\tp(bc/C')\big\}$ has at most $k$ elements.
	Since $\phi'(x,b,c)\in\tp(a/B)$ it follows that $A\nthind_CB$.

	Now we can finish the proof. Suppose there is no $a\in A$ such that
	a formula in $\tp(a/BC)$ \th-forks over~$C$, and consider any $\hat B\supseteq B$.
	By~\cite[Lemma~2.1.2(1)]{Alf Onshuus: Properties and consequences of thorn-independence}
	(the proof of which works for infinite sets as well as for finite tuples)
	there is $A'\equ_{BC}A$ such that when $\bar a'$ enumerates $A'$,
	$\tp(\bar a'/\acl(\hat BC))$ does not \th-fork over~$C$.
	By the claim this implies that $A'\mind_C\acl(\hat BC)$, so we get $A'\mind_C\hat B$.
\end{proof}

\begin{exercises}
	\begin{exercise}\label{ExcOnshuus}(strong dividing and \th-forking)
		
		Let us say that a formula $\phi(x,b,d)$ \emph{divides quite strongly}
		over~$C$ if the set
		$\big\{\phi(x,b',d) \;\big|\; {b'\models\tp(b/C)}\big\}$
		of formulas \emph{up to equivalence} is $k$-inconsistent for some $k<\omega$.
		Clearly if a formula strongly divides then it also divides quite strongly,
		though the converse does not hold in general.
		Show that the converse does hold after passing to \th-forking:
		
		If there is a tuple $c$ such that $\phi(x,b,d)$ divides quite strongly
		over~$Cc$, then $\phi(x,bd)$ is equivalent to a formula that strongly divides
		over~$Ccd$. (Hence $\phi(x,bd)$ \th-forks over~$C$.)
	\end{exercise}
\end{exercises}

\begin{notes}
	The content of this section is original inasmuch as it proves the equivalence
	of a new definition of thorn-forking to Alf Onshuus' definition.

	It should be noted that there are some ambiguities in~\cite[Definition 2.1]{Alf Onshuus: Properties and consequences of thorn-independence}
	that I passed over in silence.
	Exercise~\ref{ExcOnshuus} shows that these do not affect the definition of \th-forking.
\end{notes}


\section{Solutions for exercises}
\newenvironment{details}[1]{\noindent\textbf{Additional details for Example \ref{#1}.}}{}
\newenvironment{solution}[1]{\noindent\textbf{Solution to Exercise \ref{#1}.}}{}
\label{SectionExercises}

{\small

\begin{solution}{ExcExistence} (relations between axioms, existence and symmetry) 
	
	(i) $A\ind_CB$ implies $B\ind_CA$ by symmetry.
	Applying extension to $\hat A=AC\supseteq A$ we get $B'\equiv_{AC}B$
	such that $B'\ind_C\hat A$, i.\,e., $B'\ind_CAC$. By invariance also $B\ind_C\hat A$.
	Hence $AC\ind_CB$ by symmetry.
	
	(ii) Let $A_0\subseteq A$ be any finite subset and
	note that $A_0\ind_CC$ by local character and base monotonicity.
	So by finite character, $A\ind_CC$ holds.
	Now for any $B$ we can use extension to find $A'\equ_CA$ such that $A'\ind_CBC$,
	so $A'\ind_CB$ by monotonicity.
	
	(iii) Suppose $A\ind_CB$ and $B\subseteq\hat B$.
	By existence there is $A'\equiv_{BC}A$ such that $A'\ind_{BC}\hat B$.
	By invariance, $A'\ind_CB$.
	Using Remark~\ref{RemarkTransitivity} on the other side, which is possible because
	of symmetry, we get $A'\ind_C\hat B$.
\end{solution}

\begin{solution}{ExcLocalCharacter} (local character)

	(i) It easily follows from existence that $A\ind_AB$.
	Also, $A\ind_BB$ by existence and invariance.
	Therefore we can choose $C_1=B$ and $C_2=A$.
	
	(ii) By invariance the statement is clear for finite sets~$A$.
	Given arbitrary sets $A$ and $B$, we can find for every
	finite subset $A_0\subseteq A$ a subset $C_0\subseteq B$
	such that $A_0\ind_{C_0}B$ and $\card{C_0}<\kappa$.
	Let $C$ be the union of all these sets $C_0$.
	Then $A\ind_CB$ by finite character and base monotonicity,
	and clearly $\card C<\kappa+\card T^+$.
\end{solution}

\begin{solution}{ExcModular} (modularity and distributivity)

	(i)
	It is sufficient to show:
	If $A$, $B$ and $C=\acl C$ are s.\,t.{} $A\cap C=B\cap C=\emptyset$,
	then there is $A'\equiv_CA$ s.\,t.{} $A'\cap B=\emptyset$.
	If this were false, then by compactness there would be a counter-example
	with $A$ and $B$ finite. Towards a contradiction let $A$, $B$ and $C=\acl C$ be s.\,t.{}
	$A\cap C=B\cap C=\emptyset$, $A'\cap B\not=\emptyset$ for all $A'\equiv_CA$,
	$B$ finite and $\card A$ minimal for these properties.
	Let $A_*\subseteq A$ be a maximal subset of $A$ s.\,t.{}
	$\{A'\mid A'\equiv_CA\text{ and }A_*\subseteq A'\}$ is infinite.
	By minimality of $\card A$ there is $A_*'\equiv_CA_*$ s.\,t.{}
	$A_*'\cap B=\emptyset$. We may assume $A_*'=A_*$, so
	$A_*\cap B=\emptyset$.
	For every $b\in B$, by maximality of $A_*$ there are only finitely many
	$A'\equiv_CA$ s.\,t.{} $A_*\cup b\subseteq A'$. Hence there are only
	finitely many $A'\equiv_CA$ s.\,t.{} $A_*\subseteq A'$ and $A'\cap B\not=\emptyset$.
	Since $\{A'\mid A'\equiv_CA\text{ and }A_*\subseteq A'\}$
	is infinite, there is $A'\equiv_CA$ s.\,t.{} $A_*\subseteq A$ and
	$A'\cap B=\emptyset$, a contradiction.
	
	(ii)
	Invariance, monotonicity and normality are obvious.
	
	\emph{Finite character:} Suppose $d\in\acl(AC)\cap\acl(BC)\setminus\acl C$.
	Then $d$ is already algebraic over a finite tuple $\bar a\bar c$ with
	$\bar a\in A$ and $\bar c\in C$, and also over a finite tuple $\bar b\bar c'$
	with $\bar b\in B$ and $\bar c'\in C$.
	Hence $d\in\acl(C\bar a)\cap\acl(C\bar b)\setminus\acl C$.
	
	\emph{Transitivity:}
	Suppose $D\subseteq C\subseteq B$.
	If $\acl B\cap\acl(AC)\subseteq\acl C$
	and $\acl C\cap\acl(AD)\subseteq\acl D$,
	then $\acl B\cap\acl(AD)\subseteq\acl C\cap\acl(AD)\subseteq\acl D$.
	
	\emph{Extension:} Using Exercise~\ref{ExcExistence} (iii) this follows from (i),
	symmetry and the other axioms already shown to hold.
	
	\emph{Local character:} Given sets $A$ and $B$,
	construct sets $C_i\subseteq B$ and $D_i$ ($i<\omega$)
	as follows: $C_0=D_0=\emptyset$.
	$D_{i+1}=\acl(AC_i)\cap\acl B$.
	For every $d\in D_{i+1}$ let $\bar c_d\in B$ be a finite tuple such that
	$d\in\acl\bar c_d$. Let $C_{i+1}$ be the union over all tuples $\bar c_d$.
	Let $C=\bigcup_{i<\omega}C_i$.
	It is easy to see that $C\subseteq B$ and $\card C\leq\card T+\card A$.
	Moreover, if $d\in\acl(AC)\cap\acl(BC)$,
	then already $d\in\acl(AC_i)\cap\acl(BC)\subseteq D_{i+1}$
	for some $i<\omega$, hence $d\in\acl C_{i+1}\subseteq\acl C$.
	
	(iii)
	Let $A$ and $B\supseteq C$ be algebraically closed sets.
	First note that $B\cap\acl(AC)\supseteq\acl((B\cap A)C)$
	holds without any further assumptions.
	Now suppose $\aind$ satisfies the base monotonicity axiom.
	Then $A\aind_{A\cap B}B$ implies $A\aind_{(A\cap B)C}B$.
	Hence $B\cap\acl(AC)\subseteq\acl((B\cap A)C)$.
	
	Conversely, suppose the modular law holds,
	$A\aind_CB$ and $C\subseteq C'\subseteq B$.
	Then $\acl B\cap\acl(AC')\subseteq\acl((\acl B\cap\acl A)C')$
	by modularity and $C'\subseteq B$.
	Note that $\acl B\cap\acl A\subseteq C\subseteq C'$,
	so $\acl B\cap\acl(AC')\subseteq\acl C'$.
	Hence $A\aind_{C'}B$.
	
	(iv)~Let
	$A$, $B$, $C$ be algebraically closed sets.
	First note that $\acl((A\cap B)C)\subseteq\acl(AC)\cap\acl(BC)$
	holds without any further assumptions.
	Now suppose $\aind$ is perfectly trivial.
	Since $A\ind_{A\cap B}B$ it follows that $A\ind_{(A\cap B)C}B$
	as well, hence $AC\ind_{(A\cap B)C}BC$ (by base monotonicity,
	which can be applied on both sides due to symmetry),
	hence $\acl(AC)\cap\acl(BC)\subseteq\acl((A\cap B)C)$.
	
	Conversely, suppose the lattice is distributive,
	$A\ind_CB$ and $C'\supseteq C$.
	Then $\acl(AC')\cap\acl(BC')=\acl((\acl A\cap\acl B)C')\subseteq\acl(CC')=\acl C'$,
	hence $A\ind_{C'}B$.
\end{solution}

\begin{solution}{ExcForest1} (concerning Example~\ref{ExForest1})

	We write $[a,b]$ for the set of nodes on the path from $a$ to $b$
	if this path exists. Let $\dist(a,b)$ be the size of $[a,b]$ minus one,
	or $\infty$.
	
	$A\ind_CB\implies AC\ind_CB$ is trivial.
	
	For $A\ind_CB\implies \acl A\ind_CB$ suppose $A\ind_C B$,
	$a\in\acl A$ and $b\in B$.
	Let $a_1,a_2\in A$ be s.\,t.{} $a\in[a_1,a_2]$.
	Suppose $[a,b]$ exists.
	Let $d\in[a_1,a_2]$ be s.\,t.{} $\dist(d,b)$ is minimal.
	Then $[a_1,b]=[a_1,d]\cup[d,b]$ and $[a_2,b]=[a_2,d]\cup[d,b]$.
	Since $[a_1,b]$ and $[a_2,b]$ meet $\acl C$, it follows that
	$[d,b]$ meets $\acl C$ (otherwise $d\in\acl C$, a contradiction).
	Hence $[a,b]=[a,d]\cup[d,b]$ also meets $\acl C$.
	
	Invariance, monotonicity, finite character, base monotonicity,
	normality and anti-reflexivity are trivial.
	
	\emph{Transitivity:} Suppose $B\ind_CA$, $C\ind_DA$ and $D\subseteq C\subseteq B$.
	We need to show that $B\ind_DA$, i.\,e.,
	every path from $B$ to $A$ meets $\acl D$.
	Let $b\in B$, $a\in A$ be s.\,t.{} $[b,a]$ exists.
	Let $d\in[b,a]\cap\acl C$. Then $[d,a]\subseteq[b,a]$,
	and $[d,a]$ meets $\acl D$ by $\acl C\ind_DA$.
	Hence $[b,a]$ meets $\acl D$.
	
	\emph{Extension:}
	We prove existence instead:
	Given $A,B,C$ write $A$ as a disjoint union
	$A=(A\cap \acl C)\cup \bigcup A_i$, where each $A_i$ is the
	intersection of $A$ with a connected component of
	$\monster\setminus\acl C$.
	Let $A_i'$ be s.\,t.{} $A_i'\equiv_{\acl C}A_i$, and the connected
	component of $A_i'$ in $\acl C$ avoids $B$ and the other $A_i'$.
	Then $A'=(A\cap\acl C)\cup\bigcup A_i'\equiv_{\acl C}A$,
	and $A'\ind_CB$.
	Since $\ind$ is obviously symmetric,
	extension now follows as in Exercise~\ref{ExcExistence}.
	
	\emph{Local character:}
	For $a\in A$ s.\,t.{} a path from $a$ to $B$ exists let
	$c_a\in\acl B$ be s.\,t.{} $\dist(a,c_a)$ is minimal,
	and let $b_a$, $b_{a'}\in B$ be s.\,t.{} $c_a\in[b_a,b_{a'}]$.
	Then $C=\{c_a\mid a\in A\}\subseteq B$, $\card C\leq 2\card{A}$,
	and every path from $a\in A$ to $B$ meets $\acl C$ in~$c_a$.
\end{solution}

\begin{solution}{ExcDividingForking} (dividing and forking of formulas)

	(i)
	First suppose $\bar b\in BC$, $\models\phi(\bar a,\bar b)$,
	and $\phi(\bar x,\bar b)$ divides over~$C$.
	Let this be witnessed by $(\bar b_i)_{i<\omega}$
	such that $\bar b_i\equiv_C\bar b$ and
	$\{\phi(\bar x,\bar b_i)\mid i<\omega\}$ is $k$-inconsistent.
	It is not hard to see that we may assume that $(\bar b_i)_{i<\omega}$ is
	$C$-indiscernible, and that $\bar b_0=\bar b\in BC$.
	Thus $(\bar b_i)_{i<\omega}$ witnesses that $\bar a\ndind_CB$.
	
	For the converse, suppose $\bar a\ndind_CB$.
	This is witnessed by a sequence $(\bar b_i)_{i<\omega}$ of $C$-indiscernibles
	with $\bar b_0\in BC$ and such that there is no
	$\bar a'\equiv_{BC}\bar a$ such that
	$(\bar b_i)_{i<\omega}$ is $\bar a'C$-indiscernible.
	We may assume that $\bar b_0$ enumerates all elements of $BC$.
	(This involves extracting a sequence of indiscernibles.)
	Let $p(\bar x;\bar y)=\tp(\bar a\bar b_0/C)$.
	Then $\bigcup_{i<\omega}p(\bar x;\bar b_i)$ is inconsistent.
	(Otherwise the set would still be consistent after extending
	the sequence $(\bar b_i)_{i<\omega}$. We could realise it
	by $\bar a^*$, say, and extract an $\bar a^*C$-indiscernible
	sequence $(\bar b^*_i)_{i<\omega}$.
	The $C$-automorphism taking $(\bar b^*_i)_{i<\omega}$ to $(\bar b_i)_{i<\omega}$
	would take $\bar a^*$ to $\bar a'\equiv_{BC}\bar a$
	because $\bar a'\bar b_0\equiv_C\bar a^*\bar b^*_0\equiv_C\bar a\bar b_0$.)
	
	(ii)
	We prove only the harder direction.
	Suppose $\bar a\nfind_CB$,
	so there is $\hat B\supseteq B$ such that
	$\bar a'\ndind_C\hat B$ for all $\bar a'\equiv_{BC}\bar a$.
	By (i), $\tp(\bar a/BC)\cup\{\neg\phi(\bar x,\bar b)\mid\bar b\in\hat BC;\;\;
	\phi(\bar x,\bar b) \textsl{ divides over } C\}$ is inconsistent.
	So there are a formula $\psi(\bar x,\bar b)\in\tp(\bar a/BC)$
	and formulas $\phi_i(\bar x,\bar b_i)$, $i<k$ dividing over~$C$
	with parameters $\bar b_i\in\hat BC$
	such that $\psi(\bar x,\bar b)$ implies the disjunction
	$\bigvee_{i<k}\phi_i(\bar x,\bar b_i)$.
	
	(iii)
    Suppose $\phi(\bar x;\bar b)$ does not divide over $C$.
    Let $(\bar b_i)_{i<\kappa}$ be a Morley sequence for~$\tp(\bar b/C)$.
    Then $\{\phi(\bar x;\bar b_i)\mid i<\kappa\}$ is consistent,
    hence realised by a tuple $\bar a$, say.
    If $\kappa$ is sufficiently big, we can extract from
    $(\bar b_i)_{i<\kappa}$ a $C\bar a$-indiscernible sequence
    $(\bar b'_i)_{i<\omega}$. By Proposition~\ref{PropExMS},
    $\bar a\find_C\bar b'_0$, hence $\bar a\find_C\bar b$.
    Since $\models\phi(\bar a;\bar b)$ and $\bar a\find_C\bar b$,
    $\phi(\bar x;\bar b)$ does not fork over~$C$.
\end{solution}

\begin{solution}{ExcPropertiesDind}(additional properties of $\dind$)

	(i) Let $(\bar a_i)_{i<\omega}$ be a sequence of $B$-indiscernibles.
	For a sufficiently big cardinal $\kappa$ let $(\bar a_i)_{i<\omega}$
	be a $B$-indiscernible extension. (This exists by compactness.)
	Let $\bar b_0$ be an enumeration of $\acl B$.
	For every $i\in\kappa\setminus\{0\}$ let $\bar b_i$ be such
	that $\bar a_0\bar b_0\equiv_B\bar a_i\bar b_i$.
	If $\kappa$ was chosen big enough we can extract from the sequence
	$(\bar a_i\bar b_i)_{i<\kappa}$ a $B$-indiscernible sequence
	$(\bar a_i'\bar b_i')_{i<\omega}$.
	Now it is easy to see that we could have extended the tuples
	$\bar a_i$ in such a way by tuples $\bar b_i$ that
	$(\bar a_i\bar b_i)_{i<\omega}$ is $B$-indiscernible.
	
	We now show $\bar b_0=\bar b_1$. Otherwise there is an index $j$
	such that $\bar b_0^j\neq \bar b_1^j$. But then $\bar b_i^j\neq\bar b_{i'}^j$
	for all $i\neq i'$, a contradiction since all $\bar b_i^j$ satisfy
	the same algebraic type over~$B$. Therefore $\bar b_0=\bar b_1$,
	hence $\bar b_i=\bar b_j$ for all $i,j<\omega$.
	Now it easily follows that $(\bar a_i)_{i<\omega}$ is $\acl B$-indiscernible.
	
	(ii) Suppose $b_0\in B\setminus\acl C$.
	Then for every $\kappa$ there is a $C$-indiscernible sequence $(b_i)_{i<\omega}$ of
	distinct realisations of $\tp(b_0/B)$. (Start with a very long sequence of
	distinct realisations and extract a sequence of indiscernibles from it.)
	By $A\dind_CB$ there is $A'\equiv_{BC}A$ such that
	$(b_i)_{i<\omega}$ is $A'C$-indiscernible.
	By (i) the sequence is also $\acl(A'C)$-indiscernible.
	Hence $b_0\not\in\acl(A'C)$.
	
	(iii) Suppose $\bar b_0\in\acl(BC)$ and $(\bar b_i)_{i<\omega}$
	is a sequence of $C$-indiscernibles. First we note that we may
	assume that $\bar b_0$ actually enumerates all of $\acl(BC)$:
	For a sufficiently big cardinal	$\kappa$ we extend the sequence
	to a sequence $(\bar b_i)_{i<\kappa}$ of $C$-indiscernibles.
	Let $\bar b_0'\supseteq\bar b_0$ be a tuple enumerating $\acl(BC)$.
	Then we can find for every $i\in\kappa\setminus\{0\}$ a tuple
	$\bar b_i'\supseteq\bar b_i$ s.\,t.{} $\bar b_i'\equiv_C\bar b_i$.
	Now extract a sequence of $C$-indiscernibles from $(\bar b'_i)_{i<\kappa}$.
	
	Let $\bar b_0^*\subseteq\bar b_0$ be the subtuple of $\bar b_0$ that
	enumerates $BC$. Let $\bar b_i^*$ be the tuple corresponding
	to $\bar b_0^*$ in $\bar b_i$.
	Then there is $A'\equiv_{BC}A$ such that $(\bar b_i^*)$ is $A'C$-indiscernible.
	By compactness it is possible to extend the sequence to a sequence
	$(\bar b_i)_{i<\kappa}$ that is $C$-indiscernible and such that
	$(\bar b_i^*)_{i<\kappa}$ is $A'C$-indiscernible if $\bar b_i^*$
	is defined as before for $i\geq\omega$.
	If $\kappa$ is big enough we can extract a sequence of $A'C$-indiscernibles.
	Thus we get a sequence $\bar b'_i$ of $A'C$-indiscernibles
	such that $\bar b'_0\equiv_{A'}\bar b_0$ and $\bar b^{\prime*}_0\equiv_{A'C}\bar b^*_0$.
	Since $\bar b'_0$ enumerates $\acl(BC)$ and $\bar b_0^{\prime*}$ enumerates $BC$
	it follows that there is a permutation of $\bar b_0'$ that fixes $\bar b_0^{\prime*}$
	and such that $\bar b'_0\equiv_{A'C}\bar b_0$.
	Now we can take an $A'C$-automorphism of the big model that maps $\bar b_0'$
	to $\bar b_0$ and fixes the permutation issue.
	It takes $(\bar b'_i)_{i<\omega}$ to an $A'C$-indiscernible sequence
	$(\bar b''_i)_{i<\omega}$ such that $\bar b''_0=\bar b_0$.
	
	(iv) Suppose $A\dind_CB$. First note that $A\dind_C\acl(BC)$ by (iii).
	Now suppose $C\subseteq C'\subseteq\acl(BC)$.
	Since $\dind$ satisfies base monotonicity by Lemma~\ref{LemmaBasicDividing},
	we have $A\dind_{C'}\acl(BC)$.
	Hence $\acl(AC')\cap\acl(BC)=\acl C$ by (ii).
\end{solution}

\begin{solution}{ExcForest2} (concerning Example~\ref{ExForest2})

	We show that $A\nind_CB\implies A\nthind_CB$, so
	suppose $A\nind_CB$. Then there are nodes $a\in A$ and $b\in B$
	such that $[a,b]$ (exists and) avoids $\acl C$.
	First consider the case that there is a node $c\in C$ in the same connected component as
	$a$ and~$b$.
	Choose $c'\in[a,b]$ s.\,t.{} $\dist(c',c)$ is minimal.
	Then $c'\in\acl(aC)\cap\acl(bC)\setminus\acl C$,
	so $a\nmind_Cb$ and hence $A\nthind_CB$.
	In the second case $C$ does not meet the connected component of $a$ and~$b$.
	Let $d,e$ be two distinct neighbours of $b$.
	If, towards a contradiction, $A\thind_CB$, then also $a\thind_Cb$, hence there
	must be $d'e'\equiv_{bC}de$ such that $a\mind_Cbd'e'$.
	Since $d',e'$ are distinct neighbours of $b$, one of them, $d'$ say,
	is not in $[a,b]$. But then $b\in[a,d']$,
	so $b\in\acl(aCd')\cap\acl(bCd')\setminus\acl(Cd')$,
	hence $a\nmind_Cbd'$, a contradiction. Therefore $A\nthind_CB$ in both cases.
\end{solution}

\begin{solution}{ExcNotInteresting} (concerning Example~\ref{ExNotInteresting})

	Since $\acl A=A$ for all $A$, the lattice of algebraically closed sets is distributive.
	Hence $\aind$ is a (perfectly trivial) \sir{} by Exercise~\ref{ExcModular}.
	$\aind$ is clearly coarser than $\thind$. On the other hand $\thind$ is coarser
	than $\aind$ by Theorem~\ref{ThmThind}.
	Therefore $\aind=\thind$. Clearly $A\aind_CB$ iff $AC\cap BC=C$,
	iff $A\cap B\subseteq C$.
\end{solution}

\begin{solution}{ExcTwoSirs} (concerning Example~\ref{ExTwoSirs})

	For thorn-forking for $T$ we are in the situation of Example~\ref{ExNotInteresting}.
	In $T\eq$ we have $\acl\eq(AB)=\acl\eq A\cup\acl\eq B$, so the lattice
	of algebraically closed sets for $T\eq$ is also distributive and we can
	argue exactly as in Example~\ref{ExNotInteresting}.
\end{solution}

\begin{solution}{ExcNoExtensionNoExistence} (concerning Example~\ref{ExNoExtensionNoExistence})

	Transitivity: Suppose $B\pind_{CD}A$ and $C\pind_DA$.
	Consider the three cases for $C\pind_DA$. If $A\subseteq D$,
	then $BC\pind_DA$ follows trivially.
	If $C\subseteq D$, then $CD=D$, so $BC\pind_DA$ also follows trivially.
	Otherwise $D$ is infinite, $B\ind_{CD}A$ and $C\ind_DA$,
	hence $BC\ind_DA$, hence $BC\pind_DA$.
	Local character:
	Let $\kappa$ be suitable for local character of $\ind$.
	For $B$ and finite $A$ there is $C\subseteq B$ s.\,t.{}
	$\card C<\kappa$ and $A\ind_CB$. If $B$ is finite choose
	$C'=B$, otherwise choose $C'$ s.\,t.{} $C\subseteq C'\subseteq B$
	and $\aleph_0\leq\card{C'}<\kappa$.
	No extension, no existence: Consider the case of finite~$C$.
\end{solution}


\begin{solution}{ExcNoExtensionNoSymmetry} (strong finite character of $\mind$)

	Suppose $A\nmind_CB$.
	Then there is $D$ such that $C\subseteq D\subseteq\acl(BC)$ and an element
	$e\in\big(\acl(AD)\cap\acl(BD)\big)\setminus\acl D$.
	Let $\bar a$, $\bar b$ and $\bar c$ be enumerations of $A$, $B$ and $C$, respectively.
	
	Since $e\in\acl(AD)$, we can find a finite tuple $\bar d\in D$ and an algebraic formula
	$\alpha(u,\bar a,\bar d)$ such that $\models\alpha(d,\bar a,\bar d)$.
	Then for appropriate $k<\omega$, $e$ satisfies the formula $\alpha'(u,\bar a,\bar d)$
	defined as $\alpha(u,\bar a,\bar d)\wedge\exists_{\leq k}u'\alpha(u',\bar a,\bar d)$.
	
	Since $e\in\acl(BD)=\acl(BC)$, there is an algebraic formula
	$\beta(u,\bar b,\bar c)$ such that $\models\beta(e,\bar b,\bar c)$.
	Let $e_0,\ldots,e_{n-1}$ be all the realisations of $\beta(u,\bar b,\bar c)$ that are in~$\acl D$.
	
	Let $\chi(u,\bar d^*)$ be an algebraic formula with parameters in $D$ that is
	satisfied at least by $e_0,\ldots,e_{n-1}$. We may assume that $\bar d=\bar d^*$.
	Note that every element $e'$ that satisfies $\beta(u,\bar b,\bar c)$,
	either satisfies $\chi(u,\bar d)$ or is not algebraic over $C\bar d$ at all.
	
	Let $\delta(\bar v,\bar b,\bar c)$ be an isolating formula in
	the algebraic type $\tp(\bar d/B\cup C)$.
	Note that for any $\bar d'$ satisfying $\delta(\bar v,\bar b,\bar c)$,
	every element $e'$ that satisfies $\beta(u,\bar b,\bar c)$
	either satisfies $\chi(u,\bar d')$ or is not algebraic over $C\bar d'$ at all.
	
	Let $\phi(\bar x,\bar b,\bar c)$ be the formula defined as
	\[
		\exists u\exists\bar v\big(
			\delta(\bar v,\bar b,\bar c)
			\wedge\alpha'(u,\bar x,\bar v)
			\wedge\beta(u,\bar b,\bar c)
			\wedge\neg\chi(u,\bar v)
		\big).
	\]
	\noindent $\phi(\bar x,\bar b,\bar c)$ has the property desired:
	
	First note that $e$ and $\bar d$ witness that $\models\phi(\bar a,\bar b,\bar c)$ holds.
	On the other hand, suppose $\models\phi(\bar a',\bar b,\bar c)$ holds
	and let $e'$ and $\bar d'$ witness this, i.e.,
	\[
		\models\delta(\bar d',\bar b,\bar c)
		\wedge\alpha'(e',\bar a',\bar d')
		\wedge\beta(e',\bar b,\bar c)
		\wedge\neg\chi(e',\bar d').
	\]
	
	Let $D'=C\bar d'$. From $\delta(\bar d',\bar b,\bar c)$ we get
	$C\subseteq D'\subseteq\acl(BC)$.
	From $\models\alpha'(e',\bar a',\bar d')$ we get $e'\in\acl(\bar a'\bar d')\subseteq\acl(D'\bar a')$.
	From $\models\beta(e',\bar b,\bar c)\wedge\neg\chi(e',\bar d')$ we get
	$e'\in\acl(BC)$ and $e'\not\in\acl(C\bar d')=\acl D'$.
	Hence $e'$ witnesses $\acl(D'\bar a)\cap\acl(BD')\supsetneq\acl D'$.
\end{solution}

\begin{solution}{ExcFigureClassification} (Figure \ref{FigureClassification})

	We need to prove that the relations defined in Examples~\ref{ExNoLCNoExistence},
	\ref{ExNoExtensionNoExistence}, \ref{ExNoExtensionNoSymmetry} and~\ref{ExNoLCNoSymmetry}
	satisfy strong finite character.
	For Example~\ref{ExNoLCNoExistence} this is trivial.
	In Examples~\ref{ExNoExtensionNoSymmetry} and~\ref{ExNoLCNoSymmetry}
	we are dealing with $\mind$, so we can just use Exercise~\ref{ExcNoExtensionNoSymmetry}.
	
	Example~\ref{ExNoExtensionNoExistence}:
	We can choose for $\ind$ a \sir{} satisfying strong finite character.
	For infinite $C$ we have $A\pind_CB\iff A\ind_CB$, so we can just use
	strong finite character of~$\ind$.
	For finite $C$, suppose $A\npind_CB$.
	Choose $a\in A\setminus C$, $b\in B\setminus C$,
	and let $\bar c$ be an enumeration of $C$.
	Let $\phi(\bar a,\bar b,\bar c)$ express that
	neither $a$ nor $b$ is among the elements of the tuple $\bar c$.
\end{solution}

\begin{solution}{ExcStrongFiniteCharacter} (alternative definition for strong finite character)

	For one direction, suppose $\ind$ satisfies extension, and
	for any sequence of variables $\bar x$ and any sets $B$ and $C$,
	the set $\big\{\tp(\bar a/BC) \;\big|\; \bar a\nind_CB\big\}$ is an open
	subset of $\S_{\bar x}(BC)$. We will show that $\ind$ satisfies
	strong finite character.
	So suppose $A\nind_CB$. Let $\bar a$ be an enumeration of~$A$.
	The basic open sets of $\S_{\bar x}(BC)$ are those of the form
	$\big\{\tp(\bar a/BC) \;\big|\; \models\phi(\bar x,\bar b,\bar c)\big\}$,
	where $\phi(\bar x,\bar y,\bar z)$ is a formula without parameters
	and $\bar b\in B$ and $\bar c\in C$ are finite tuples.
	Let $\phi(\bar x,\bar y,\bar z)$ and $\bar b\in B$, $\bar c\in C$
	be such that
	$\tp(\bar a/BC)\in\big\{\tp(\bar a/BC) \;\big|\; \models\phi(\bar a,\bar b,\bar c)\big\}
	\subseteq\big\{\tp(\bar a/BC) \;\big|\; \bar a\nind_CB\big\}$.
	Then clearly $\models\phi(\bar a,\bar b,\bar c)$ and
	$\bar a'\nind_CB$ for all $\bar a'$ satisfying $\models\phi(\bar a',\bar b,\bar c)$.
	By extension, $\bar a'\nind_C\bar b$ for all $\bar a'$ satisfying $\models\phi(\bar a',\bar b,\bar c)$.
	Since only finitely many variables from $\bar x$ really occur
	in $\phi(\bar x,\bar y,\bar z)$, we may replace $\bar x$ by a finite
	subtuple $\bar x_0$ and $\bar a$ by a finite subtuple $\bar a_0$.%

	Conversely, suppose $\ind$ satisfies monotonicity and strong finite character.
	We will show that for any sequence of variables $\bar x$ and any sets $B$ and $C$,
	the set $\big\{\tp(\bar a/BC) \;\big|\; \bar a\nind_CB\big\}$ is an open
	subset of $\S_{\bar x}(BC)$.
	So suppose $\tp(\bar a/BC)\in\big\{\tp(\bar a/BC) \;\big|\; \bar a\nind_CB\big\}$.
	Let $A$ be the set of elements of $\bar a$.
	By monotonicity and strong finite character there are
	tuples $\bar e\in A$, $\bar b\in B$, $\bar c\in C$ and
	a formula $\phi(\bar u,\bar y,\bar z)$ without parameters
	such that $\models\phi(\bar e,\bar b,\bar c)$, and
	$\bar e'\nind_CB$ for all $\bar e'$ satisfying
	$\models\phi(\bar e',\bar b,\bar c)$.
	We can write $\bar e=(a_{i_0},a_{i_1},\ldots,a_{i_{k-1}})$.
	Let $\psi(\bar x)\equiv\phi(x_{i_0},x_{i_1},\ldots,x_{i_{k-1}})$.
	Then $\models\psi(\bar a,\bar b,\bar c)$, and
	$\bar a'\nind_CB$ for all $\bar a'$ such that $\models\psi(\bar a',\bar b,\bar c)$.
	So $\tp(\bar a/BC)\in\big\{\tp(\bar a'/BC) \;\big|\; \models\phi(\bar a',\bar b,\bar c)\big\}
	\subseteq\big\{\tp(\bar a'/BC) \;\big|\; \bar a'\nind_CB\big\}$.
\end{solution}

\begin{solution}{ExcForkStar} ($\Delta$-forking)

	First suppose $\tp(\bar a/BC)$ $\Delta$-forks over~$C$ for some $\Delta\subseteq\Omega$.
	Then $\tp(\bar a/BC)$ implies some disjunction $\bigvee_{i<n}\phi^i(\bar x;\bar b^i)$,
	where each of the formulas $\phi(\bar x;\bar b^i)$ \  $(\phi^i,\psi^i)$-divides over $C$
	and $(\phi^i,\psi^i)\in\Delta\subseteq\Xi$.
	(Note that we admit $\bar b^i\not\in BC$.)
	Let $\hat B=BC\bar b^{<n}$.
	Then every $\bar a'\equiv_{BC}\bar a$ realises $\tp(\bar a/BC)$,
	so $\models\bigvee_{i<n}\phi^i(\bar a';\bar b^i)$, and so
	$\bar a'\nomind_C\bar b^i$ for some $i<n$.
	Hence $\bar a'\nomind_C\hat B$.
	Therefore $\bar a\nind[$\Omega\ast$]_C\hat B$.
	
	Conversely, suppose $\tp(\bar a/BC)$ does not $\Delta$-fork over~$C$
	for any finite $\Delta\subseteq\Omega$.
	Then for any $\hat B\supseteq B$ the partial type
	\[
		\tp(\bar a/BC)\cup\big\{
			\neg\phi(\bar x;\bar b)
		\;\big|\;
			\bar b\in\hat B,\;
			(\phi,\psi)\in\Omega\restrict\bar x,\;
			\text{$\phi(\bar x;\bar b)$ \  $(\phi,\psi)$-divides over $C$}
		\big\}
	\]
	is consistent. Let $\bar a'$ realise this type.
	Then clearly $\bar a'\equiv_{BC}\bar a$ and $\bar a'\omind_C\hat B$.
\end{solution}

\begin{solution}{ExcDividing} (more on dividing and forking of formulas)\ 

	(i) First suppose $\phi(\bar x;\bar b)$ divides over~$C$,
	so there is a number $k<\omega$ and a sequence $(\bar b_i)_{i<\omega}$
	such that $\bar b_i\equiv_C\bar b$ and
	$\{\phi(\bar x;\bar b_i)\mid i<\omega\}$ is $k$-inconsistent.
	Consider the formula
	$\psi(\bar x_{<k})\equiv\neg\exists\bar x\bigwedge_{i<k}\phi(\bar x;\bar y_i)$,
	which is clearly a $k$-inconsistency witness for $\phi(\bar x;\bar y)$.
	Now each $\bar b_i$ realises $\tp(\bar b/C)$,
	and $\psi(\bar b_{i_0},\ldots,\bar b_{i_{k-1}})$
	holds for all $i_0<\ldots<i_{k-1}<\omega$ by $k$-inconsistency.
	Therefore $\phi(\bar x;\bar b)$ $(\phi,\psi)$-divides over~$C$.
	Conversely, suppose that $\phi(\bar x;\bar b)$ $(\phi,\psi)$-divides over~$C$,
	where $\psi(\bar y_{i<k})$ is any inconsistency witness for $\phi(\bar x;\bar y)$.
	Then there is a sequence $(\bar b_i)_{i<\omega}$ such that
	$\bar b_i\equiv_CB$ for all $i<\omega$ and
	$\psi(\bar b_{i_0},\ldots,\bar b_{i_{k-1}})$ holds for all
	$i_0<\ldots<i_{k-1}<\omega$.
	Since $\psi(\bar y_{i<k})$ is an inconsistency witness for $\phi(\bar x;\bar y)$
	this implies that $\bigwedge_{i<k}\phi(\bar x;\bar y_i)$ is inconsistent
	for all $i_0<\ldots<i_{k-1}<\omega$. In other words:
	$\{\phi(\bar x;\bar b_i)\mid i<\omega\}$ is $k$-inconsistent.
	Therefore $\phi(\bar x;\bar b)$ divides over~$C$.
	
	(ii) $\phi(\bar x;\bar b)$ forks over~$C$
	iff $\phi(\bar x;\bar b)$ implies a finite disjunction
	$\bigvee_{i<n}\phi^i(\bar x;\bar y^i)$ of formulas
	$\phi^i(\bar x;\bar y^i)$ each of which divides over~$C$,
	or, which is equivalent by Exercise~\ref{ExcDividing},
	each of which $(\phi^i,\psi^i)$-divides over~$C$ for some
	inconsistency witness $\psi^i$ for~$\phi^i$.
\end{solution}

\begin{solution}{ExcOmindNotTransitiveNotNormal} (concerning Example~\ref{ExOmindNotTransitiveNotNormal})

	We prove that $A\nomind_CB$ iff the condition
	$\card{(A\cap B)\setminus C}\geq 2$ holds.
	$A\nomind_CB$ iff there are $a,a'\in A$ and $b,b'\in B$ s.\,t.{}
	$\models\phi(aa',bb')$ and a sequence $(b_ib'_i)_{i<\omega}$
	s.\,t.{} $b_ib_i'\equiv_Cbb'$ and $\models\psi(b_ib_i',b_jb_j')$
	for all $i<j<\omega$.
	This is the case iff there are $a,a'\in A\cap BC$ s.\,t.{} $a\not=a'$
	and a sequence $(b_ib_i')_{i<\omega}$
	s.\,t.{} $b_i\not=b_j$, $b_j\not=b_j'$ and
	$b_ib_i'\equiv_Caa'$ for all $i<j<\omega$.
	This is equivalent to existence of $a,a'\in A\cap BC$ s.\,t.{}
	$a\not=a'$, $a\not\in C$ and $a'\not\in C$.
	But that just means $\card{(A\cap BC)\setminus C}\geq 2$.
\end{solution}

\begin{solution}{ExcTreeProperty} (tree property)

	First note that given any formula $\phi(\bar x;\bar y)$ and $k<\omega$,
	the formula $\psi(\bar y_{<k})\equiv\neg\exists\bar x\bigwedge_{i<k}\phi(\bar x;\bar y_i)$
	is the most general $k$-inconsistency witness for~$\phi$ in the sense
	that whenever $\psi'(\bar y_{<k})$ is a $k$-inconsistency witness for~$\phi$
	and $(\bar b_i)_{i<\omega}$ is a sequence such that we have
	$\models\phi'(\bar b_{i_0},\dots,\bar b_{i_{k-1}})$
	for all $i_0<\dots<i_{k-1}<\omega$,
	then also $\models\phi(\bar b_{i_0},\dots,\bar b_{i_{k-1}})$
	for all $i_0<\dots<i_{k-1}<\omega$.
	Now note that $\D_{\phi,\psi}(\emptyset)=\infty$
	iff the unique element $\xi\in\{(\phi,\psi)\}^\omega$ is a dividing pattern,
	iff the type $\divseq_{\emptyset}^\xi$ is consistent.
	But the tree described in the exercise is just a partial realisation
	of $\divseq_{\emptyset}^\xi$.
\end{solution}

\begin{solution}{ExcOmindStar} (symmetry of $\omind$) 

	Symmetry of $\omind$ implies condition (5) of Theorem~\ref{TheoremOmindStar},
	so $\ind[$\Omega\ast$]$ is an independence relation.
	Concerning $\omind$, it now follows that $\omind$ satisfies existence.
	Using this we can show that $\ind[$\Omega$]$ satisfies extension
	(and hence $\omind=\ind[$\Omega\ast$]$, and therefore $\omind$ is also an
	independence relation):
	Suppose $A\omind_CB$ and $\hat B\supseteq B$.
	Let $\hat B'\equiv_{BC}\hat B$ be such that $\hat B'\omind_{BC}A$.
	Since $B\omind_CA$ we can apply transitivity, to get $\hat B'\omind_CA$
	and hence $A\omind_C\hat B'$.
\end{solution}

\begin{solution}{ExcMSymmetric} (M-symmetry)

	(i) Note that $\acl(AC)\cap B \supseteq \acl(C(A\cap B))$
	always holds for $C\subseteq B$.
	Now first suppose $\operatorname M(A,B)$ holds and $C$ is a set s.\,t.{}
	$A\cap B\subseteq C\subseteq B$.
	Then $\acl(AC)\cap\acl(BC)=\acl(A\acl C)\cap\acl B=\acl(\acl C(A\cap B))=\acl C$.
	Conversely, suppose $A\mind_{A\cap B}B$ and $C\subseteq B$.
	Let $C'=C(A\cap B)$. Then $\acl(AC)\cap B=\acl(AC')\cap(BC')=\acl C'=\acl(C(A\cap B))$.
	
	(ii) This follows from (i) since $A\mind_CB$ holds iff $\acl(AC)\mind_C\acl(BC)$,
	iff $\acl(AC)\cap\acl(BC)=\acl C$ and $\operatorname M(\acl(AC),\acl(BC))$.
\end{solution}


\begin{solution}{ExcCanonicalReduct} (independence in a reduct)
	
	(i) Let $\bar a$ enumerate~$A$.
	Let $(\bar a_i)_{i<\omega}$ be a $\ind$-Morley sequence for $\tp(\bar a/BC)$.
	Note that $\acl(C\bar a_{<k})\cap\acl(C\bar a_{\geq k})=\acl C = C$
	(for $k<\omega$) by anti-reflexivity of~$\ind$.
	Since
	$C\subseteq\acl'(C\bar a_{<k})\cap\acl'(C\bar a_{\geq k})\subseteq
	\acl(C\bar a_{<k})\cap\acl(C\bar a_{\geq k}) = C$
	we have $\acl'(C\bar a_{<k})\cap\acl'(C\bar a_{\geq k})=C$.
	Applying Lemma~\ref{LemAclMorley} we get $\bar a\pind_CB$.
	($\acl'$ denotes algebraic closure computed in~$T'$.)
	
	(ii) $\find$ is a \sir{} for $T\eq$ and for $T'{}\eq$.
	By elimination of hyperimaginaries and Corollary~\ref{CorSimpleForkingEHI}
	$\find$ is in fact a canonical independence	relation for $T'{}\eq$.
	Therefore we can apply~(i).
\end{solution}

\begin{solution}{ExcWcb} (weak canonical bases)
	
	(i) For $C\subseteq\acl B$, enumerated as $\bar c$, say, there is~$\bar c'$
	such that $\bar a\bar c\equiv_B\bar a'\bar c'$.
	Since $\bar c\equiv_C\bar c'$, $\bar c'$ also enumerates~$C$.
	Hence $\bar a\ind_CB\iff\bar a'\ind_CB$ by invariance.
	Both $\wcb(\bar a/\acl B)$ and $\wcb(\bar a'/\acl B)$ are defined as the
	smallest such~$C$, so they agree.
	
	(ii) The smallest algebraically closed set $C\subseteq\acl B$
	satisfying $\bar a'\ind_CB$ clearly satisfies $\bar a\ind_CB$,
	so $\wcb(\bar a/B)\subseteq C$.
	
	(iii) (3) implies (2) because $\wcb(\bar a/C)\subseteq\acl C$.\\
	(2) implies (1): By $\bar a\ind_{\wcb(\bar a/BC)}\acl(BC)$, $\wcb(\bar a/BC)\subseteq\acl C\subseteq\acl(BC)$
	and base monotonicity we have $\bar a\ind_{\acl C}\acl(BC)$,
	hence $\bar a\ind_CB$.\\
	(1) implies (3): Suppose $\bar a\ind_CB$.
	For algebraically closed $D\subseteq\acl(BC)$ we have
	the equivalence $\bar a\ind_DBC\iff\bar a\ind_DC$.
\end{solution}

\begin{solution}{ExcOneBased} (1-based theories) 

	(i)
	Recall from Exercise~\ref{ExcModular} that the lattice is modular iff
	$\aind$ is an independence relation.
	Now suppose $T$ is 1-based.
	Then $A\aind_CB$ implies $\acl(AC)\cap\acl(BC)=\acl C$,
	so $A\find_CB$. Since $\find$ is a \sir{}, $A\find_CB$ also implies $A\aind_CB$,
	so $\aind=\find$ is an independence relation.
	Conversely, suppose $\aind$ is an independence relation.
	Since $\aind$ is strict and $\find$ is the coarsest \sir{}
	because it is canonical, $A\aind_{\acl A\cap\acl B}$ implies $A\find_{\acl A\cap\acl B}$,
	so $T$ is 1-based.
	
	(ii)
	If $T$ is 1-based and perfectly trivial, then $\aind=\find$,
	so the lattice is distributive by Exercise~\ref{ExcModular}.
	Conversely, if the lattice is distributive, then it is also modular,
	hence $T$ is 1-based and $\aind=\find$. And $\aind$ is
	perfectly trivial, again by Exercise~\ref{ExcModular}.
\end{solution}


\begin{solution}{ExcOnshuus} (strong dividing and \th-forking)

	Suppose $\phi(x,b,d)$ divides quite strongly over~$Cc$.
	Consider the definable equivalence relation
	$\epsilon(yu,y'u')\equiv\forall x(\phi(x,y,u)\leftrightarrow\phi(x,y',u'))$
	and the formula $\psi(x,v)\equiv\exists y\exists u(v=(yu/\epsilon) \wedge \phi(x,y,u))$.
	Then $\phi(x,bd)$ is equivalent to $\psi(x,e)$ where $e=bd/\epsilon$.
	It is sufficient to show that $\psi(x,e)$ strongly divides over $Ccd$.
	
	By assumption the set
	$\big\{\phi(x,b',d) \;\big|\; {b'\models\tp(b/Cc)}\big\}$
	of formulas up to equivalence is $k$-inconsistent for some~$k$.
	Hence the same is true for the smaller set
	$\big\{\phi(x,b',d) \;\big|\; {b'\models\tp(b/Ccd)}\big\}$
	of formulas up to equivalence which we can also write as
	$\big\{\phi(x,b'd') \;\big|\; {b'd'\models\tp(bd/Ccd)}\big\}$.
	Hence the set
	$\big\{\psi(x,e') \;\big|\; {e'\models\tp(e/Ccd)}\big\}$
	of formulas up to equivalence is also $k$-inconsistent.
	But for this last set the qualification `up to equivalence'
	is unnecessary because every $e'$ is a canonical parameter.
	Therefore $\psi(x,e)$ strongly divides over $Ccd$.
\end{solution}

\par}


\section{Open questions}\label{SectionQuestions}

\begin{question}
	If $\find$ satisfies existence, does it follow that $\find=\dind$?
\end{question}

\begin{question}
	If $\find=\dind$, does it follow that $\thind=\mind$?
\end{question}

\begin{question}\label{QEveryIndSFC}
	Does every independence relation have strong finite character?
\end{question}

\begin{question}
	Is every independence relation of the form $\omsind$\;?
\end{question}

I would have liked to provide counterexamples for these four questions,
but I did not find any.

\begin{question}\label{QTwoSirs}
	Is there a theory with more than one \emph{strict} independence relation on $T\eq$?
\end{question}

It would be a good thing if not:
Given any simple theory, Shelah-forking independence
$\find$ is a strict independence relation on $T\eq$.
It follows from Theorem~\ref{ThmThind} that thorn-forking independence is also a strict
independence relation on $T\eq$.
If $\find$ and $\thind$ agree, then we are very close
to elimination of hyperimaginaries.

\begin{question}
	Is every rosy theory the reduct of a theory with
	elimination of hyperimaginaries and such that $\mind$\,
	is symmetric?
\end{question}

This is perhaps not important, but it is the second half of a question I have
been asking myself privately for years. The first half
(`Does symmetry of $\mind$ imply that $\mind$\, is an independence relation?')
was answered by Theorem~\ref{ThmMsymmetricSir}.

\begin{question}
	Does `hyperimaginary thorn-forking' make sense?
	In other words:
	How much of Section~\ref{SectionThornForking} can be carried out for $T\heq$ with
	bounded closure $\bdd$ instead of algebraic closure $\acl$?
\end{question}

In a simple theory, Shelah-forking can be extended to and has canonical bases in $T\heq$.
Therefore Shelah-forking agrees with hyperimaginary thorn-forking in this case.
This suggests that perhaps hyperimaginary thorn-forking is a more natural
notion than thorn-forking. On the other hand it is not the right notion of
independence in o-minimal theories.
If hyperimaginary thorn-forking does not make sense in general,
then this might be considered a heuristic argument for elimination of hyperimaginaries
(at least in a weak sense) in all simple theories.
I think Clifton Ealy may have some partial results.

\begin{question}
	For $T\heq$ define $A\ind[b]_CB\iff \bdd(AC)\cap\bdd(BC)=\bdd C$.
	Does $\ind[b]$ always satisfy existence (at least over sufficiently saturated models)?
	Suppose $T$ is such that $\ind[b]$ satisfies existence and base
	monotonicity. Does it follow that $\ind[b]$ satisfies finite character?
\end{question}

The previous question may be too hard. This is a simpler test question.

\subsubsection*{An o-minimal theory without weak canonical bases}

Up to a very late draft of this thesis I asked if $\thind$ is canonical for $T\eq$ for every
o-minimal theory. All o-minimal theories are rosy, and I did not even know a rosy
theory such that $\thind$ is not canonical for $T\eq$.
I am grateful to Anand Pillay for directing me to~\cite{Loveys + Peterzil: Linear o-minimal structures},
which contains a counterexample (Example~4.5) that I present here in a slightly
untangled form.\index{example!no weak canonical bases}

We will construct two o-minimal theories $T$ and $T'$ with elimination of imaginaries
such that $T'$ is a reduct of $T$, $\thind$ is a canonical independence relation for~$T$,
but $\thind$ for $T'$ does not satisfy the intersection property. $T$ will be interpretable
in the theory of $(\R; +, <)$.

Let $\R_<$ and $\R_>$ be two copies of the set $\R$ of real numbers.
Consider the following structure~$\mathcal R$:
The domain of $\mathcal R$ consists of the disjoint union of $\R_<$,
$\R_>$ and a new point $*$.
There are constants for the two points $-1$ and $1$ in $\R_<$ as well as for
the points $-1$ and $1$ in~$\R_>$ and the point~$*$.
The structure has a linear order $<$ extending the usual order on
$\R_<$ and $\R_>$ and such that $\R_<<*<\R_>$.
There is also a partial binary function $+$ that is defined on pairs from
$\R_<$ and on pairs from $\R_>$ and is interpreted as addition in $\R_<$ or $\R_>$, respectively.
Finally, there is a partial function $f$ defined on $\R_<$:
the obvious order-preserving isomorphism $\R_<\rightarrow\R_>$.
(The only purpose of the point $*$ is to make the definable sets
$\R_<$ and $\R_>$ definable open intervals. Otherwise $\mathcal R$ would
not be o-minimal.)

Let $T$ be the theory of~$\mathcal R$.
Since $T$ is essentially just the theory of $\R$ as an ordered group,
it is straightforward to check that
$T$ is o-minimal and that $T$ has elimination of imaginaries and a modular lattice of algebraically
closed sets. Hence $\thind$ is a canonical independence relation for $T\eq$.

Now consider $T'$, the theory of the following variant $\mathcal R'$ of~$\mathcal R$:
Instead of $f$ we have a partial function $f'$
which is the restriction of $f$ to the interval $(-1,1)$ of~$\R_<$.
Everything else is as in the definition of $\mathcal R$.
Note that $T'$ is essentially a reduct of~$T$.

Let $\overline{\mathcal R}$ be a big elementary extension of $\mathcal R$
and $\overline{\mathcal R}'$ its `reduct' to the signature of~$T'$.
Let $\overline{\mathcal R}_<=\{a\in\overline{\mathcal R}\mid a<*\}\supseteq\R_<$ and
$\overline{\mathcal R}_>=\{a\in\overline{\mathcal R}\mid a>*\}\supseteq\R_>$.
The (algebraic) closure of a set $A$ in $\overline{\mathcal R}'$
turns out to be the smallest set $C\supseteq A$ containing
the five constants and
such that $C\cap\overline{\mathcal R}_<$ and $C\cap\overline{\mathcal R}_>$ are divisible subgroups of
$\overline{\mathcal R}_<$ or $\overline{\mathcal R}_>$, respectively, and which is closed under $f'$ and $f'{}^{-1}$.
Since every closed subset of a model of $T'$ is an elementary submodel,
$T'$ has elimination of imaginaries by~\cite[Proposition~3.2]{Anand Pillay: Some remarks on definable equivalence relations in O-minimal structures}.

Yet $\thind$ for $T'$ does not have the intersection property:
Choose an element $a$ such that $\R_<<a<*$, and let $b=f(a)$.
Note that $\dim(ab)=1$ when evaluated in~$T$, but $\dim(ab)=2$ when evaluated in~$T'$.
Choose $\epsilon_1,\epsilon_2$ in the interval $(-1,1)$ of $\overline{\mathcal R}_<$
such that $\dim(\epsilon_1\epsilon_2)=2$, and set $a_i=a+\epsilon_i$, $b_i=f(a_i)$.

From now on all calculations will be in~$T'$.
Now $b=b_1-f'(a_1-a)\in\cl(a_1b_1a)$ where $\cl$ is (algebraic) closure in $T'$,
and it is only a matter of linear algebra to check that we have
$b\not\in\cl(a_1b_1a_2b_2)$.

Hence $\dim(ab/a_1b_1a_2b_2)=\dim(ab/a_1b_1)=\dim(ab/a_2b_2)=1$,
and therefore $ab\thind_{a_1b_1}a_1b_1a_2b_2$
and also $ab\thind_{a_2b_2}a_1b_1a_2b_2$.
Yet $\cl(a_1b_1)\cap\cl(a_2b_2)=\cl\emptyset$,
and $ab\nthind_\emptyset a_1b_1a_2b_2$ because $b\in\cl(a_1b_1a)$ and $b\not\in\cl(a)$.
Therefore $ab\nthind_{\cl(a_1b_1)\cap\cl(a_2,b_2)}a_1b_1a_2b_2$ witnesses that the
intersection property fails for $\tp(ab/\cl(a_1b_1a_2b_2))$.

Thus $T'$ is in fact an o-minimal theory for which $\thind$ is not a canonical independence relation,
and $T'$ is the reduct of another o-minimal theory for which $\thind$ is canonical.
Therefore $\thind$ for $T\eq$ need not be canonical for rosy theories, in fact not even for
o-minimal theories; and having a canonical independence relation is not preserved under
reducts.

\clearpage\phantomsection
\addcontentsline{toc}{chapter}{Bibliography, Index}

	\clearpage\phantomsection	
	\printindex
\end{document}